\documentclass[11pt]{article}

\usepackage{amsmath,amssymb,amsfonts,enumitem,amsthm,accents,pstricks,pst-node,pstricks-add,mathrsfs,xy}
\xyoption{all}

\usepackage{color}
\usepackage{makeidx}
\usepackage{graphicx}
\usepackage{setspace}
\usepackage{extarrows}

\numberwithin{equation}{section}

\usepackage{tocloft}

\newtheorem{theorem}{Theorem}[section]
\newtheorem{lemma}[theorem]{Lemma}
\newtheorem{conjlemma}[theorem]{Conjectural Lemma}
\newtheorem{corollary}[theorem]{Corollary}
\newtheorem{proposition}[theorem]{Proposition}
\newtheorem{conjecture}[theorem]{Conjecture}

\theoremstyle{definition}
\newtheorem{definition}[theorem]{Definition}
\newtheorem{remark}[theorem]{Remark}
\newtheorem{example}[theorem]{Example}

\def\ov#1{\overline{#1}}
\def\tn#1{\textnormal{#1}}
\def\mf#1{\mathfrak{#1}}
\def\wt#1{\widetilde{#1}}
\def\wh#1{\widehat{#1}}
\def\scz{\scriptsize}
\def\vr{\varrho}
\def\ll{\left\langle}
\def\rr{\right\rangle}
\def\mc{\mathcal}
\def\lra{\longrightarrow}
\def\dbar{\bar\partial}

\def\ve{\varepsilon}

\newcommand\uvec[1]{\underaccent{\vec}{#1}}
\newcommand{\ucev}[1]{\reflectbox{\ensuremath{\uvec{\reflectbox{\ensuremath{#1}}}}}}

\newcommand{\abs}[1]{\left\vert #1 \right\vert}
\newcommand{\lrp}[1]{\left( #1 \right)}
\newcommand{\lrc}[1]{\left\{ #1 \right\}}

\def\scz{\scriptsize}

\def\bEq#1{\begin{equation}\label{#1}}
\def\eEq{\end{equation}}

\def\bsEq{\begin{equation*}}
\def\esEq{\end{equation*}}

\def\bDf#1{\begin{definition}\label{#1}}
\def\eDf{\end{definition}}

\def\bTh#1{\begin{theorem}\label{#1}}
\def\eTh{\end{theorem}}

\def\bCn#1{\begin{conjecture}\label{#1}}
\def\eCn{\end{conjecture}}

\def\bLm#1{\begin{lemma}\label{#1}}
\def\eLm{\end{lemma}}

\def\bCLm#1{\begin{conjlemma}\label{#1}}
\def\eCLm{\end{conjlemma}}

\def\bRm#1{\begin{remark}\label{#1}}
\def\eRm{\end{remark}}

\def\bEx#1{\begin{example}\label{#1}}
\def\eEx{\end{example}}

\def\bPr#1{\begin{proposition}\label{#1}}
\def\ePr{\end{proposition}}

\def\bCr#1{\begin{corollary}\label{#1}}
\def\eCr{\end{corollary}}

\def\bFg#1{\begin{figure}\label{#1}}
\def\eFg{\end{figure}}

\def\bPf{\begin{proof}}
\def\ePf{\end{proof}}

\def\bIt{\begin{itemize}[leftmargin=*]}
\def\eIt{\end{itemize}}

\def\bEn{\begin{enumerate}[label=$(\arabic*)$,leftmargin=*]}
\def\eEn{\end{enumerate}}

\def\bEnalph{\begin{enumerate}[label=$(\alph*)$,leftmargin=*]}
\def\eEnalph{\end{enumerate}}

\def\AK{\tn{AK}}
\def\Symp{\tn{Symp}}

\def\coker{\tn{coker}}

\def\nd{\tn{d}}

\def\Def{\tn{Def}}

\def\id{\tn{id}}
\def\ker{\tn{ker}}

\def\Obs{\tn{Obs}}
\def\ord{\tn{ord}}

\def\pt{\tn{pt}}

\def\ob{\tn{ob}}

\def\cJ{\mc{J}}
\def\cM{\mc{M}}
\def\cO{\mc{O}}

\def\cP{\mc{P}}

\def\cN{\mc{N}}
\def\cR{\mc{R}}

\def\cZ{\mc{Z}}

\def\cL{\mc{L}}

\def\E{\mathbb E}
\def\R{\mathbb R}
\def\C{\mathbb C}
\def\Z{\mathbb Z}

\def\P{\mathbb P}

\def\T{\mathbb T}
\def\N{\mathbb N}
\def\D{\mathbb D}
\def\L{\mathbb L}
\def\K{\mathbb K}
\def\CK{\mathbb{C}\mathbb{K}}
\def\V{\mathbb V}
\def\mfi{\mf{i}}
\def\mfj{\mf{j}}

\def\aut{\mf{Aut}}

\def\mfs{\mf{s}}

\def\la{\lambda}
\def\La{\Lambda}
\def\ep{\epsilon}
\def\De{\Delta}
\def\de{\delta}
\def\om{\omega}
\def\Om{\Omega}
\def\si{\sigma}
\def\Si{\Sigma}
\def\al{\alpha}

\def\Ga{\Gamma}
\def\ze{\zeta}
\def\na{\nabla}
\def\eset{\emptyset}

\makeindex
\begin{document}
\title{Pseudoholomorphic Curves Relative to a Normal Crossings\\ Symplectic Divisor: Compactification}
\author{Mohammad Farajzadeh-Tehrani}
\date{September 17 2017, Last updated \today}
\maketitle

\begin{abstract}
Inspired by the log Gromov-Witten (or GW) theory of Gross-Siebert/Abramovich-Chen, we introduce a geometric notion of \textit{log $J$-holomorphic curve} relative to a simple normal crossings  symplectic divisor defined in \cite{FMZ1}. 
Every such moduli space is characterized by a second homology class, genus, and contact data.
For certain almost complex structures, we show that the moduli space of stable log $J$-holomorphic curves of any fixed type is compact and metrizable with respect to an enhancement of the Gromov topology. 
In the case of smooth symplectic divisors, our compactification is often smaller than the relative compactification and there is a projection map from the latter onto the former. 
The latter is constructed via expanded degenerations of the target. Our construction does not need any modification of (or any extra structure on) the target.
Unlike the classical moduli spaces of stable maps, these log moduli spaces are often virtually singular. We describe an explicit toric model for the normal cone (i.e. the space of gluing parameters) to each stratum in terms of the defining combinatorial data of that stratum. 
In \cite{FT2}, we introduce a natural set up for studying the deformation theory of log (and relative) curves and obtain a logarithmic analog of the space of Ruan-Tian perturbations for these moduli spaces. In a forthcoming paper, we will prove a gluing theorem for smoothing log curves in the normal direction to each stratum. With some modifications to the theory of Kuranishi spaces, the latter will allow us to construct a virtual fundamental class for every such log moduli space and define relative GW invariants without any restriction.
\end{abstract}
\tableofcontents

%--------------------------------------------------------------
\section{Introduction}\label{intro_s}

Studying pairs of a smooth variety $X$ and a normal crossings (or \textbf{NC}) divisor\footnote{curves and divisors are sub-varieties of dimension 1 and codimension 1 over the ground field, respectively.} $D\!\subset\!X$ has a rich history in complex algebraic geometry.  For example, studying such pairs is central to the Minimal Model Program and to the construction of moduli spaces in algebraic geometry. By a celebrated theorem of Hironaka (1964), given a singular variety $Y$, there is a smooth ``blowup" $X$ of $Y$ such that the pre-image of the singular locus of $Y$ is an NC divisor $D\!\subset\!X$. Therefore, study of such pairs is also important toward the study of singularities. Curves are (Poincar\'e) dual objects to divisors. Moduli spaces of curves in $X$ that intersect $D$ in some particular ways are fundamental tools for understanding the geometry of $(X,D)$.   \\

\noindent
In the last $40$ years, analogue of these notions have been defined in the symplectic category and have led to significant advances in our understanding of symplectic manifolds.  In 1980s, Gromov combined the rigidity of algebraic geometry with the flexibility of the smooth category and initiated the use of $J$-holomorphic curves as a generalization of holomorphic curves in symplectic geometry. The use of $J$-holomorphic curve techniques has led to numerous connections with algebraic geometry, string theory, and to the appearance of symplectic divisors (as  the dual objects) in various contexts. The latter includes relations with complex line bundles \cite{D}, relative Gromov-Witten (or \textbf{GW}) theory \cite{IP1,LR, BP7}, degeneration formulas for GW invariants \cite{IP2,LR, BP5, FZ}, topological study of singularities \cite{Mc}, symplectic cohomology and mirror symmetry of complements $X \backslash D$ \cite{A,GP}, and classification of symplectic log Calabi-Yau 4-manifolds \cite{LM}. A smooth symplectic divisor is simply a symplectic submanifold of real codimension two. Topological notions of NC symplectic divisors and varieties were recently introduced by M.~McLean, A.~Zinger, and the author in \cite{FMZ1,FMZ2,FMZ3}.  \\

\noindent
While most applications of $J$-holomorphic curves in symplectic topology have so far concerned smooth symplectic manifolds or pairs $(X,D)$ of a smooth manifold and a smooth symplectic divisor, recent developments in symplectic topology and the existing rich structures in algebraic geometry (some of which are listed above) suggest the need for constructing and studying moduli spaces of $J$-holomorphic curves relative to an arbitrary NC symplectic divisor from the analytical perspective. In this paper we introduce an explicit and efficient compactification of moduli spaces of $J$-holomorphic curves relative to an arbitrary simple normal crossings (\textbf{SNC}) symplectic divisor. In upcoming papers \cite{FT2,FT3},  we will setup the analytic framework needed for constructing a (virtual) fundamental class and define relative GW invariants. In particular, in \cite{FT2}, we will define a notion of semi-positive pair that allows a direct construction of relative GW invariants via perturbed $J$-holomorphic maps as in \cite{RT}. In \cite{FT1}, based on these log moduli spaces, we outline an explicit degeneration formula that relates the GW invariants of smooth fibers to the GW invariants of central fiber, in a semistable degeneration with an SNC central fiber. It is worth mentioning that even in the case of smooth divisors, our compactification is different and smaller than the well-known relative compactification in  \cite{JLi1, IP1,LR}. \\

\noindent
We begin by setting up the most commonly used notation and recalling some of the known facts about the classical and relative moduli spaces of closed $J$-holomorphic curves. Therefore, experts may skip to Section~\ref{RelSC_ss} where the main question is explained.

\subsection{Classical stable maps and GW invariants}
For a smooth manifold $X$, $g,k\!\in\!\N$, $A\!\in\!H_2(X,\Z)$,  and an almost complex structure\footnote{i.e. $J$ is a real-linear endomorphism of $TX$ lifting the identity map satisfying $J^2\!=\!-\tn{id}_{TX}$.} $J$ on $X$, a (nodal) $k$-marked genus $g$ degree $A$  $J$-holomorphic map into $X$  is a tuple $\big(u,\Si,\mfj,z^1,\ldots,z^k\big)$, where
\bIt
\item $(\Si,\mfj)$ is a connected nodal Riemann surface of arithmetic genus $g$ with $k$ distinct ordered  marked points $z^1,\ldots,z^k$ away from the nodes, 
\item $u\colon (\Si,\mfj)\lra(X,J)$ is a continuous and component-wise smooth map satisfying the Cauchy-Riemann equation 
\bEq{CR_e}
\dbar u \!=\frac{1}{2}(\nd u \!+\!J\nd u\circ \mfj )\!=\!0
\eEq
on each smooth component, and
\item the map $u$ represents the homology class $A$. 
\eIt
Two such tuples 
$$
(u,\Si,\mfj,z^1,\ldots,z^k)\quad\tn{and}\quad (u',\Si',\mfj',w^1,\ldots,w^k)
$$ 
are \textbf{equivalent} if there exists a biholomorphic isomorphism $h\colon\!(\Si,\mfj)\!\lra\!(\Si',\mfj')$ such that $h(z^a)\!=\!w^a$, for all $a\!=\!1,\ldots,k$, and $u\!=\!u'\circ h$. Such a tuple is called stable if the group of self-automorphisms is finite. Let $\ov\cM_{g,k}(X,A,J)$ (or simply $\ov\cM_{g,k}(X,A)$ when $J$ is fixed in the discussion) denote the space of equivalence classes of stable $k$-marked genus $g$ degree $A$ $J$-holomorphic maps into $X$. Such an equivalence class is called a marked $J$-holomorphic \textit{curve}.\\

\noindent
By a celebrated theorem\footnote{and its subsequent refinements; see before Theorem~\ref{Gromov_Th}.} of Gromov \cite[Thm 1.5.B]{Gr}, for every smooth closed (i.e. compact and without boundary) symplectic manifold $(X,\om)$, $g,k,A$ as above, and an almost complex structure $J$ compatible\footnote{i.e. $\om(\cdot,J\cdot)$ is a metric.} with $\om$ (or taming $\om$), the moduli space $\ov\cM_{g,k}(X,A,J)$ has a natural sequential convergence topology, called the \textbf{Gromov topology}, which is compact, Hausdorff, and furthermore metrizable. The symplectic structure only gives an energy bound which is needed for establishing the compactness, and the precise choice of that, up to deformation, is not important. If $\ov\cM_{g,k}(X,A)$ has an oriented orbifold structure of the expected real dimension
\bEq{exp-dim_e}
2\big(c_1^{TX}(A)+(n-3)(1-g)+k\big),
\eEq
GW invariants are obtained by the integration of appropriate cohomology classes against its fundamental class. 
These numbers are independent of $J$ and only depend on the deformation equivalence class of $\om$.
These allow the formulation of symplectic analogues of enumerative questions from algebraic geometry, as well-defined invariants of symplectic manifolds.
However,  in general, such moduli spaces can be highly singular. 
This issue is known as the \textit{transversality problem}. 
Fortunately, it has been shown  (e.g. see\footnote{it is beyond the scope of this paper to list all the related literature.} \cite{LT,FO, MS94, HWZ, MW,JP}) that $\ov\cM_{g,k}(X,A)$ still carries a rational homology class, called \textbf{virtual fundamental class} (or \textbf{VFC}); integration of cohomology classes against VFC gives rise to GW-invariants.

\subsection{Relative stable maps}
Given a symplectic manifold $(X,\om)$ and a closed submanifold $D\!\subset\!X$, we say $D\!\subset\!X$ is a symplectic submanifold if $\om|_{D}$ is a symplectic structure.  A (smooth) \textbf{symplectic divisor} is a symplectic submanifold of real codimension $2$. For such $D$ (or a smooth divisor in complex algebraic geometry), \textbf{relative} GW theory  (virtually) counts $J$-holomorphic curves in $X$ with a fixed contact order $\mfs\!\equiv\! (s_1,\cdots,s_k)\!\in\!\N^k$ with $D$. In this theory, we require $J$ to be also \textbf{compatible} with $D$ in the following sense. First, we require $D$ to be $J$-holomorphic, i.e. $J(TD)\!=\!TD$. This implies, for example, that every $J$-holomorphic map to $X$ from a smooth domain is either mapped into $D$ or intersects $D$ positively in a finite set of points. Furthermore, we need to at least require $J$ to be integrable to the first order in the normal direction to $D$  in the sense that
\bEq{intInnormal_e0}
N_J(v_1,v_2)\in T_xD \qquad \forall~x\!\in\!D,~v_1,v_2\!\in\!T_xX,\quad\tn{where }
\eEq
$$
N_J \!\in\! \Gamma(X,\Om^2_X\otimes TX),\quad N_J(u,v)\!\equiv\![u,v]\!+\!J[u,Jv]\!+\!J[Ju,v]\!-\![Ju,Jv]\quad\forall u,v\!\in\! TX,
$$
is the Nijenhueis tensor of $J$. This ensures that certain operators are complex linear, see (\ref{DuNXV_e}), and certain sequence of almost complex structures on the normal bundle $\cN_XD$ converges to a standard one, see Lemma~\ref{JLemma_e}. The space $\cJ(X,D,\om)$ of $\om$-tame and $D$-compatible almost complex structures $J$ on $X$ is again non-empty and contractible. For every $J\!\in\!\cJ(X,D,\om)$ and $\mfs\!\equiv\! (s_1,\cdots,s_k)\!\in\!\N^k$, with 
\bEq{DegreeDecomp_e}
\sum_{a=1}^k s_a = A\cdot D,
\eEq
let $\cM_{g,\mfs}(X,D,A)\!\subset\!\cM_{g,k}(X,A)$ (in the stable range) denote the sub-space of $k$-marked  degree $A$ genus $g$ $J$-holomorphic curves $[u,\Si,\mfj,z^1,\ldots,z^k]$ such that $\Si$ is smooth and $u$ has a tangency of order $s_a$ at $z^a$ with $D$. In particular, by (\ref{DegreeDecomp_e}),
$$
u^{-1}(D)\!\subset\!\{z^1,\ldots,z^k\}.
$$
The subset of marked points $z^a$ with $s_a\!=\!0$ corresponds to the classical marked points of the classical GW theory with image away from $D$. The relative compactification $\ov\cM_{g,\mfs}^{\tn{rel}}(X,D,A)$ of $\cM_{g,\mfs}(X,D,A)$ constructed in \cite{JLi1}, in the algebraic case, and in \cite{IP1,LR}, in the symplectic case, includes stable nodal curves with components mapped into $X$ or an \textbf{expanded degeneration}\footnote{a normal crossings variety made of $X$ and finite copies of the $\P^1$-bundle $\P_XD=\P(\cN_XD\oplus \cO_D)$ over $D$.}  of that so that the contact order $\mfs$ still makes sense; we will review this construction in Section~\ref{RelComp_ss}.

\subsection{$J$-holomorphic maps relative SNC divisors}\label{RelSC_ss}
In \cite{FMZ1,FMZ2}, with M.~Mclean and A.~Zinger, we defined topological notions of symplectic normal crossings divisor and variety and 
showed that they are equivalent, in a suitable sense, to the desired rigid notions. For $N\!\in\!\N$, let 
$$
[N]=\{1,\ldots,N\};
$$ 
in particular $[0]\!=\!\eset$.
A \textbf{simple normal crossings (or SNC) symplectic divisor} $D\!=\bigcup_{i\in [N]} D_i$ in $(X,\om)$ is a transverse union 
of smooth symplectic divisors $\{D_i\}_{i\in [N]}$  in $X$ such that all the strata 
$$
D_I\equiv \bigcap_{i\in I} D_i\qquad \forall I\!\subset\![N]
$$
are symplectic, and the symplectic orientation of $D_I$ coincides with its ``intersection" orientation for all $I\!\subset\![N]$; see \cite[Dfn~2.1]{FMZ1}. For  
$$
J\!\in\!\cJ(X,D,\om)=\bigcap_{i\in [N]}\cJ(X,D_i,\om),
$$  
we similarly  define $\cM_{g,\mfs}(X,D,A)$ (in the stable range) to be the space of equivalence classes of degree $A$ $J$-holomorphic maps from a $k$-marked genus $g$ connected smooth domain $\Si$ into $X$ of contact order $\mfs$ with $D$, for which   
$$
\mfs\!\equiv\!\big(s_a\equiv (s_{ai})_{i\in [N]}\big)_{a\in [k]} \!\in\!(\N^N)^k,
$$
each vector $s_a$ records the intersection numbers of the $a$-th marked point $z^a$ with the divisors $\{D_i\}_{i\in [N]}$, and 
\bEq{Equiv-mfs-conditions_e}
u^{-1}(D)\!\subset \{z^1,\ldots,z^k\},\quad\tn{or equivalently}\quad A\cdot D_i\!=\!\sum_{a=1}^k s_{ai}\quad \forall~i\!\in\![N].
\eEq
\noindent
Because of the tangency conditions, it follows from (\ref{exp-dim_e}) that the expected real dimension of $\cM_{g,\mfs}(X,D,A)$ is equal to
\bEq{dlog_e}
2\Big(c_1^{TX}(A)+(n-3)(1-g)+ k-A\cdot D\Big)=2\Big(c_1^{TX(-\log D)}(A)+(n-3)(1-g)+ k\Big),
\eEq
where $TX(-\log D)$ is the log tangent bundle associated to the deformation equivalence class of $(X,D,\om)$, defined in \cite[(8)]{FMZ2}. In the holomorphic case, the \textbf{log tangent sheaf} is the sheaf of holomorphic tangent vector fields in $TX$ whose restriction to each $D_i$ is tangent to $D_i$. The definition in the symplectic case is similar but depends\footnote{the deformation equivalence class of complex vector bundle $TX(-\log D)$ is independent of the auxiliary data.} on some auxiliary data. The similarity between the left-hand side of (\ref{dlog_e}) and (\ref{exp-dim_e}) shows the importance of considering the log tangent bundle in the study of relative moduli spaces. \textbf{The main goal is:}\vskip.1in

\noindent
($\star$)~\textit{to construct a natural geometric compactification $\ov\cM_{g,\mfs}(X,D,A)$ of $\cM_{g,\mfs}(X,D,A)$ so that the definition of the contact vector $\mfs$ naturally extends to every element of $\ov\cM_{g,\mfs}(X,D,A)$, and $\ov\cM_{g,\mfs}(X,D,A)$ is (virtually) smooth enough to admit a natural class of cobordant  Kuranishi structures of the expected real dimension (\ref{dlog_e}).}

\noindent
We refer to \cite{FF,MW} for the technical terms in ($\star$).
If $D$ is smooth, the well-known relative compactification $\ov\cM_{g,\mfs}^{\tn{rel}}(X,D,A)$ has (or it is expected\footnote{cf. \cite{FZ} for an overview of the analytical approaches of \cite{IP1,LR}.} to have) these nice properties.\\

\noindent
In the algebraic category, every (algebraic) NC variety $D\!\subset\!X$ defines a natural ``fine saturated log structure'' on $X$; see \cite{ACGHOSS} for a review of log geometry and log moduli spaces associated to NC pairs $(X,D)$. Then the log GW theory of \cite{AC,GS} constructs a good compactification with a perfect obstruction theory for every fine saturated log variety $X$. 
Unlike in \cite{JLi1}, the algebraic log compactification does not require any expanded degeneration of the target. Instead, it uses the extra log structure on $X$ (and various log structures on the domains) to keep track of the contact data for the curves that have image inside the support of the log structure (i.e. $D$).\\

\noindent 
Since the classical GW invariants are invariants of the deformation equivalence class of the underlying symplectic structure,  it is interesting and important to generalize the results of \cite{AC,GS} to (or find an analogue of that for) the symplectic category, i.e. to construct log GW invariants as invariants of the symplectic deformation equivalence class of $(X,D)$. With such a construction, the flexibility of symplectic topology can be used in certain situations to define log GW invariants as an actual count of $J$-holomorphic curves with tangency conditions, at the expense of deforming $J$ or the Cauchy-Riemann equation (to avoid working with VFC); see \cite{RT,FT2}. Moreover, in the case of moduli spaces of holomorphic curves with boundary on Lagrangian submanifolds, it is sometimes easier to work with an analytical construction of moduli spaces of $J$-holomorphic maps.
\\

\noindent
On the analytical side, in \cite{BP2,BP7,BP8} and several other related papers, Brett Parker uses his enriched almost K\"ahler category of ``exploded manifolds", defined in \cite{BP1}, to construct such a compactification relative to an almost K\"ahler NC divisor and address ($\star$). His approach can be considered as a direct translation/generalization of the algebraic log GW theory involving some non-Hausdorff spaces, analytical sheaves, and a richer cohomology theory \cite{BP4}. His approach has close ties to Tropical Geometry. In \cite{I},  Eleny Ionel approaches ($\star$), by considering  expanded degenerations similar to \cite{IP1}.
Nevertheless, the main motivation behind the log GW theory of Gross-Siebert-Abramovich-Chen, exploded theory of Parker, and the current paper is that considering spaces and maps enriched with certain log structures is a better idea for  addressing ($\star$) in the general case. In particular, all these logarithmic approaches lead to similar ``Degeneration Formulas" (\cite{ACGS} call it  an ``invariance property") relating the moduli spaces in smooth fibers and the SNC central fiber of an arbitrary semistable degeneration; see \cite{ ACGS, BP5, FT2}.

\subsection{Log compactification and the main result}
In this paper, for an arbitrary SNC symplectic divisor $D\!\subset\!(X,\om)$ and certain $J\!\in\!\cJ(X,D,\om)$, we construct a ``minimal geometric compactification" 
\bEq{MainRelSpace_e}
\ov\cM_{g,\mfs}^{\tn{log}}(X,D,A)
\eEq
that does not require any modification of the target (or the nodal domains).  For its connection to the algebraic log maps, and the appearance of various log structures\footnote{such as the use of log tangent bundle in the deformation theory of log $J$-holomorphic curves.} throughout the construction, we call our maps/curves: \textbf{log $J$-holomorphic maps/curves}. \\

\noindent
For $J\!\in\!\cJ(X,D,\om)$, a (nodal) log $J$-holomorphic map into $(X,D=\bigcup_{i\in [N]} D_i)$ of contact type 
$$
\mfs\!\equiv\!\big(s_a\equiv (s_{ai})_{i\in [N]}\big)_{a=1}^k \!\in\!(\Z^N)^k,
$$
with the marked nodal domain $(\Si,\mfj,\vec{z})\!=\!\bigcup_{v\in \V}(\Si_v,\mfj_v,\vec{z}_v)$, is collection of tuples
$$
u_{\log}\equiv \big((u_v\colon\!\Si_v\!\lra\! D_{I_v}, \vec{z}_v), ([\ze_{v,i}])_{i\in I_v}\big)_{v\in \V}
$$
over smooth components of $\Si$ such that 
\bIt
\item $u\!\equiv\!(u_v)_{v\in \V}\colon\!(\Si,\mfj,\vec{z})\!\lra\! (X,J)$ is a $k$-marked $J$-holomorphic nodal map in the classical sense,
\item for each $v\!\in\!\V$, $I_v\!\subset\![N]$ is the maximal subset such that $\tn{Im}(u_v)\!\subset\!D_{I_v}$,
\item for each $v\!\in\!\V$ and any $i\!\in\!I_v$, $[\ze_{v,i}]$ is the $\C^*$-equivalence class\footnote{$\C^*$ acts by multiplication on the set of meromorphic sections.} of a non-trivial meromorphic section $\ze_{v,i}$ of the holomorphic\footnote{since $\dim_\C \Si_v\!=\!1$, the pull-back line bundle $u_v^*\cN_{X}D_i$ is holomorphic.} line bundle $u^* \cN_{X}D_i$,
\item the contact order vectors in $\Z^N$, defined in (\ref{Ordx_e1}) and  (\ref{Ordx_e2}), are the opposite of each other at the nodal points, 
\item every point in $\Si$ with a non-trivial contact vector is either a marked point or a nodal point, and the contact order vector at $z^a$ is the pre-determined vector $s_a\!\in\!\Z^N$,
\item there exists a vector-valued function $s\colon\!\V\!\lra\!\R^N$ such that $s_v\!=\!s(v)\!\in\! \R_{+}^{I_v}\!\times \{0\}^{[N]-I_v}$ for all $v\!\in\!\V$, and $s_v\!-\!s_{v'}$ is a positive multiple of the contact order vector of any nodal point on $\Si_v$ connected to $\Si_{v'}$, for all $v,v'\!\in\!\V$,
\item\label{GObs_it} and, certain group (a complex torus) element associated to $u_{\log}$, defined in (\ref{PLtoG_e}), is equal to $1$;
\eIt
see Definition~\ref{LogMap_dfn} for more details. 
Two marked log maps are \textbf{equivalent} if one is a ``reparametrization" of the other.
A marked log map is stable if it has a finite ``automorphism group". For $g,k\!\in\!\N$, $A\!\in\!H_2(X,\Z)$, and $\mfs\!\in\!(\Z^N)^{k}$, \textit{we denote the space of equivalence classes of stable  $k$-marked degree $A$ genus $g$  log maps of contact type $\mfs$ by} 
$$
\ov\cM_{g,\mfs}^{\log}(X,D,A). 
$$
Such an equivalence class is called a log \textit{curve}.
There is natural forgetful map
$$
\ov\cM_{g,\mfs}^{\log}(X,D,A)\lra \ov\cM_{g,k}(X,A),\qquad \big((u_v\colon\!\Si_v\!\to\! D_{I_v},\vec{z}_v), ([\ze_{v,i}])_{i\in I_v}\big)_{v\in \V}\lra \big(u_v\colon\!\Si_v\!\to\! X,\vec{z}_v\big)_{v\in \V}\,.
$$
Given $\mfs\!\in\!(\Z^N)^k$, it turns out  that for every $k$-marked stable nodal curve $f$ in $\ov\cM_{g,k}(X,A)$, there exists at most finitely many log curves $f_{\log}\!\in\!\ov\cM^{\log}_{g,\mfs}(X,D,A)$ (with distinct decorations on the dual graph) lifting $f$; see Lemma~\ref{UniqueLog_lm}. Furthermore, $f_{\log}$ is stable if and only if $f$ is stable (and the automorphism groups are often the same). \\

\noindent
In the integrable case and in comparison with the algebraic approach, we conjecture the following statement:
\begin{conjecture}\label{conj}
In the complex algebraic setting, for any  choice of combinatorial data $\beta=(g,\mfs,A)$ and the natural log structure on $X$ associated to $D$, there is a stratified finite-to-one surjective map from the underlying space of the log moduli space $\cM(X/\tn{pt},\beta)$ in \tn{\cite{ACGS}} to $\ov\cM_{g,\mfs}^{\log}(X,D,A)$ that is one-to-one over the main stratum $\cM_{g,\mfs}(X,D,A)$.
\end{conjecture}

\noindent
\textit{In particular, this conjecture says that the group element (\ref{PLtoG_e}) mentioned in the last bullet condition above is the only  non-combinatorial obstruction for lift-ability of a nodal map (with correct combinatorial properties) to a log map (with the canonical log structures on  $X$ corresponding to $D$)}.
It is likely that we need to allow certain ``non-saturated" curves in $\cM(X/\tn{pt},\beta)$ for the conjecture to be true, or the projection map will not be surjective.
The projection map conjectured above behaves like a normalization map between varieties (e.g. unfolding self-intersections).  Based on a comparison of the coefficients of  the degeneration formula in \cite{ACGS} with our degeneration formula outlined in \cite{FT1}, we think that the degree of the projection map on each stratum should be the multiplicity $m_\Gamma$ in (\ref{TorusofGl_e}). \vskip.1in

\noindent
Similarly, in comparison with the Brett Parker approach in \cite{BP2}, under certain assumptions on the almost complex structure $J$, we expect the following statement. 
\begin{conjecture}
With respect to the exploded structure associated to an almost K\"ahler SNC divisor $D\subset X$, for any choice of combinatorial data $\beta=(g,\mfs,A)$ the ``smooth part" map gives a  finite-to-one surjective map from the moduli stack in \tn{\cite{BP2}} to $\ov\cM_{g,\mfs}^{\log}(X,D,A)$.
\end{conjecture}

\noindent 
We postpone a careful comparison of the moduli spaces constructed in this paper and those arising from \cite{AC,GS} and \cite{BP2} to a future paper.  \\

\noindent
Approaching $(\star)$, we face some new challenges that are not present in the case of the classical and relative stable maps. Unlike the smooth case, it is not a priori clear whether every SNC symplectic divisor $D\!\subset\!(X,\om)$ admits a compatible almost complex structure. Furthermore, even if $\cJ(X,D,\om)\!\neq\!\emptyset$, it is not clear whether it is contractible (or even connected).  In order to address this issue, in \cite{FMZ1}, we consider the space\footnote{this space is denoted by $\tn{Symp}^+(X,D)$ in \cite{FMZ1}.}  $\tn{Symp}(X,D)$ of all symplectic forms on $X$ such that a given transverse configuration $D=\bigcup_{i\in [N]}D_i$ is an SNC symplectic divisor in $(X,\om)$. Consequently, instead of focusing on a particular $\om$, we consider the connected component of symplectic forms in $\tn{Symp}(X,D)$ deformation equivalent to $\om$. With $\cJ(X,D,\om)$ as before, let 
$$
\cJ(X,D)= \bigcup_{\om \in \tn{Symp}(X,D)} \cJ(X,D,\om)
$$
be the space of all $D$-compatible pairs $(\om,J)$. We then define a space of almost K\"ahler auxiliary data $\AK(X,D)$ consisting of tuples $(\om,\cR,J)$ where $\om\!\in\!\tn{Symp}(X,D)$, $\cR$ is an ``$\om$-regularization" for $D$ in $X$, and $J$ is $\om$-tame and  $\cR$-compatible (which we will simply call it ``$(\cR,\om)$-compatible") almost complex structure on $X$; see Section~\ref{CptLog_ss} or \cite[p8]{FMZ1}. Roughly speaking, a regularization is a compatible set of symplectic identifications of neighborhoods of $\{D_I\}_{I\subset [N]}$ in their normal bundles with neighborhoods of them in $X$; \cite[Dfn~2.12]{FMZ1}. A regularization serves as a replacement for holomorphic defining equations in holomorphic manifolds. These regularizations are also the auxiliary data that we need to define the log tangent bundle $TX(-\log D)$. For every $(\om,\cR,J)\!\in\!\AK(X,D)$, we have $(\om,J)\!\in\!\cJ(X,D)$. Therefore, $\AK(X,D)$ is essentially a nice subset of $\cJ(X,D)$ consisting of those almost complex structures that are of  some specified type in a sufficiently small neighborhood of $D$. These special almost complex structures are similar to the almost complex structures with translational symmetry considered in \cite{LR} and in SFT \cite{EGH}.  By \cite[Thm 2.13]{FMZ1}, the forgetful map 
\bEq{AKtoSymp_e}
\AK(X,D)\lra \Symp(X,D),\qquad  (\om,\cR,J)\!\lra\!\om,
\eEq
is a weak homotopy equivalence. This implies that any invariant of the deformation equivalence classes in $\AK(X,D)$ is an invariant of the symplectic deformation equivalence class of $(X,D,\om)$. In particular, by restricting to the subclass $\AK(X,D)$, the last statement in $(\star)$ follows from constructing Kuranishi structures for families. \\

\noindent
The main goal of this paper is to prove the following compactness result, addressing the first part of ($\star$). We will address the rest in subsequent papers. We will briefly outline our approach to the deformation theory and gluing in Sections~\ref{EXDim_ss} and \ref{Gluing_ss}.

\begin{definition}
A continuous function $f\colon\! M\!\lra\! N$ between two topological spaces is a \textbf{local embedding} if for all $x\!\in\!M$ there is an open neighborhood $U\!\ni\!x$ such that $f|_U\colon\! U\!\lra\! N$ is an embedding.
\end{definition}

\noindent
By Smirnov's theorem, every paracompact, Hausdorff, and locally metrizable space is metrizable.
Therefore, if $f\colon\! M\!\lra\! N$ is a local embedding from a compact Hausdorff space $M$ to a compact metrizable space $N$ then $M$ is metrizable. 
\bTh{Compactness_th}
Assume $X$ is a compact symplectic manifold, $D\!=\!\bigcup_{i\in [N]} D_i\!\subset\! X$ is an SNC symplectic divisor.  If $(\om,\cR,J)\!\in\!\tn{AK}(X,D)$ or if $(X,D,\om,J)$ is K\"ahler, then
for every $A\!\in\!H_2(X,\Z)$, $g,k\!\in\!\N$,  and $\mfs\!\in\!\big(\Z^N\big)^k$, the Gromov sequential convergence topology on $\ov\cM_{g,k}(X,A)$ lifts to a compact Hausdorff sequential convergence topology
on $\ov\cM^{\log}_{g,\mfs}(X,D,A)$ such that the natural forgetful map 
\bEq{FogetLog_e}
\iota \colon \ov\cM^{\log}_{g,\mfs}(X,D,A)\lra \ov\cM_{g,k}(X,A)
\eEq
is a local embedding. In particular, $\ov\cM^{\log}_{g,\mfs}(X,D,A)$ is metrizable.
If $g\!=\!0$, then (\ref{FogetLog_e}) is a global embedding.
\eTh

\noindent
In other words, the open sets of $\ov\cM^{\log}_{g,\mfs}(X,D,A)$ are  the components of the intersection of open sets in $\ov\cM_{g,k}(X,A)$ with the image of $\ov\cM^{\log}_{g,\mfs}(X,D,A)$.  

\bRm{GeneralJ_rmk}
Except for the proof of Proposition~\ref{VertexOrder_prp}, every other statement in the proof of Theorem~\ref{Compactness_th}  is stated and proved for arbitrary $(\om,J)\!\in\!\cJ(X,D)$. We expect the local statement of Proposition~\ref{VertexOrder_prp}, and thus Theorem~\ref{Compactness_th}, to be true for a larger class of almost K\"ahler structures weakly homotopy equivalent to $\tn{Symp}(X,D)$ that includes both \tn{AK}(X,D) and the space of K\"ahler structures. If $D$ is smooth, a significantly simpler version of  Proposition~\ref{VertexOrder_prp} is sufficient for proving Proposition~\ref{LGSmoothV_prp}, and thus  Theorem~\ref{Compactness_th} for arbitrary $(\om,J)\!\in\!\cJ(X,D)$; see Remark~\ref{onSTNS1_rmk}. Nevertheless, by  the argument around (\ref{AKtoSymp_e}), the sub-class $\tn{AK}(X,D)$ is ideal for defining GW-type invariants  and the holomorphic case is sufficient for most of interesting examples/calculations. 
\eRm

\bRm{Gross_rmk}
While $\ov\cM^{\log}_{g,\mfs}(X,D,A)$ is defined for arbitrary $\mfs\!\in\!(\Z^N)^{k}$ satisfying the second identity in (\ref{Equiv-mfs-conditions_e}) and the compactness result holds for every such $\mfs$, 
the resulting moduli spaces do not have some of the nice properties unless $\mfs\!\in\!(\N^N)^{k}$; e.g. the (virtually) main stratum $\cM_{g,\mfs}(X,D,A)$ would be empty if any of $s_{ai}$ is negative. For $\mfs\!\in\!(\N^N)^k$, by Lemma~\ref{ExpectedGamma_lmm},  the expected dimension of $\ov\cM^{\log}_{g,\mfs}(X,D,A)$ is equal to (\ref{dlog_e}), and the only stratum with the top expected dimension is $\cM_{g,\mfs}(X,D,A)$. As pointed out to the author by M. Gross,
the case where $s_{ai}$ could be negative is called  ``punctured curves" in the work in progress \cite{ACGS2}. One feature of these punctured curves is that the moduli spaces may not carry a VFC, as even in the unobstructed case the moduli space may have irreducible components of different dimension.
\eRm

\noindent
If $D$ is smooth, we show in Proposition~\ref{ReltoLog_prp} that there is a surjective projection map 
$$
\ov\cM^{\tn{rel}}_{g,\mfs}(X,D,A)\!\lra\!\ov\cM^{\log}_{g,\mfs}(X,D,A).
$$ 
This is as expected, since our notion of log $J$-holomorphic curve involves more $\C^*$-quotients on the set of meromorphic sections than in the relative case. 
In the algebraic case, \cite[Thm~1.1]{AMW} shows that an algebraic analogue of this projection map induces an equivalence of the virtual fundamental classes.
We expect the same to hold for invariants/VFCs arising from our log moduli spaces.\\

\noindent
Approaching the rest of $(\star)$, the transversality issue aside, log moduli spaces constructed in this paper are often virtually singular in the sense that the (virtual) normal cone of each stratum is not necessarily an orbibundle. More precisely,  $\ov\cM^{\log}_{g,\mfs}(X,D,A)$ admits a stratification
$$
\ov\cM^{\log}_{g,\mfs}(X,D,A)=\bigcup_{\Gamma} \cM_{g,\mfs}(X,D,A)_\Gamma
$$
where $\Gamma$ runs over all the possible ``decorated dual graphs"; see Definition~\ref{DecoratedGamma_dfn}. For any $f$ in $\cM_{g,\mfs}(X,D,A)_\Gamma$,  the natural process of describing a neighborhood of $f$ in $\ov\cM^{\log}_{g,\mfs}(X,D,A)$ is by, first, describing a neighborhood $U$ of $f$ in $\cM_{g,\mfs}(X,D,A)_\Gamma$ and then extending that by a ``gluing" theorem of smoothing the nodes to a neighborhood of the form $U \times \cN_\Gamma'$ for $f$ in  $\ov\cM^{\log}_{g,\mfs}(X,D,A)$, where $\cN_\Gamma'$ is a neighborhood of the origin in an affine sub-variety $\cN_\Gamma\!\subset\!\C^m$, for some $m\!\in\!\N$. In this situation, we say that $\cN_\Gamma$ is the \textbf{normal cone} to $\cM_{g,\mfs}(X,D,A)_\Gamma$ or it is the \textbf{space of gluing parameters}. In the case of classical stable maps, $\cN_\Gamma$ is isomorphic to $\C^\E$ where $\E$ is the set of the edges of $\Gamma$ (or nodes of the nodal domain). Unlike in the classical case, for the log (or relative) moduli spaces,  $\cN_\Gamma$ could be reducible, and the normalization of $\cN_\Gamma$ might be singular as well; see Example~\ref{d1Rd22Pt_ex}. Nevertheless, we show that $\cN_\Gamma$ is (isomorphic to some finite copy of) an affine  toric verity that can be explicitly described in terms of $\Gamma$.  More precisely, let $\V$ and $\E$ be the set of vertices and edges of $\Gamma$, receptively. For each $v\!\in\!\V$, $I_v\!\subset\! [N]$ is the maximal subset such that image of the $v$-th component of $f$ lies in $D_{I_v}$. Similarly, for each $e\!\in\!\E$, $I_e\!\subset\! [N]$ is the maximal subset such that image of the $e$-th node lies in $D_{I_e}$.   In (\ref{DtoT_e}), associated to every such $\Gamma$, we construct a $\Z$-linear map 
\bEq{ZMap_e}
\vr\colon \D(\Gamma)=\Z^E\oplus \bigoplus_{v\in \V} \Z^{I_v} \lra\T(\Gamma)= \bigoplus_{e\in \E} \Z^{I_e},
\eEq
such that $\cN_\Gamma$ is isomorphic to (some finite copy of) the toric variety associated to a maximal convex rational polyhedral cone in $\tn{Ker}(\vr)\otimes \R$. Moreover, the group element mentioned in the last bullet condition of Page~\pageref{GObs_it} (i.e. in the definition of a log map) is an element of the Lie group $\mc{G}(\Gamma)$ with the Lie algebra $\tn{Coker}(\vr)\otimes \C$. In other words, $\tn{Ker}(\vr)$ gives the deformation space in the normal direction and $\tn{Coker}(\vr)$ gives an obstruction for the smoothability of such maps.

\subsection{Outline and acknowledgements} 
In Section~\ref{dbar_ss}, we review the definition and properties of $\dbar$-operators.
The $\dbar$-operator $\dbar_{\cN_XD}$ on the normal bundle $\cN_XD$ described in Lemma~\ref{ConnectionToDbar_lmm} plays a key role in defining the basic building blocks of  relative and log maps. In Section~\ref{DDG_ss}, we set up our notation for the decorated dual graph of nodal maps. The $\Z$-linear map (\ref{ZMap_e}) is defined in terms of such decorated dual graphs. In Section~\ref{LogMapRelSmooth_ss}, we define the moduli spaces of log $J$-holomorphic curves and provide several examples to highlight their features. This is done in two steps: first, in Definition~\ref{PreLogMap_dfn}, we define a straightforward notion of pre-log map. Then in Definition~\ref{LogMap_dfn}, we impose two non-trivial conditions on a such a pre-log map to define a log map. Proof of Theorem~\ref{Compactness_th} relies on Gromov's compactness result for the underlying stable maps. In Section~\ref{CptGromov_ss}, we review the Gromov compactness theorem and set up the notation for the proof of Theorem~\ref{Compactness_th}. In Section~\ref{CptLog_ss}, we state a log enhancement of the Gromov compactness theorem. Proof of the main result is done in multiple steps in Sections~ \ref{LocalTh_ss} and \ref{CptLogSmooth_ss}. The main step of the proof is Proposition~\ref{VertexOrder_prp} that compares the limiting behavior of the rescaling and gluing parameters. 
In the case of smooth divisors, we compare the relative and the log compactifications of the same combinatorial type in Section~\ref{RelToRel_ss}.  We review the construction of relative compactification in Section~\ref{RelComp_ss}. In Section~\ref{EXDim_ss}, we outline a Fredholm set up for studying the deformation  theory of log $J$-holomorphic maps and make some conclusions. This setup is extended to perturbed log maps and discussed in details in \cite{FT2}. 
In Section~\ref{Gluing_ss}, we explicitly describe the space of gluing parameters of any fixed type $\Gamma$ and identify it with an explicit affine toric variety. \\

\noindent
I am indebted to A. Zinger for many years of related collaborations that sorted out my thoughts toward this paper.  
I would like to thank K.~Fukaya and J.~Morgan for supporting my research at Simons Center and for many fruitful  conversations about the details of Kuranishi structures.
I~am thankful to G.~Tian for supporting my research and for many related conversations about virtual fundamental class.
Finally, I am also thankful to  Q.~Chen, D.A.~Cox, M.~Gross, H-J.~Hein, M.~Liu, D.~McDuff, M.~McLean, B.~Parker, D.~Pomerleano, D.~Ranganathan, H.~Ruddat, B.~Siebert, J.~Starr, and the referees for answering my questions and for their helpful comments. My research on this topic is supported by the NSF grant DMS-2003340.

%------------------------------------------------------------------------------------------------------
\section{Log pseudoholomorphic maps}\label{LogM_s}

In this section, we construct the moduli spaces of log $J$-holomorphic curves relative to an arbitrary SNC symplectic divisor defined in \cite{FMZ1}. 
This is done by first introducing a notion of pre-log $J$-holomorphic map which only involves a matching condition of contact orders at the nodes. We then define a $\Z$-linear map between certain $\Z$-modules associated to the dual graph of such a pre-log map which encodes the essential deformation/obstruction data for defining and studying log maps.\\

\noindent
Let us start by some well known facts about almost complex structures.
Let $(X,\om)$ be a smooth symplectic manifold and $J$ be an $\om$-tame almost complex structure on $X$.
Let $\na$ be the Levi-Civita connection of the metric $\ll u, v\rr=\frac{1}{2}(\om(u,Jv)+\om(v,Ju))$ and 
\bEq{CLConn_e}
\wt\na_{v} \ze\!=\!\na_v \ze - \frac{1}{2} J(\na_v J) \ze\!=\! \frac{1}{2}\big(\na_v \ze-J \na_v (J\ze)\big)\qquad \forall~v\!\in\!TX,~\ze\!\in\!\Gamma(X,TX)
\eEq
be the associated Hermitian connection. The Hermitian connection $\wt\na$ coincides with 
$\na$ if and only if $(X,\om,J)$ is K\"ahler, i.e. $\na J \!\equiv\! 0$. The torsion $T$ of the modified $\C$-linear connection
\bEq{CLConn_e2}
\wh\na_{v} \ze\!=\!\wt\na_v \ze - A(\ze)v, \quad A(\ze)= \frac{1}{4} \big(\nabla_{J\ze}J+ J\nabla_{\ze}J\big) \qquad \forall~v\!\in\!TX,~\ze\!\in\!\Gamma(X,TX)
\eEq
is related to the Nijenhueis tensor (\ref{intInnormal_e0}) by
\bEq{Torsion_e}
T_{\wh\na}(v,w)=-\frac{1}{4}N_J(v,w)\qquad\forall~v,w\!\in\!TX.
\eEq
If $J$ is $\om$-compatible, $\wt\na$ coincides with $\wh\na$. See \cite[Ch~3.1 and Appendix C]{MS2} for details.

%---------------------------------------------------
\subsection{Almost complex structures and $\dbar$-operators}\label{dbar_ss}

Suppose $M$ is a smooth manifold, $\mfi_M$ is an almost complex structure on $M$, and  $(L,\mfi_L)\!\lra\!M$ is a complex vector bundle.
Let
\bEq{FormDecomposition_e}
\Om_{M,\mfi_M}^{1,0}\equiv \{ \eta\!\in\!T^*M\otimes_\R\C \colon \eta\circ \mfi_M\! =\! \mfi \circ \eta\}\quad\tn{and}\quad 
\Om_{M,\mfi_M}^{0,1}\equiv \{ \eta\!\in\!T^*M\otimes_\R\C \colon \eta\circ \mfi_M \!=\! -\mfi \circ \eta\}
\eEq
be the bundles of $\C$-linear and $\C$-antilinear $1$-forms on $M$, where $\mfi$ is the unit imaginary number in~$\C$.
Given a smooth function $f\colon\!M\!\lra\!\C$, (\ref{FormDecomposition_e}) gives a decomposition of $\nd f$ into $\C$-linear and $\C$-antilinear parts $\partial f$ and $\dbar f$, respectively.
A \textbf{$\bar\partial$-operator} on $(L,\mfi_L)$ is a complex linear operator 
\bEq{barPartialL_e2}
\bar\partial\colon \Gamma(M,L)\lra \Gamma (M, \Om_{M,\mfi_M}^{0,1}\otimes_\C L)
\eEq
such that 
$$
\bar\partial (f\ze)\!=\! \ov\partial f\otimes \ze\!+\! f \bar\partial \ze\qquad \forall f\!\in\!C^\infty(M,\C),~\ze\!\in\!\Gamma(M,L).
$$
Given a complex linear connection $\nabla$ on $(L,\mfi_L)\!\lra\!(M,\mfi_M)$, the $(0,1)$-part 
\bEq{ConnToNa_e}
\nabla^{(0,1)}=\frac{1}{2}(\nabla+ \mfi_{L} \nabla \circ \mfi_M)
\eEq
of $\nabla$ is a $\dbar$-operator which we denote by $\dbar_{\nabla}$. Every $\dbar$-operator is the associated $\dbar$-operator of some $\C$-linear connection $\na$ as above. The connection, however, is not uniquely determined. Every two connections $\nabla$ and $\nabla'$ differ by a global $\tn{End}(L)$-valued $1$-form $\al$, i.e. $\nabla'\!=\!\nabla\!+\!\al$. If $\nabla$ and $\nabla'$ are complex linear connections on $(L,\mfi_L)$ with $\nabla'=\nabla\!+\!\al$, then 
$$
\bar\partial_{\nabla'} \!=\! \bar\partial_{\nabla} \!+\! \al^{(0,1)},
$$
where $\al^{(0,1)}$ is the $(0,1)$-part of $\al$ in the decomposition (\ref{FormDecomposition_e}). In particular, $\bar\partial_{\nabla'} \!=\! \bar\partial_{\nabla}$ whenever $\al$ is of $(1,0)$ type. \\

\noindent
By \cite[Lmm 2.2]{Z}, corresponding to every $\bar\partial$-operator (\ref{barPartialL_e2}), there exists a unique almost complex structure $J \!=\! J_{\bar\partial}$ on the total space of $L$ such that 
\bEn
\item\label{HoloProj_l} the projection $\pi\colon\!L\!\lra\!M$ is a $(\mfi_M,J)$-holomorphic map (i.e. $\nd\pi\!+\!\mfi \,\nd\pi\, J\!=\!0$), 
\item\label{FiberJ_l} the restriction of $J$ to the vertical tangent bundle $TL^{\tn{ver}}\cong \pi^*L\!\subset\!TL$ agrees with $\mfi_L$,
\item\label{HoloSec_l} and the map $\ze\colon\! M\!\lra\!L$ corresponding to a section $\ze\!\in\!\Gamma(M,L)$ is $(J,\mfi_M)$-holomorphic if and only if $\bar\partial \ze\!=\!0$. 
\eEn

\noindent
Suppose $(X,\om)$ is a symplectic manifold, $D$ is a symplectic submanifold, and $J$ is an $\om$-tame almost complex structure on $X$ such that $J(TD)\!=\!TD$. The last condition implies that $J$ induces a complex structure~$\mfi_{\cN_XD}$ on 
(the fibers~of) the normal bundle
\bEq{Firstpi_e}
\pi\!: \cN_XD\equiv TX|_D\big/TD\lra D.
\eEq
Under the isomorphism
$$
\cN_XD\cong TD^{\perp}\!=\!\{u\in TX|D\colon~\ll u,v\rr\!=\!0\quad\forall~v\!\in\!TD\},
$$
$\mfi_{\cN_XD}$ is the same as the restriction to $TD^\perp$ of $J$.
Let $J_D$ denote the restriction of $J$ to $TD$.

\bLm{ConnectionToDbar_lmm}
Suppose $(X,\om)$ is a symplectic manifold, $D$ is a symplectic submanifold, $J$ is an $\om$-tame almost complex structure on $X$ such that $J(TD)\!=\!TD$, and $\wh\na$ is the $\C$-linear connection associated to $(\om,J)$ in (\ref{CLConn_e2}). Then the $\dbar$-operator 
$$
\dbar_{\wh\na}\equiv \wh\na^{(0,1)}\colon \Gamma\big(X,TX\big)\lra \Gamma\big(X,\Om_{X,J}^{0,1}\!\otimes_{\C}\!TX\big)
$$
in (\ref{ConnToNa_e}) descends to a $\dbar$-operator 
\bEq{dbarcNXV_e}
\dbar_{\cN_XD}\colon \Gamma\big(D,\cN_XD\big)\lra \Gamma\big(D,\Om_{D,J_D}^{0,1}\!\otimes_{\C}\!\cN_XD\big)
\eEq
on $(\cN_XD,\mfi_{\cN_XD})\!\lra\!(D,J_D)$.
\eLm

\bPf
We need to show that $\dbar_{\wh\na}$ maps $\Gamma\big(D,TD\big)$ to $\Gamma\big(D,\Om_{D,J_D}^{0,1}\!\otimes_{\C}\!TD\big)$.
Let $\na$ and $\na^D$ be the Levi-Civita connections of the metics associated to $(\om,J)$ and $(\om|_{TD},J_D)$ on $X$ and $D$, respectively. Then 
$$
\na\ze = \na^{D} \ze + \na^{N}\ze\qquad \forall~\ze\!\in\!\Gamma(D,TD),
$$
such that 
$$
\na^{N}\ze\!\in\!\Gamma\big(D,\Om_{D}^1\!\otimes \!TD^\perp\big).
$$ 
Similarly, let $\wt\na$ and $\wt\na^D$ be the Chern connections on $TX$ and $TD$ associated to $\na$ and $\na^D$ as in  (\ref{CLConn_e}), respectively. 
It follows from  (\ref{CLConn_e}) that 
\bEq{nawtnawhna_e}
\wt\na \ze = \wt\na^D \ze + \wt\na^N \ze\qquad \forall~\ze\!\in\!\Gamma(D,TD),
\eEq
where
$$
\wt\na^N\ze\!=\!\frac{1}{2}\big(\na^{N} \ze - J\na^{N} (J\ze)\big)\!\in\!\Gamma\big(D,\Om_{D}^1\!\otimes_\C \!TD^\perp\big).
$$  
Let $\wh\nabla$ and $\wh\nabla^D$ be the modifications of $\wt\na$ and $\wt\nabla^D$  as in (\ref{CLConn_e2}), respectively. By  (\ref{CLConn_e2}) and (\ref{nawtnawhna_e}), we also have 
\bEq{nawtnawhna_e2}
\wh\na_\xi \ze = \wh\na_\xi^D \ze + {\wh\na}_\xi^N \ze\qquad \forall~\xi,\ze\!\in\!\Gamma(D,TD),
\eEq
where 
$$
\aligned
&\wh\na^N_\xi\ze\!=\!\wt\na^N_\xi\ze\!-\! A^N(\ze)\xi \in\!\Gamma\big(D,\Om_{D}^1\!\otimes_\C \!TD^\perp\big)\\ 
& A^N(\ze)\xi= \frac{1}{4} \big(\nabla^N_{J\ze}J+ J\nabla^N_{\ze}J\big)\xi, 
\quad (\nabla^{N}_{\ze}J)\xi:=\nabla^N_{\ze}(J\xi)-J\nabla^N_{\ze}\xi.
\endaligned
$$  
From (\ref{Torsion_e}), (\ref{nawtnawhna_e2}), and  
$$
N_J(\ze,\xi)\!=\!N_{J_D}(\ze,\xi)\!\in\!TD\qquad\forall~\xi,\ze\!\in\!\Gamma(D,TD),
$$
we conclude that 
$$
\aligned
\wh\na^N_\xi\ze-\wh\na^N_\ze\xi&= \big(\wh\na_\xi\ze-\wh\na_\ze\xi\big)-\big(\wh\na^D_\xi\ze-\wh\na^D_\ze\xi\big)=\\
&\big(\wh\na_\xi\ze-\wh\na_\ze\xi-[\ze,\xi]\big)-\big(\wh\na^D_\xi\ze-\wh\na^D_\ze\xi-[\ze,\xi]\big)=T_{\wh\na}(\ze,\xi)-T_{\wh\na^D}(\ze,\xi)=0;
\endaligned
$$
i.e. 
$$
\wh\na^N_\xi\ze = \wh\na^N_\ze\xi\qquad \forall \ze,\xi\!\in\!\Gamma(D,TD).
$$
From the last identity we get 
$$
\wh\nabla^N_\xi \ze+ J \wh\nabla^N_{J\xi} \ze=\wh\nabla^N_\ze \xi+ J \wh\nabla^N_{\ze} J\xi=\wh\nabla^N_\ze \xi- \wh\nabla^N_{\ze} \xi=0\qquad \forall~\ze,\xi\!\in\!\Gamma(D,TD).
$$
Therefore,
$$
\aligned
\wh\na^{(0,1)}_\xi\ze=&\frac{1}{2}\big(\wh\nabla_\xi \ze+ J \wh\nabla_{J\xi} \ze\big)=\frac{1}{2}\big(\wh\nabla^D_\xi \ze+ J \wh\nabla^D_{J\xi} \ze\big)+ \frac{1}{2}\big(\wh\nabla^N_\xi \ze+ J \wh\nabla^N_{J\xi} \ze\big)=\\
&\frac{1}{2}\big(\wh\nabla^D_\xi \ze+ J \wh\nabla^D_{J\xi} \ze\big)=\wh\na^{D,(0,1)}_\xi\ze\!\in\!\Gamma(D,TD)\qquad \forall~\ze,\xi\!\in\!\Gamma(D,TD).
\endaligned
$$
\ePf
\begin{remark}
The term $A(\ze)v$ in (\ref{CLConn_e2}) is $\C$-linear in $\ze$ and $\C$-antilinear in $v$. It vanishes if $J$ is $\om$-compatible.
Therefore,
$$
\wh\na^{0,1}_{\xi}\ze=\wt\na^{0,1}_{\xi}\ze - A(\ze)\xi.
$$
\end{remark}

%---------------------------------------------------
\subsection{Decorated dual graphs}\label{DDG_ss}

Let $\Gamma\!=\!\Gamma(\V,\E,\L)$ be a graph with the set of vertices $\V$, edges $\E$, and legs $\L$; the latter, also called flags or roots, are half edges that have a vertex at one end and are open at the other end. Let $\uvec{\E}$ be the set of edges with an orientation. Given an oriented edge $\uvec{e}\!\in\!\uvec{\E}$, let $\ucev{e}$ denote the same edge $e$ with the opposite orientation. For each $\uvec{e}\!\in\!\uvec{\E}$, let $v_1(\uvec{e})$ and $v_2(\uvec{e})$ in $\V$ denote the starting and ending points of the arrow, respectively. For $v,v'\!\in\!\V$, let $\E_{v,v'}$ denote the subset of edges between the two vertices and $\uvec{\E}_{v,v'}$ denote the subset of oriented edges from $v$ to $v'$. For every $v\!\in\!\V$, let $\uvec{\E}_{v}$ denote the subset of oriented edges starting from $v$.\\

\noindent
A \textbf{genus labeling} of $\Gamma$ is a function $g\colon\!\V\!\lra\! \N$. An \textbf{ordering of the legs} of $\Gamma$ is a bijection $a\colon \!\L\!\lra\! \{1,\ldots, |\L|\}$.
If a decorated graph $\Gamma$ is connected, the \textbf{arithmetic} genus of $\Gamma$ is 
\bEq{genus_e}
g=g_\Gamma= \sum_{v\in \V} g_v\! +\! \tn{rank} ~H_1(\Gamma,\Z),
\eEq 
where $H_1(\Gamma,\Z)$ is the first homology group of the underlying topological space of $\Gamma$. Figure~\ref{labled-graph_fg}-left illustrates a labeled graph with $2$ legs. \\

\begin{figure}
\begin{pspicture}(8,1.3)(11,3)
\psset{unit=.3cm}

\pscircle*(35,6){.25}\pscircle*(39,6){.25}
\pscircle*(35,9){.25}\pscircle*(41.5,9){.25}\pscircle*(39,9){.25}
\psline[linewidth=.05](39,6)(35,9)
\psline[linewidth=.05](39,6)(39,9)\psline[linewidth=.05](39,6)(41.5,9)
\psarc[linewidth=.05](33,7.5){2.5}{-36.9}{36.9}\psarc[linewidth=.05](37,7.5){2.5}{143.1}{216.0}
\psline[linewidth=.05](41.5,9)(40.5,10.5)\psline[linewidth=.05](41.5,9)(42.5,10.5)
\rput(40.5,11.2){1}\rput(42.5,11.2){2}
\rput(35,5){$g_4$}\rput(39.5,5){$g_5$}
\rput(34.5,10){$g_1$}\rput(42.3,8.5){$g_3$}
\rput(39,10){$g_2$}

\pscircle*(55,6){.25}\pscircle*(59,6){.25}
\pscircle*(55,9){.25}\pscircle*(61.5,9){.25}\pscircle*(59,9){.25}
\psline[linewidth=.05](59,6)(55,9)
\psline[linewidth=.05](59,6)(59,9)\psline[linewidth=.05](59,6)(61.5,9)
\psarc[linewidth=.05](53,7.5){2.5}{-36.9}{36.9}\psarc[linewidth=.05](57,7.5){2.5}{143.1}{216.0}
\psline[linewidth=.05](61.5,9)(60.5,10.5)\psline[linewidth=.05](61.5,9)(62.5,10.5)
\rput(60.5,11.2){1}\rput(62.5,11.2){2}
\rput(55,5){$(g_4,A_4)$}\rput(59.5,5){$(g_5,A_5)$}
\rput(54,10){$(g_1,A_1)$}\rput(63.8,8.5){$(g_3,A_3)$}
\rput(58.3,10){$(g_2,A_2)$}
\end{pspicture}
\caption{On left, a labeled graph $\Ga$ representing elements of $\ov\cM_{g,2}$. On right, a labeled graph $\Ga$ representing elements of $\ov\cM_{g,2}(X,A)$.}
\label{labled-graph_fg}
\end{figure}

\noindent 
Such decorated graphs $\Gamma$ characterize different topological types of nodal marked surfaces 
$$
(\Si,\vec{z}\!=\!(z^1,\ldots,z^k))
$$ 
in the following way.
Each vertex $v\!\in\!\V$ corresponds to a smooth\footnote{We mean a smooth closed oriented surface.} component $\Si_v$ of $\Si$ with genus $g_v$. Each edge $e\!\in\!\E$ corresponds to a node $q_e$ obtained by connecting $\Sigma_v$ and $\Sigma_{v'}$ at the points $q_{\uvec{e}}\!\in\!\Si_v$ and $q_{\scz\ucev{e}}\!\in\!\Si_{v'}$, where $e\!\in\!\E_{v,v'}$ and $\uvec{e}$ is an orientation on $e$ with $v_1(\uvec{e})\!=\!v$. The last condition uniquely specifies $\uvec{e}$ unless $e$ is a loop connecting $v$ to itself. Finally, each leg $l\!\in\!\L$ connected to the vertex $v_l$ corresponds to a marked point $z^{a_l}\!\in \Sigma_{v_l}$ disjoint from the connecting nodes. If $\Si$ is connected, then $g_\Gamma$ is the arithmetic genus of $\Si$. Thus we have
\bEq{nodalcurve_e}
(\Si,\vec{z})\! =\! \coprod_{v\in \V}(\Si_v,\vec{z}_v,{q}_v)/\sim, \quad  q_{\uvec{e}}\!\sim\! q_{\scz\ucev{e}}\quad\forall~e\!\in\!\E,
\eEq
where 
$$
\vec{z}_v\!=\!\vec{z}\cap \Sigma_v\quad\tn{and}\quad \quad {q}_v=\{q_{\uvec{e}}\colon \uvec{e}\!\in\!\uvec{\E}_v\}\qquad \forall~v\!\in\!\V.
$$
In this situation, we say $\Gamma$ is the \textbf{dual graph} of $(\Si,\vec{z})$.
We treat $q_v$ as an un-ordered set of marked points on $\Si_v$. If we fix an ordering on the set $q_v$, we denote the ordered set by $\vec{q}_v$.
\\

\noindent
A complex structure $\mfj$ on $\Sigma$ is a set of complex structures $(\mfj_v)_{v\in \V}$ on its components. By a (complex) marked nodal curve, we mean a  marked nodal real surface together with a complex structure $(\Si,\mfj,\vec{z})$. 
Figure~\ref{labled-curve_fg} illustrates a nodal curve with $(g_1,g_2,g_3,g_4,g_5)=(0,2,0,1,0)$ corresponding to Figure~\ref{labled-graph_fg}-left. \\

\noindent 
Similarly, for nodal marked surfaces mapping into a topological space $X$, we consider similar decorated graphs where the vertices carry an additional \textbf{degree labeling} 
$$
A\colon \V\lra H_2(X,\Z),\quad v\lra A_v,
$$
recording the homology class of the image of the corresponding component. Figure~\ref{labled-graph_fg}-right illustrates a dual graph associated to a marked nodal map over the graph on the left. \\

\begin{figure}
\begin{pspicture}(8,1.7)(11,4.5)
\psset{unit=.3cm}

\pscircle*(52.4,8.2){.25}\pscircle*(52.4,11.8){.25}
\pscircle(50,10){3}\psellipse[linestyle=dashed,dash=1pt](50,10)(3,1)
\pscircle(54,7){2}\psellipse[linestyle=dashed,dash=1pt](54,7)(2,.77)
\pscircle*(55,8.6){.25}\rput(56,9){$z^1$}
\pscircle*(55,5.4){.25}\rput(56,5.4){$z^2$}

\pscircle(45,10){2}\psellipse[linestyle=dashed,dash=1pt](45,10)(2,.77)
\pscircle*(47,10){.25}\pscircle*(43.4,8.8){.25}\pscircle*(43.4,11.2){.25}
\psarc(41.8,7.6){2}{270}{71}
\psarc(41.8,12.4){2}{-71}{90}
\psarc(41.8,10){4.4}{90}{270}\psarc(42.8,10){.6}{120}{240}
\psarc(40,8.7){1.5}{60}{120}
\psarc(40,11.2){1.5}{225}{315}

\psarc(54,13){2}{60}{300}
\psarc(58,13){2}{240}{120}
\psarc(56,16.464){2}{240}{300}
\psarc(56,9.546){2}{60}{120}
\psarc(54,11.7){1.5}{60}{120}
\psarc(54,14.2){1.5}{225}{315}
\psarc(58,11.7){1.5}{60}{120}
\psarc(58,14.2){1.5}{225}{315}

\end{pspicture}
\caption{A nodal curve in $\ov\cM_{4,2}$.}
\label{labled-curve_fg}
\end{figure}

\noindent
Assume $D\!=\!\bigcup_{i\in [N]} D_i \!\subset\! X$ is an SNC symplectic divisor, $(\om,J)\!\in\!\cJ(X,D)$, and $(\Si,\mfj)$ is a connected smooth complex curve. Then every $J$-holomorphic map $u\colon\! (\Si,\mfj)\!\lra\!(X,J)$ has a well-defined \textbf{depth} $I\!\subset\![N]$, which is the maximal subset of $[N]$ such that $\tn{Image}(u)\!\subset\!D_I$. In particular, any map $u$ intersecting $D$ in a discrete set is of depth $I\!=\!\eset$. We say a point $x\!\in\!\Si$ is of \textbf{depth $I$}, if $D_I$ is the minimal stratum containing $u(x)$. Let $\cP(N)$ be the set of subsets of $[N]$. The dual graph of $(u,\Si)$ carries additional labelings
\bEq{LabelI_e}
I\colon \V,\E\lra \cP(N), \qquad v\lra I_v\quad \forall v\!\in\!\V, \qquad e\lra I_e\quad \forall e\!\in\!\E,
\eEq
recording the depths of smooth components and nodes of $\Si$. 

%---------------------------------------------------
\subsection{Log moduli spaces}\label{LogMapRelSmooth_ss}
Assume $D\!= \!\bigcup_{i\in [N]} D_i \!\subset\! X$ is an SNC symplectic divisor, $(\om,J)\!\in\!\cJ(X,D)$, and $u\colon\!(\Si,\mfj)\!\lra\!(X,J)$ is a $J$-holomorphic map of depth $I\!\subset\![N]$ with smooth domain. Then, for every $i\!\in\![N]\!-\!I$, the function
\bEq{Ordx_e1}
\ord_{u}^i\colon\Si\lra \N, \quad \ord_{u}^i(x)=\tn{ord}_x(u,D_i),
\eEq
recording the contact order of $u$ with $D_i$ at $x$ is well-defined. 
For every $i\!\in\!I$, let $u^*\dbar_{\cN_XD_i}$ be the pull-back of $\dbar$-operator $\dbar_{\cN_XD_i}$ associated to $(J,D_i)$ in (\ref{dbarcNXV_e}). Since every $\dbar$-operator over a complex curve is integrable,  $u^* \dbar_{\cN_XD_i}$ defines a holomorphic structure on $u^*\cN_XD_i$; see \cite[Rmk C.1.1]{MS2}.
The holomorphic line bundles
$$
\big(u^*\cN_XD_i, u^* \dbar_{\cN_XD_i}\big)\qquad \forall~i\!\in\!I
$$ 
play a key role in definition of the log moduli space below. Let $\Om_{\tn{mero}}(\Si,u^*\cN_XD_i)$ be the space of non-trivial meromorphic sections of $u^*\cN_XD_i$ with respect to $u^* \dbar_{\cN_XD_i}$; $\C^*$ acts on $\Om_{\tn{mero}}(\Si,u^*\cN_XD_i)$ by multiplication. We denote the $\C^*$-equivalence class of a section 
$$
\ze\!\in\!\Om_{\tn{mero}}(\Si,u^*\cN_XD_i)
$$ by $[\ze]$. The function 
\bEq{Ordx_e2}
\ord_{[\ze]}\colon\Si\lra \Z, \quad \ord_{[\ze]}(x)= \tn{ord}_x(\ze),
\eEq
recording the vanishing order of $\ze$ at $x$ (which is negative if $\ze$ has a pole at $x$) is well-defined.\\

\noindent
A \textbf{log $J$-holomorphic tuple} $(u,[\ze], \Si,\mfj,w)$ consists of a smooth (closed) connected curve $(\Si,\mfj)$, $\ell$ distinct points $w\!=\!\{w^1,\ldots,w^\ell\}$ on $\Si$,  a $(J,\mfj)$-holomorphic map $u\colon\!(\Si,\mfj)\!\lra\!(X,J)$ of depth $I\!\subset\![N]$, and 
\bEq{merosection_e}
[\ze]\equiv \big([\ze_i]\big)_{i\in I}\in    \prod_{i\in I} \Big(\Om_{\tn{mero}}(\Si,u^*\cN_XD_i)/\C^*\Big)
\eEq
such that 
\bEq{Intersection-Condition_e}
\ord_{u,[\ze]}(x)\!\neq \!0~~\Rightarrow~~x\!\in\!w, \quad \forall~x\!\in\!\Si,
\eEq
where the vector-valued \textbf{order function}
$$
\ord_{{u,[\ze]}}(x)\!=\!\big((\ord_{u}^i(x))_{i\in [N]-I},(\ord_{[\ze_i]}(x))_{i\in I}\big) \!\in\!\Z^N \qquad \forall~x\!\in\!\Si
$$ is defined via (\ref{Ordx_e1}) and (\ref{Ordx_e2}).
\noindent
In particular, if $u$ is of degree $A\!\in\!H_2(X,\Z)$, then (\ref{Intersection-Condition_e}) implies
\bEq{Intersection-Condition-on-z_e}
(A\cdot D_i)_{i\in [N]}=\sum_{w^a\in w} \ord_{u,[\ze]}(w^a) \in \Z^{N}.
\eEq

\bRm{UniqueLog_rmk}
For every $J$-holomorphic map $u\colon (\Si,\mfj)\!\lra\!(X,J)$ with smooth domain, $\ell$ distinct points $w^1,\ldots,w^\ell$  in $\Si$, and $s_1,\ldots, s_\ell\!\in\!\Z$, if $\tn{Im}(u)\!\subset\!D_i$, up to $\C^*$-action there exists at most one meromorphic section $\ze_i\!\in\!\Om_{\tn{mero}}(\Si,u^*\cN_XD_i)$ with zeros/poles of order $s_a$ at $w^a$ (and nowhere else).
\eRm

\bDf{PreLogMap_dfn} 
Let $D\!= \!\bigcup_{i\in [N]} D_i \!\subset\! X$ be an SNC symplectic divisor, $(\om,J)\!\in\!\cJ(X,D)$, and 
$$
C\!\equiv\!(\Si,\mfj,\vec{z})=\bigg(\coprod_{v\in \V} C_v\equiv(\Si_v,\mfj_v,\vec{z}_v,q_v)\bigg)/\sim,\quad   q_{\uvec{e}}\!\sim\! q_{\scz\ucev{e}}\quad\forall~\uvec{e}\!\in\!\uvec{\E},
$$ 
be a $k$-marked connected nodal  curve with smooth components $C_v$ and dual graph $\Gamma\!=\!\Gamma(\V,\E,\L)$ as in (\ref{nodalcurve_e}). A \textbf{pre-log} $J$-holomorphic map of contact type $\mfs\!=\!(s_a)_{a=1}^k\!\in\!\big(\Z^N\big)^k$ from $C$ to $X$ is a collection 
\bEq{fplogSetUp_e}
f\equiv \big(f_{v}\!\equiv\! (u_v,[\ze_v],C_v)\big)_{v\in \V}
\eEq
such that  
\bEn
\item for each $v\!\in\!\V$, $(u_v,[\ze_v]=([\ze_{v,i}])_{i\in I_v},\Si_v,\mfj_v,z_v \cup q_v)$ is a log $J$-holomorphic tuple,
\item\label{MatchingPoints_l} $u_v(q_{\uvec{e}})\!=\!u_{v'}(q_{\scz\ucev{e}})\!\in\! X$ for all $\uvec{e}\!\in\! \uvec{\E}_{v,v'}$;
\item\label{MatchingOrders_l} $s_{\uvec{e}}\!\equiv\! \ord_{u_v,\ze_v}(q_{\uvec{e}})\!=\! -\ord_{u_{v'},\ze_{v'}}(q_{\scz\ucev{e}})\!\equiv\!-s_{\scz\ucev{e}}$ for all $v,v'\!\in\!\V$ and $\uvec{e}\!\in\! \uvec{\E}_{v,v'}$;
\item\label{OatMkd_l} $\tn{ord}_{u_v,\ze_v}(z^a)=s_a$ for all $v\!\in\!\V$ and $z^a\!\in\!z_v$.
\eEn
\eDf

\noindent
\textit{In other words, a pre-log map is a nodal $J$-holomorphic map with a bunch of meromorphic sections on each smooth component, opposite contact orders at the nodes, and prescribed contact orders at the marked points.}

\noindent
\bRm{IevsIv_rmk}
For every $v\!\in\!\V$ and $\uvec{e}\!\in\!\uvec{\E}_v$, let
\bEq{seDf_e}
s_{\uvec{e}}=(s_{\uvec{e},i})_{i\in [N]}=\big((\ord_{u_v}^i(q_{\uvec{e}}))_{i\in [N]-I_v},(\ord_{[\ze_{v,i}]}(q_{\uvec{e}}))_{i\in I_v}\big) \!\in\!\Z^{N}
\eEq
be the contact order data at the nodal point $q_{\uvec{e}}\!\in\!\Si_v$.
For $\uvec{e}\!\in\! \uvec{\E}_{v,v'}$, if $u_v$ and $u_{v'}$ have image in $D_{I_{v}}$ and $D_{I_{v'}}$, respectively, by Condition~\ref{MatchingPoints_l} above, we have
$$
u(q_e)\!=\!u_v(q_{\uvec{e}})\!=\!u_{v'}(q_{\scz\ucev{e}})\in D_{I_{v}}\cap D_{I_{v'}}= D_{I_{v}\cup I_{v'}};
$$
i.e. $I_{e}\!\supset\!I_{v}\!\cup\! I_{v'}$. If $i\!\in\![N]\!\setminus\!I_{v}\!\cup\! I_{v'}$, by (\ref{Ordx_e1}), we have 
$$
s_{\uvec{e},i},s_{\scz\ucev{e},i}\geq 0.
$$
Therefore, by Condition~\ref{MatchingOrders_l} above, they are both zero, i.e.
 \bEq{Ieunion_e}
 I_{e}\!=\!I_{v}\!\cup\! I_{v'}\quad\tn{and}\quad s_{\uvec{e}}\in \Z^{I_{e}} \!\times\! \{0\}^{[N]-I_{e}} \subset \Z^{N}\qquad \forall~\uvec{e}\!\in\! \uvec{\E}_{v,v'}.
 \eEq
 \eRm
 \vskip.1in
 
 \noindent
The dual graph $\Gamma$ of every pre-log map in Definition~\ref{LogMap_dfn} carries an additional decoration $s_{\uvec{e}}\in \Z^{N}$, for all $\uvec{e}\!\in\!\uvec{\E}$,
which records the contact order of $(u_v,[\ze_v])$ at the nodal point $q_{\uvec{e}}\!\in\!\Si_v$ for every $\uvec{e}\!\in\!\uvec{\E}_v$; see Figure~\ref{Freeline_fig}. The set $\L$ of legs of $\Gamma$ is  also decorated with the vector-valued contact order function 
$$
\tn{ord}\colon\!\L \!\lra\!\Z^N, \quad l \lra s_l,
$$
recording the contact vector at the marked point $z^{a_l}$ corresponding to $l$.\\

\noindent
Two pre-log maps $(u,[\ze],C)\!\equiv\! \big(u_v,[\ze_v],C_v\big)_{v\in \V}$ and $(\wt{u},[\wt{\ze}],\wt{C})\!\equiv\!\big(\wt{u}_v,[\wt{\ze}_v],\wt{C}_v\big)_{v\in \V}$  with isomorphic decorated dual graphs $\Gamma$ as in Definition~\ref{PreLogMap_dfn} are equivalent if there exists a biholomorphic identification 
\bEq{NodalEquivalence_e}
(h\colon \wt{C}\lra C)\equiv \big(h_v\colon (\wt\Si_v,\wt\mfj_v)\!\lra\!(\Si_{h(v)},\mfj_{h(v)})\big)_{v\in \V}
\eEq \vskip-.1in
\noindent
such that
$$
h(\wt{z}^a)\!=\!z^a\quad \forall a\!=\!1,\ldots,k,\quad u\circ h\!=\!\wt{u},\quad [h_v^*\ze_{h(v),i}]=[\wt{\ze}_{v,i}]\quad \forall~v\!\in\!\V,~i\!\in\!I_v.
$$
A pre-log map $f$ is \textbf{stable} if the group of self-equivalences $\aut(f)$ is finite. 
By Remark~\ref{UniqueLog_rmk}, a pre-log map is stable if and only if the underlying nodal marked $J$-holomorphic map is stable. Clearly, the automorphism group of a pre-log map is a subgroup of the automorphism group of the underlying nodal marked $J$-holomorphic map.
Example~\ref{DiffAut} below illustrates some rare cases when the two groups are different. \textit{The equivalence class of a pre-log map is called a pre-log \textit{curve}. 
For every such $\Gamma$, we denote the space of $k$-marked degree $A$ pre-log $J$-holomorphic curves with dual graph $\Gamma$ and contact pattern $\mfs$ by 
\bEq{PreLog_e}
\cM^{\tn{plog}}_{g,\mfs}(X,D,A)_\Gamma.
\eEq}
If $\Gamma$ has only one vertex $v$ with $I\!=\!I_v$ then 
$$
\cM_{g,\mfs}(X,D,A)_I\equiv \cM^{\tn{plog}}_{g,\mfs}(X,D,A)_\Gamma
$$ 
is simply the space of equivalence classes of genus $g$ degree $A$ $k$-marked log $J$-holomorphic tuples with an ordering on the marked points and contact type $\mfs$.  \\

\noindent
In $g\!=0$, the forgetful map 
\bEq{ForgetLog_e}
\cM_{0,\mfs}(X,D,A)_{I} \lra \cM_{0,k}(D_I,A), \quad [u,[\ze], \Si,\mfj, \vec{z}]\lra [u, \Si,\mfj, \vec{z}]
\eEq
into the (virtually) main stratum of moduli space of $k$-marked degree $A$   $J$-holomorphic curves into $D_I$ gives an identification of two sets. That is because for every degree $d\!\in\!\Z$ holomorphic line bundle $\cL\to \P^1$, every set of distinct points $z^1,\ldots,z^k\!\in\!\P^1$, and every set of integers $m_1,\ldots,m_k$ such that $m_1+\cdots+m_k=d$, up to the action of $\C^*$, there always exists exactly one meromorphic section of $\cL$ with poles/zeros of order $m_i$ at $z^i$. In the higher genus case, however, the (virtual) normal bundle of this embedding is the direct sum of $I$ copies of dual of Hodge bundle (i.e. tangent space of $\tn{Pic}^0(\Si)$ at the trivial line bundle); see Lemma~\ref{HodgeObs_lmm}. 

\bEx{VsmoothinX_ex}
If $D$ is smooth, i.e. $N\!=\!1$, a (pre-)log map with smooth domain of depth $\eset$ is just a $J$-holomorphic map $u$ with image not into $D$,  $u^{-1}(D)\!\subset\!\vec{z}$, and 
$$
\mfs\!=\!\big(\tn{ord}_{z^a}(u,D) \big)_{a\in [k]} \in \N^{k}
$$
as in the definition of the relative moduli spaces in (\ref{RelativeOrder_e}). Thus there exists a one-to-one correspondence between the virtually main stratum of the moduli space of relative $J$-holomorphic curves of contact order $\mfs$, and the space of depth $\eset$ (pre-)log curves of the same contact pattern. Also, a depth $\{1\}$ (pre-)log $J$-holomorphic curve with smooth domain is represented by a $J$-holomorphic map $u\colon\!(\Si,\mfj)\!\lra\!(D,J|_{TD})$ and  a meromorphic  section $\ze$ of $u^*\cN_XD$ such that $\vec{z}$ includes the set of zeros and poles of $\ze$ and 
$$
\mfs\!=\!\big(\tn{ord}_{z^a}(\ze) \big)_{a\in [k]} \in \Z^{k}
$$
as in the definition of the relative moduli spaces. The definitions, however, become different if we consider maps with nodal domain.     
\eEx

\vskip.1in
\noindent
For some decorated dual graphs $\Gamma$, the expected dimension of $\cM^{\tn{plog}}_{g,\mfs}(X,D,A)_\Gamma$, calculated via (\ref{VirDimRel_e}) and the matching conditions at the nodes, could be bigger than or equal to the expected dimension of the (virtually) main stratum $\cM_{g,\mfs}(X,D,A)$ (something that we do not want to happen); see the following example.
In order for a nodal pre-log curve to be in the limit of the (virtually) main stratum, there are other global combinatorial and non-combinatorial obstructions that we are going to describe next. Of course, as in the classical case, we might get pre-log curves satisfying these conditions that do not belong to the closure of the main stratum.

\bEx{32+01inP^2_ex}
Let $X\!=\!\P^2$ with projective coordinates $[x_1,x_2,x_3]$ and $D\!=\!D_1\!\cup\!D_2$ (thus $N\!=\!2$) be a transverse union of two hyperplanes (lines). For \vskip-.1in
$$
g\!=\!0,\quad \mfs\!=\!((3,2),(0,1))\!\in\!(\N^{2})^2,\quad\tn{and}\quad A\!=\![3]\!\in\!H_2(X,\Z)\!\cong\!\Z,
$$ 
$\cM_{0,\mfs}(X,D,[3])$ is a manifold of complex dimension $4$. If $D_1\!=\!(x_1\!=\!0)$ and $D_2\!=\!(x_2\!=\!0)$, every element in $\cM_{0,\mfs}(X,D,[3])$ is equivalent to a holomorphic map of the form
\bEq{32maps_e}
[z,w]\!\lra [z^3, z^2w, a_3z^3+a_2z^2w + a_1zw^2 + a_0 w^3].
\eEq
Let $\Gamma$ be the dual graph with three vertices $v_1,v_2,v_3$, and two edges $e_1,e_2$ connecting $v_1$ to $v_3$ and $v_2$ to $v_3$, respectively. Furthermore, choose the orientations $\uvec{e}_1$ and $\uvec{e}_2$ to end at $v_3$ and assume 
$$
I_{v_1}\!=\!I_{v_2}\!=\!\eset,\quad I_{v_3}\!=\!\{1,2\},\quad s_{\uvec{e}_1}\!=\!(2,1),\quad s_{\uvec{e}_2}\!=\!(1,1),\quad A_{v_1}\!=\![2],\quad A_{v_2}\!=\![1];
$$
see Figure~\ref{Freeline_fig}. Note that $u_{v_3}$ is map of degree $0$ from a sphere with three special points, two of which are the nodes connecting $\Si_{v_3}$ to $\Si_{v_1}$ and $\Si_{v_2}$ and the other one is the first marked point $z^1$ with contact order $(3,2)$. The second marked point with contact order $(0,1)$ lies on $\Si_{v_1}$. A simple calculation shows that $\cM^{\tn{plog}}_{0,\mfs}(X,D,[3])_\Gamma$ is also a manifold of complex dimension $4$. Image of $u_2$ (blue) could be any line different from $D_1$ and $D_2$ passing through $D_{12}$, and every such $u_1$ is equivalent to a holomorphic map of the form
$$
[z,w]\!\lra [z^2, zw, a_2z^2+a_1zw + a_0w^2].
$$
\qed
\begin{figure}
\begin{pspicture}(-3,-.5)(11,2)
\psset{unit=.3cm}
\psline(0,0)(8,0)\psline(0,0)(0,8)

\psline[linecolor=blue](0,0)(8,5)
\rput(10,5.5){\small{$\tn{Im}(u_{v_2})$}}

\psarc[linecolor=red](10,0){10}{140}{180}
\psarc[linecolor=red](3,0){3}{0}{180}
\psarc[linecolor=red](8,0){2}{140}{180}
\rput(3.5,7){\small{$\tn{Im}(u_{v_1})$}}

\rput(9,9){\small{$X$}}
\rput(9,0){\small{$D_1$}}\rput(0,9){\small{$D_2$}}

\psline[linewidth=.1](20,1)(27,1)\psline[linewidth=.1](20,1)(20,8)
\psline[linewidth=.15]{->}(27,1)(24,1)\psline[linewidth=.15]{->}(20,8)(20,5)
\rput(25.5,2){\small{$\uvec{e}_2$}}\rput(21,6.5){\small{$\uvec{e}_1$}}
\pscircle*(20,1){.25}\pscircle*[linecolor=blue](27,1){.25}\pscircle*[linecolor=red](20,8){.25}
\rput(19,2){\small{$v_3$}}\rput(28,0){\small{$v_2$}}\rput(19,8){\small{$v_1$}}
\psline[linewidth=.1](20,1)(18.7,-.3)\rput(18,-1){\small{$(3,2)$}}
\psline[linewidth=.1](20,8)(21.3,9.3)\rput(22.7,9.3){\small{$(0,1)$}}

\end{pspicture}
\caption{A $2$-marked genus $0$  nodal degree $3$ pre-log map in $\P^3$ relative to two lines. The blue curve is a line. The red curve is a conic. They are connected by a ghost bubble that maps to $D_{12}$.}
\label{Freeline_fig}
\end{figure}
\eEx

\noindent
Corresponding to the decorated dual graph $\Gamma\!=\!\Gamma(\V,\E,\L)$ of a pre-log map as in Definition~\ref{PreLogMap_dfn}  and an arbitrary orientation $O\!\equiv\!\{\uvec{e}\}_{e\in \E} \subset \uvec{\E}$ on the edges, we define a homomorphism of $\Z$-modules
\bEq{DtoT_e}
\D=\D(\Gamma)\equiv \Z^\E\oplus \bigoplus_{v\in \V} \Z^{I_v}  \stackrel{\vr=\vr_{O}}{\xrightarrow{\hspace*{1.5cm}}} 
\T=\T(\Gamma)\equiv\bigoplus_{e\in \E} \Z^{I_e}
\eEq
in the following way. For every $e\!\in\!\E$, let 
\bEq{vre_e}
\vr(1_{e})\!=\!s_{\uvec{e}}\!\in\!\Z^{I_e},
\eEq
where $1_{e}$ is the generator of $\Z^e$ in $\Z^{\E}$ and $\uvec{e}$ is the chosen orientation on $e$ in $O$. In particular, $\vr(1_{e})\!=\!0$ for any $e$ with $I_e\!=\!\eset$. 
Similarly, for every $v\!\in\!\V$ and $i\!\in\! I_{v}$, let $1_{v,i}$ be the generator of the $i$-th factor in $\Z^{I_v}$, and define 
\bEq{vrv_e}
\vr(1_{v,i})\!=\!\xi_{v,i}\!\in\! \bigoplus_{e\in \E} \Z^{I_e}
\eEq 
to be the vector which has $1_{e,i}\!\in\!\Z^{I_e}\!\subset\!\Z^N$ in the $e$-th factor, if $v\!=\!v_1(\uvec{e})$ and $e$ is not a loop, it has $-1_{e;i}\!\in\!\Z^{I_e}$ in the $e$-th factor, if $v\!=\!v_2(\uvec{e})$ and $e$ is not a loop, and is zero otherwise. This is well-defined by the first equality in (\ref{Ieunion_e}). Let 
\bEq{cLGamma_e}
\Lambda\!=\!\Lambda(\Gamma)\!=\!\tn{image}(\vr),\quad \K\!=\!\K(\Gamma)=\tn{Ker}(\vr)\quad\tn{and}\quad  
\CK\!=\!\CK(\Gamma)=\T/\Lambda=\tn{coker}(\vr).
\eEq
By Definition~\ref{PreLogMap_dfn}.\ref{MatchingOrders_l}, the $\Z$-modules $\Lambda$, $\K$, and $\CK$ are independent of choice of the orientation $O$ on $\E$ and are invariants of the decorated graph $\Gamma$. In particular,
\bEq{kernel_e}
\K=\big\{\big((\la_e)_{e\in \E},(s_v)_{v\in \V}\big)\!\in\!\ \Z^\E\oplus \bigoplus_{v\in \V} \Z^{I_v} \colon s_{v}\!-\!s_{v'}\!=\!\la_es_{\uvec{e}}\quad\forall~v,v'\!\in\!\V,~\uvec{e}\!\in\!\uvec{\E}_{v',v}\big\}.
\eEq
In this equation, via the first identity in (\ref{Ieunion_e}) and the inclusion $\Z^{I_v}\!\cong\!\Z^{I_v}\!\times\!\{0\}^{I_e-I_v}\!\subset\!\Z^{I_e}$, we think of $s_{v}$ as a vector also in $ \Z^{I_e}$, for all $e\!\in\!\E_v$.
For any field of characteristic zero $F$, let 
\bEq{FVersion_e}
\D_F=\D\otimes_\Z\! F,~~~ \T_F=\T\otimes_\Z \!F,~~~ \La_F=\La\otimes_\Z\!F,~~~\K_F=\K\otimes_\Z\!F,~~~\tn{and}~~~\CK_F=\CK\otimes_\Z\!F
\eEq
be the corresponding $F$-vector spaces and $\vr_F\colon\!\D_F\!\lra\!\T_F$ be the corresponding $F$-linear map.
Via the exponentiation map, let 
$$
\tn{exp}(\Lambda_\C) \subset  \prod_{e\in \E} (\C^*)^{I_e}
$$
be the subgroup corresponding to the sub-Lie algebra $\La_\C\!\subset\!\T_\C$, and denote the quotient group by 
$$
\mc{G}\!=\!\mc{G}(\Lambda)\!=\!\Big(  \prod_{e\in \E} (\C^*)^{I_e}\Big)/\tn{exp}(\Lambda_\C) =\frac{  \prod_{e\in \E} (\C^*)^{I_e}}{\tn{exp}(\rho_\C)\Big( (\C^*)^\E \times \prod_{v\in \V} (\C^*)^{I_v}\Big)}
$$ 

\noindent
In the following, we will construct a map
\bEq{PLtoG_e}
\cM^{\tn{plog}}_{g,\mfs}(X,D,A)_\Gamma \stackrel{\ob_\Gamma}{\lra} \mc{G}(\Gamma)
\eEq
that will be used in the definition of log moduli spaces.\\

\noindent
Given a pre-log map $f\equiv \big(f_{v}\!\equiv\! (u_v,[\ze_v],C_v)\big)_{v\in \V}$ as in Definition~\ref{PreLogMap_dfn}, fix an arbitrary set of representatives  
\bEq{zeRep_e}
\ze_v=(\ze_{v,i})_{i\in I_v}\in \Om_{\tn{mero}}(\Si_v, u_v^*\cN_{X}D_{I_v})\qquad \forall~v\!\in\!\V.
\eEq
For each $v\!\in\!\V$ and $\uvec{e}\!\in\!\uvec{\E}_v$, let $z_{\uvec{e}}$ be an arbitrary holomorphic coordinate in a sufficiently small disk $\De_{\uvec{e}}$ around the nodal point $(z_{\uvec{e}}\!=\!0)\!=\!q_{\uvec{e}}\!\in\!\Si_v$.  By (\ref{Ordx_e2}), for every  $v\!\in\!\V$, $\uvec{e}\!\in\!\uvec{\E}_v$, and $i\!\in\!I_v$, in a local holomorphic trivialization 
$$
u^*\cN_XD_i|_{\De_{\uvec{e}}}\approx \cN_XD_i|_{u(q_e)}\!\times\!\De_{\uvec{e}},
$$ 
we have 
\bEq{LocalCoord_e}
\ze_{v,i}(z_{\uvec{e}})= z_{\uvec{e}}^{s_{\uvec{e},i}} \wt{\ze}_{v,i}(z_{\uvec{e}})
\eEq
such that 
$$
0\!\neq\!  \wt{\ze}_{v,i}(0) \! \equiv\!\eta_{\uvec{e},i}\!\in\!\cN_{X}D_i|_{u(q_e)}
$$ 
is independent of the choice of  the trivialization. Similarly, by \cite[(6.1)]{FZ}, 
for every  $v\!\in\!\V$, $\uvec{e}\!\in\!\uvec{\E}_v$, and $i\!\in\!I_e\!-\!I_v$, the map $u_v$ 
has a well-defined $s_{\uvec{e},i}$-th derivative 
\bEq{normalDerivative}
\eta_{\uvec{e},i}\!\in\!\cN_{X}D_i|_{u(q_e)}
\eEq
(with respect to the coordinate $z_{\uvec{e}}$) in the normal direction to $D_i$ at the nodal marked point $q_{\uvec{e}}$.\\

\noindent
With the choice of orientation $O\!\equiv\!\{\uvec{e}\}_{e\in \E}\! \subset\! \uvec{\E}$ on the edges as before,
since $\eta_{\uvec{e},i}\!\neq\!0$ for all $\uvec{e}\!\in\!\uvec{\E}$ and $i\!\in\!I_e$,
the tuples
\bEq{etae_e}
\eta_{e}=\big(\eta_{\uvec{e},i}/\eta_{\scz\ucev{e},i}\big)_{i\in I_e}\!\in\!(\C^*)^{I_e}\qquad \forall~\uvec{e}\in O
\eEq
give rise to an element 
\bEq{Totaleta_e}
\eta\!\equiv \!(\eta_e)_{e\in \E}\!\in\! \prod_{e\in \E} (\C^*)^{I_e}.
\eEq
The action of the subgroup $\tn{exp}(\Lambda_\C)$ on $\eta$ corresponds to rescalings of (\ref{zeRep_e}) and change of coordinates in (\ref{LocalCoord_e}); i.e the class $\ob_\Gamma(f)\!=\![\eta]$ of $\eta$ in 
$$
\mc{G}\!=\! \Big(\prod_{e\in \E} {\C^*}^{I_e}\Big)/\tn{exp}(\Lambda_\C)
$$ 
is independent of the choice of representatives in (\ref{zeRep_e}) and local coordinates in (\ref{LocalCoord_e}). If $f$ and $f'$ are equivalent with respect to a reparametrization $h\colon\!\Si'\!\lra\!\Si$ as in (\ref{NodalEquivalence_e}), the associated group elements $\eta$ and $\eta'$, respectively, would be the same with respect to any $h$-symmetric choice of holomorphic coordinates $\{z_{\uvec{e}}\}_{\uvec{e}\in\uvec{\E}}$.
Therefore, 
\bEq{[eta]_e}
\ob_\Gamma([f])\!:=\![\eta]\!\in\!\mc{G}
\eEq
is well-defined. By definition, $\ob_\Gamma([f])=1$ if and only if  there exists a choice of representatives  $\{\ze_{v,i}\}_{i\in I_v, v \in \V}$ and local coordinates 
$\{z_{\uvec{e}}\}_{\uvec{e}\in \E}$ such that 
$$
\eta_{\uvec{e}}=\eta_{\scz\ucev{e}}\qquad \forall~e\!\in\! \E.
$$

\bDf{LogMap_dfn}
Let $D\!= \!\bigcup_{i\in [N]} D_i \!\subset\! X$ be an SNC symplectic divisor and $(\om,J)\!\in\!\cJ(X,D)$. A \textbf{log $J$-holomorphic map} is a pre-log $J$-holomorphic map  $f$ with the decorated dual graph $\Gamma$ such that 
\bEn
\item\label{Tropical_l} there exist functions 
$$
s \colon\! \V\!\lra\!\R^N, \quad v\!\lra\!s_v,\qquad\tn{and}\qquad \la \colon\! \E\!\lra\!\R_+, \quad e\!\lra\!\la_e,
$$
such that 
\bEnalph
\item\label{Stratum_l} $s_v\!\in\!\R_{+}^{I_v}\!\times\!\{0\}^{[N]-I_v}$ for all $v\!\in\!\V$,
\item\label{Direction_l} $s_{v_{2}(\uvec{e})}\!-\!s_{v_1(\uvec{e})}\!=\!\la_{e} s_{\uvec{e}}$  for every $\uvec{e}\!\in\!\uvec{\E}$;
\eEnalph 
\item\label{GObs_e} and $\ob_\Gamma(f)\!=\!1\!\in\!\mc{G}(\Gamma)$.
\eEn
\eDf
\noindent 
Condition~\ref{Tropical_l}\ref{Direction_l} is well-defined because of Definition~\ref{PreLogMap_dfn}.\ref{MatchingOrders_l}. 
If \ref{GObs_e} holds, we say that the pre-log map $f$ is $\mc{G}$-\textbf{unobstructed}. Condition~\ref{GObs_e} is independent of the choice of orientation $O$ on $\E$ used to define $\tn{ob}_\Gamma$.

\bRm{TropicalCondition_rmk}
A nodal map in the relative compactification (when $D$ is smooth) with image in an expanded degeneration $X[m]$ comes with a partial ordering of the smooth components of the domain, such that the components mapped into $X$ have order $0$ and those mapped into the $r$-th copy of $\P_{X}D$ are of order $r$; see Section~\ref{RelComp_ss}. In the compactification process, a component sinking faster into $D$ results in a component with higher order.
From our perspective, the vector-valued function $s\colon\!\V\!\lra\!\R^N$ in Condition~\ref{Tropical_l} is a generalization of this partial ordering to the SNC case with $\R^N$ instead of $\Z$; see Lemma~\ref{Partial-Order_lmm}.
From the tropical perspective of \cite[Dfn~2.5.3]{ACGS}, Condition~\ref{Tropical_l} is equal to the existence of a tropical map from a tropical curve associated to $\Gamma$ into $\R_{\geq 0}^N$.
 This condition puts a big restriction on the set of contact vectors $s_{\uvec{e}}$. For example if $I_v,I_{v'}\!=\!\eset$, then for any other $v''\!\in\!\V$ and oriented edges $\uvec{e}\!\in\!\uvec{\E}_{v,v''}$ and $\uvec{e}'\!\in\!\uvec{\E}_{v',v''}$, the contact vectors $s_{\uvec{e}}$ and $s_{\uvec{e}'}$ should be positively proportional. Condition (2) has no explicit equivalent in \cite{AC,GS,BP2,I} but it is related to the slope condition at each node in \cite{I}.
 \eRm
 
 \bRm{Various-s_rmk}
The discussion above includes $\R^N$-valued functions, all of them denoted by $s$, on the set of vertices, oriented edges, and legs of a decorated dual graph $\Gamma$ that play different roles and should not be confused. 
The contact orders $\mfs\!=\!(s_1,\ldots,s_k)$ at the legs (marked points) are fixed for a moduli space (they are independent of $\Gamma$) and define a function $s\colon\!\L\!\lra\!\Z^N$. The contact orders $(s_{\uvec{e}})_{\uvec{e}\in \uvec{\E}}$ at nodal points define a function $s\colon\!\uvec{\E}\!\lra\!\Z^N$ and are part of the decoration of $\Gamma$.
Finally the function $s \colon\! \V\!\lra\!\R^N$ (and $\la \colon\! \E\!\lra\!\R_+$) is not part of the defining data of a log map. We only require the latter to exist in order for a pre-log map to define a log map.
 \eRm
 
 \bEx{Backto32+01inP^2_ex}
 Example~\ref{32+01inP^2_ex} does not satisfy Definition~\ref{LogMap_dfn}.\ref{Tropical_l}.
 Since $I_{v_1}\!=\!I_{v_2}\!=\!\eset$, we should have $s_{v_1}\!=\!s_{v_2}\!=(0,0)$. Then
 Condition~\ref{Tropical_l}\ref{Direction_l} requires $s_{\uvec{e}_1}\!=\!(2,1)$ and $s_{\uvec{e}_1}\!=\!(1,1)$ to be positive multiples of $s_{v_3}$, which is impossible. A straightforward calculation shows that the line component $u_{v_2}$ in any limit of (\ref{32maps_e}) with a component $u_{v_1}$ as in Figure~\ref{Freeline_fig} should lie in $D_1$. Then the function $s\colon\!\V\!\lra\!\R^2$ given by $s_{v_1}\!=\!(0,0)$, $s_{v_2}\!=\!(1,0)$, and $s_{v_3}\!=\!(2,1)$ satisfies Definition~\ref{LogMap_dfn}.\ref{Tropical_l}.
 \eEx
 
 \noindent
 The following definition lists the combinatorial properties of an admissible decorated dual graph.
 
 \bDf{DecoratedGamma_dfn}
For a fixed SNC symplectic divisor $D\!=\! \bigcup_{i\in [N]} D_i$ in $X$, given $g,k\!\in\!\N$, $A\!\in\!H_2(X,\Z)$, and $\mfs\!\in\!(\Z^N)^k$, we denote by $\tn{DG}(g,\mfs,A)$ to be the set of (stable) connected dual graphs $\Gamma\!=\!\Gamma(\V,\E,\L)$ with $k$ legs and
\bEnalph
\item a genus decoration of total genus $g$, 
\item a degree decoration of total degree $A$, 
\item an ordering $a\colon\!\L\!\lra\!\{1,\ldots,k\}$,
\item set decorations $I\colon \V,\E\lra \cP(N)$ satisfying $I_e\!=\!I_v\!\cup\!I_{v'}$ for all $v,v'\!\in\!\V$ and $e\!\in\!\E_{v,v'}$,
\item\label{vecDec_it} and a vector decoration on the set of oriented edges $\uvec{\E}$, $\uvec{e}\!\to\!s_{\uvec{e}}\!\in\!\Z^{I_e}\!\subset\!\Z^N$, satisfying 
$$
s_{\uvec{e}}\!+\!s_{\scz\ucev{e}}\!=\!0\qquad \forall~\uvec{e}\!\in\!\uvec{\E},
$$
\eEnalph
such that Condition \ref{Tropical_l} of Definition~\ref{LogMap_dfn} holds and 
\bEq{DegreeTos_e}
(A_v\cdot D_i)_{i\in [N]} = \sum_{\uvec{e}\in \uvec{\E}_v} s_{\uvec{e}}+\sum_{\substack{l\in \L\\ v_l=v}} s_{l}\qquad \forall~v\!\in\!\V.
\eEq 
\eDf

\noindent
$\tn{DG}(g,\mfs,A)$ is the set of possible combinatorial  types of stable  connected genus $g$ $k$-marked degree $A$ log curves of contact type $\mfs$. 
Note that the defining conditions of $\tn{DG}(g,\mfs,A)$ do not capture Definition~\ref{LogMap_dfn}.\ref{GObs_e}; the latter is a non-combinatorial condition. 
Example~\ref{2dinP3toP2P2_ex} below illustrates a legitimate $\Gamma$ such that the moduli space of pre-log curves of type $\Gamma$ has an expected dimension larger than the expected dimension of the (virtually) main stratum. Then, imposing Condition~\ref{GObs_e} of Definition~\ref{LogMap_dfn} would reduce the dimension to less than the expected dimension of the (virtually) main stratum. \\

\noindent
For every $\Gamma\!\in\!\tn{DG}(g,\mfs,A)$ define 
\bEq{GammaLogStratum_e}
\cM_{g,\mfs}(X,D,A)_\Gamma=\tn{ob}_\Gamma^{-1}(1)\!\subset\!\cM^{\tn{plog}}_{g,\mfs}(X,D,A)_\Gamma
\eEq
to be the stratum of log $J$-holomorphic curves of type $\Gamma$. \textit{We then define the moduli space of genus $g$ degree $A$ \textbf{stable nodal log $J$-holomorphic curves} of contact type $\mfs$ to be the union}
\bEq{LogCpt_e}
\ov\cM^{\log}_{g,\mfs}(X,D,A)\equiv \bigcup_{\Gamma\in \tn{DG}(g,\mfs,A)} \cM_{g,\mfs}(X,D,A)_\Gamma,\\
\eEq
\begin{figure}
\begin{pspicture}(4,1.3)(11,4.5)
\psset{unit=.3cm}

\pscircle*(35,10){.25}\rput(35.1,9){\small$v_1$}\rput(34.9,8){\tiny$I_{v_1}\!=\!\{1\}$}\rput(34.9,7){\tiny$g_{v_1}\!=\!\frac{d(d-1)}{2}$}\rput(34.9,5.7){\tiny$A_{v_1}\!=\![d]$}
\pscircle*(45,10){.25}\rput(45.1,9){\small$v_2$}\rput(45.1,8){\tiny$I_{v_2}\!=\!\{2\}$}\rput(45.1,7){\tiny$g_{v_2}\!=\!\frac{d(d-1)}{2}$}\rput(45.1,5.7){\tiny$A_{v_2}\!=\![d]$}
\psarc(40,6.55){6}{35}{145}\psarc{<-}(40,6.55){6}{90}{145}
\psarc(40,4){7.75}{50}{130}\psarc{<-}(40,4){7.75}{90}{130}
\rput(40,10.3){$\vdots$}\rput(41.8,10){\tiny$d$~edges}
\psarc(40,13.45){6}{215}{325}\psarc{->}(40,13.45){6}{215}{270}
\psarc(40,16){7.75}{230}{310}\psarc{->}(40,16){7.75}{230}{270}

\rput(40,13.5){\tiny$s_{\uvec{e}}\!=\!(-1,1)$}

\psline(35,10)(32,12)
\psline(35,10)(32,11)
\rput(32.2,10.3){$\vdots$}
\psline(35,10)(32,8)
\psline(35,10)(32,9)

\psline(45,10)(48,12)
\psline(45,10)(48,11)
\rput(47.8,10.3){$\vdots$}
\psline(45,10)(48,8)
\psline(45,10)(48,9)

\rput(27,11){\tiny $2d$ points of}
\rput(28,10){\tiny contact order $(1,0)$}
\rput(51,11){\tiny $2d$ points of}
\rput(52.5,10){\tiny contact order $(0,1)$}

\end{pspicture}
\caption{A decorated graph in $\tn{DG}(g\!=\!(d-1)^2,\mfs,A\!=\![2d])$, corresponding to two generic degree $d$ curves in $D_1$ and $D_2$ intersecting at $d$ points along $D_{12}$.}
\label{dd-decomp_fg}
\end{figure}
\bEx{2dinP3toP2P2_ex}
Let 
$$
\aligned
&X\!=\!\P^3,\quad D_1\cup D_2\!=\!\P^2\cup \P^2,\quad A\!=\![2d]\!\in\!H_2(X,\Z)\!\cong\! \Z,\quad g\!=\!(d-1)^2, \\
&\mfs=\big((1,0),\ldots,(1,0),(0,1),\ldots,(0,1)\big)\!\in\!\big(\Z^2\big)^{4d},
\endaligned
$$
and $\Gamma\!\in\!\tn{DG}(g,\mfs,A)$ be the decorated dual graph illustrated in Figure~\ref{dd-decomp_fg}. Note that the function $s\colon\!\V\!\lra\!\R^2$ given by $s_{v_1}\!=\!(1,0)$ and $s_{v_2}\!=\!(0,1)$ satisfies Definition~\ref{LogMap_dfn}.\ref{Tropical_l}. Every element of $\cM^{\tn{plog}}_{g,\mfs}(X,D,A)_\Gamma$ is supported on two generic degree $d$ plane curves in $D_1$ and $D_2$ intersecting at $d$ points along $D_{12}$. By (\ref{VirDimRel_e}) and Definition~\ref{PreLogMap_dfn}.\ref{MatchingPoints_l}, the expected $\C$-dimension of $\cM_{g,\mfs}(X,D,A)$ and $\cM^{\tn{plog}}_{g,\mfs}(X,D,A)_\Gamma$ are $8d$ and $9d\!-\!2$, respectively.\\

\noindent
Orient each edge such that $v_1(\uvec{e}_i)\!=\!v_1$ for all $i\!=\!1,\ldots, d$.  Then $\Lambda\!=\!\vr(\D)$ in (\ref{cLGamma_e}) is generated by the vectors $s_{\uvec{e}_1},\ldots,s_{\uvec{e}_d}$, $\xi_{v_1}\!=\!\xi_{v_1,1}$, and $\xi_{v_2}\!=\!\xi_{v_2,2}$, such that the only relation is 
$$
\xi_{v_1}+\xi_{v_2}+(s_{\uvec{e}_1}+\ldots+s_{\uvec{e}_d})=0.
$$
We conclude that the obstruction group $\mc{G}(\Gamma)$ is complex $(d\!-\!1)$-dimensional. Therefore, the subset of log curves
$$
\cM_{g,\mfs}(X,D,A)_\Gamma\!\subset \cM^{\tn{plog}}_{g,\mfs}(X,D,A)_\Gamma
$$
is of the expected $\C$-dimension $(9d\!-\!2)\!-\!(d-1)\!=\!8d\!-\!1\!<\!8d$.\qed
\eEx

\bRm{UniqueLog_rmk2}
By Remark~\ref{UniqueLog_rmk}, for every $k$-marked stable nodal curve $f$ in $\ov\cM_{g,k}(X,A)$ with dual graph $\Gamma$, fixing $\mfs\!\in\!(\Z^N)^k$ and the vector decoration $\{s_{\uvec{e}}\}_{\uvec{e}\in \uvec{\E}}$  in Definition~\ref{DecoratedGamma_dfn}.\ref{vecDec_it}, there exists at most one log curve $f_{\log}\!\in\!\ov\cM^{\log}_{g,\mfs}(X,D,A)$ with orders $s_i$ at $z^i$ and $s_{\uvec{e}}$ at $q_{\uvec{e}}$ lifting $f$. Furthermore, $f_{\log}$ is stable if and only if $f$ is stable.\eRm

\bLm{UniqueLog_lm}
Given $f\!\in\!\ov\cM_{g,k}(X,A)$ with the dual graph $\Gamma$ and $\mfs\!\in\!(\Z^N)^k$, the set of possible vector decorations $\{s_{\uvec{e}}\}_{\uvec{e}\in \uvec{\E}}$ as in Definition~\ref{DecoratedGamma_dfn} satisfying (\ref{DegreeTos_e}), and thus the set of possible log lifts of $f$ is finite. 
\eLm

\bPf
Since $s_{\uvec{e},i}\!=\!-s_{\scz\ucev{e},i}$ for all $e\!\in\! \E$ and $i\! \in\! I_e$, it is sufficient to show that the set of possible values for $\{s_{\uvec{e},i}\}_{\uvec{e}\in \uvec{\E}, i\in I_e}$ is bounded from above. Fix $i\in [N]$. If $v\in \V$ and $e\in \E_{v}$ such that $i\!\in\! I_e-I_v$, then $s_{\uvec{e},i}$ and $s_{\scz\ucev{e},i}$  are uniquely determined by the tangency order of $u_v$ with $D_i$. Therefore, we can restrict to the subset  $\V_i\subset \V$ of all vertices $v$ such that $i \!\in\! I_v$, and the edges between them. Given a decoration $\{s_{\uvec{e}}\}_{\uvec{e}\in \uvec{\E}}$ as in Definition~\ref{DecoratedGamma_dfn} satisfying (\ref{DegreeTos_e}), let $\uvec{\E}_i$ be the subset of oriented edges $\uvec{e}$ such that $e \in \E_{v,v'}$ for some $v,v'\in \V_i$ and $s_{\uvec{e},i}> 0$.  Let $\Gamma_i$ be the oriented graph made of $\V_i$ and the oriented edges in $\uvec{\E}_i$. By Condition (1b) in Definition~\ref{LogMap_dfn}, $\Gamma_i$ does not have any oriented loop. Therefore, $\uvec{\E}_i$ defines a partial order on $\V_i$. Let $v\in \V_i$ be a maximal vertex. There is no oriented edge in $\uvec{\E}_i$ pointing toward $v$. Therefore, for every $\uvec{e}\in \uvec{\E}_v$ either $s_{\uvec{e},i}=0$ or $s_{\uvec{e},i}>0$. The identity
$$
A_v\cdot D_i= \sum_{\uvec{e}\in \uvec{\E}_v} s_{\uvec{e},i}+\sum_{\substack{l\in \L\\ v_l=v}} s_{l,i}
$$
puts an upper bound on  $\{s_{\uvec{e},i}\}_{\uvec{e}\in \uvec{\E}_v}$. Moving down in the partial order on $\V_i$ we get upper-bounds on other $s_{\uvec{e},i}$.

\ePf

\bLm{g0Embedding}
For every genus $0$ $k$-marked stable nodal map $f$ in $\ov\cM_{0,k}(X,A)$ with dual graph $\Gamma$ and a fixed $\mfs$, there exists at most one vector decoration $\{s_{\uvec{e}}\}_{\uvec{e}\in \uvec{\E}}$ as in Definition~\ref{DecoratedGamma_dfn}.\ref{vecDec_it} satisfying (\ref{DegreeTos_e}). 
In particular, the forgetful map 
$$
\ov\cM^{\log}_{0,\mfs}(X,D,A) \lra \ov\cM_{0,k}(X,A)
$$
is an embedding (of sets).
\eLm

\bPf
Assume there are two different decorations $\{s_{\uvec{e}}\}_{\uvec{e}\in \uvec{\E}}$ and $\{s'_{\uvec{e}}\}_{\uvec{e}\in \uvec{\E}}$ as in Definition~\ref{DecoratedGamma_dfn} satisfying (\ref{DegreeTos_e}). There is some $i\!\in\! [N]$ such that $\{s_{\uvec{e},i}\}_{\uvec{e}\in \uvec{\E}}$ and $\{s'_{\uvec{e},i}\}_{\uvec{e}\in \uvec{\E}}$ are different. 
Since $g\!=\!0$, $\Gamma$ is a tree and the subset of edges $\Om\!\subset\!\E$ where $s_{\uvec{e},i}\!\neq\!s'_{\uvec{e},i}$ determines a sub-tree of that. In particular, there exists a vertex $v\!\in\!\V$ that is connected to only one edge $e'\!\in\!\Om$. Orient $e'$ so that $v$ is the starting point. Then, by (\ref{DegreeTos_e}),
$$
A_v\cdot D_i= s_{\uvec{e}',i}+ \sum_{\uvec{e}\in \uvec{\E}_v-\{\uvec{e}'\} } s_{\uvec{e},i}+\sum_{\substack{l\in \L\\ v_l=v}} s_{l,i}\neq  s'_{\uvec{e}',i}+ \sum_{\uvec{e}\in \uvec{\E}_v-\{\uvec{e}'\} } s'_{\uvec{e},i}+\sum_{\substack{l\in \L\\ v_l=v}} s'_{l,i}=A_v\cdot D_i\,;
$$
that is a contradiction.
\ePf

\noindent
Example~\ref{Twolifts} below describes a situation where $f$ has different lifts but the automorphism groups of $f$ and its lifts are the same.
Example~\ref{DiffAut} describes a situation where $f$ has different lifts and some of them have smaller automorphism groups.

\bEx{Twolifts}
Let $X\!=\!\P^2$, $D\!=\!\P^1$ be a hyperplane (line), and $p_1,p_2,p_3,p_4$ be four distinct points in $D$. Let $u_{v_1}\colon\!\Si_{v_1}\!\stackrel{\cong}{\lra}\!D$ be a degree one map and $\ze_{v_1}$ be a meromorphic section of $u_{v_1}^*\cN_XD$ with two poles of order $1$ and $2$ at $q_{\uvec{e}_1}=u_{v_1}^{-1}(p_1)$ and $q_{\uvec{e}_2}=u_{v_1}^{-1}(p_2)$, respectively, and a zero of order $4$ at $z^1=u_{v_1}^{-1}(p_3)$. Similarly, let $u_{v_2}\colon\!\Si_{v_2}\!\stackrel{\approx}{\lra}\!D$ be a degree one map and $\ze_{v_2}$ be a meromorphic section of $u_{v_2}^*\cN_XD$ with two zeros of order $1$ and $2$ at $q_{\scz\ucev{e}_1}=u_{v_2}^{-1}(p_1)$ and $q_{\scz\ucev{e}_2}=u_{v_2}^{-1}(p_2)$, respectively, and a pole of order $2$ at $q_{\uvec{e}_3}=u_{v_2}^{-1}(p_4)$. Finally, let $u_{v_3}\colon\!\Si_{v_3}\!\lra \!X$ be a smooth conic with a tangency of order $2$ with $D$ at $q_{\scz\ucev{e}_3}=u_{v_3}^{-1}(p_4)$.
The tuple 
$$
u_{\log}\equiv (u_{v_3}, (u_{v_2},\ze_{v_2}),(u_{v_1},\ze_{v_1}))
$$
with the nodal $1$-marked domain
$$
(\Si,z^1)\!=\! (\Si_{v_1},z^1,q_{\uvec{e}_1},q_{\uvec{e}_2}) \sqcup (\Si_{v_2}, q_{\scz\ucev{e}_1},q_{\scz\ucev{e}_2},q_{\uvec{e}_3} ) \sqcup (\Si_{v_3}, q_{\scz\ucev{e}_3} ) / q_{\uvec{e}}\sim q_{\scz\ucev{e}} ~~~\forall~e\!\in\!\{e_1,e_2,e_3\}
$$
defines an element of $\ov\cM^{\log}_{1,(4)}(X,D,[4])$. Let $u'_{\log}$ be a similar tuple with the roles of $p_1$ and $p_2$ reversed, i.e. $u\!=\!(u_{v_1},u_{v_2},u_{v_3})$ remains the same but $\ze_{v_1}$ and $\ze_{v_2}$ exchange their orders at the pre-images of $p_1$ and $p_2$. Therefore, $[u_{\log},\Si,z^1]$ and $[u'_{\log},\Si,z^1]$ are different lifts of the same $1$-marked stable curve $[u,\Si,z^1]$ in $\ov\cM_{1,1}(X,[4])$. Note that $e_1$ and $e_2$ form a loop in $\Gamma$. In this example, 
the two vector decorations corresponding to $(u_{\log},\Si,z^1)$ and $(u'_{\log},\Si,z^1)$ yield isomorphic decorated dual graphs $\Gamma$. In other words, the forgetful map 
$$
\cM_{1,(4)}(X,D,[4])_\Gamma\!\lra\!\ov\cM_{1,1}(X,[4])
$$
is a double-covering of its image.\qed
\eEx

\bEx{DiffAut} 
Assume $u\colon\!\Si\!\lra\!D\!\subset\!X$ is a stable map, where $\Si$ is the genus $1$ nodal curve made of two copies of $\P^1$, say $\P^1_1$ and $\P^1_2$, attached at $0$ and $\infty$, and $u_i\!=\!u|_{\P^1_i}\colon\P^1_i\lra D$, for $i\!=\!1,2$, is a double-covering of some rational curve $C_i\!\subset\!D$, with $u_i(z^{-1})\!=\! u_i(z)$; i.e. 
$$
u_{i}(0)\!=\!u_i(\infty)\!=\!x\!\in\! C_1\!\cap\!C_2\!\subset \!D.
$$ 
Furthermore, assume $\cN_XD|_{C_1}\!=\!\cO(2)$ and $\cN_XD|_{C_2}\!=\!\cO(-2)$. The automorphism group of the stable map $f\!=\!(u,\Si)$ is $\Z_2$. Since $u_1^*\cN_XD\!=\!\cO(4)$ and $u_2^*\cN_XD\!=\!\cO(-4)$, there are $2$ possible ways to lift $f$ to a log map $f_{\log}\!\in\!\ov\cM^{\log}_{1,\eset}(X,D,2(C_1\!+\!C_2))$. The holomorphic section $\ze_1$ of $u_1^*\cN_XD$ can be chosen to have zeros of orders $(3,1)$, $(2,2)$, or $(1,3)$ at $(0,\infty)$. In the middle case, the automorphism group of $f_{\log}$ is $\Z_2$.  In the remaining two cases, the two lifts are equivalent with respect to the reparametrization map 
$$
h\colon \Si\lra\Si, \quad h|_{\P^1_i}(z)\!=\!z^{-1},\quad i\!=\!1,2,
$$ 
and their equivalence class defines a single element of $\ov\cM^{\log}_{1,\eset}(X,D,2(C_1\!+\!C_2))$ with the trivial  automorphism group.
\eEx

\noindent
In Section~\ref{Cpt_s}, for $J$ as in the statement of Theorem~\ref{Compactness_th}, we will lift the Gromov convergence topology to a compact sequential convergence topology on (\ref{LogCpt_e}) such that the forgetful map (\ref{FogetLog_e}) is a continuous local embedding. It follows that the lifted topology is also metrizable. If $g\!>\!0$, globally, (\ref{FogetLog_e}) behaves like an immersion.
If $\mfs\!\in\!(\N^N)^k$, by Lemma~\ref{ExpectedGamma_lmm} bellow, $\ov\cM^{\log}_{g,\mfs}(X,D,A)$ is a compact space of the expected real dimension 
\bEq{ExpectedDim_e}
2\Big(c_1^{TX(-\log D)}(A) + (\dim_\C X-3)(1-g) + k\Big).
\eEq
In subsequent papers we will construct Kuranishi-type charts of dimension (\ref{ExpectedDim_e}) around every point of $\ov\cM^{\log}_{g,\mfs}(X,D,A)$.\\

\noindent
The following example describes the log compactification of the moduli space of lines in $\P^2$ relative to a transverse union of two hyperplanes (lines).
The same example is studied in \cite{BP3}, where Parker uses tropical geometry to describe Ionel's compactification in \cite{I} and compare it with his construction.

\bEx{2LinesP2_ex}
Let $X\!=\!\P^2$ with projective coordinates $[x_1,x_2,x_3]$, $D_1\!=\!(x_1\!=\!0)$, $D_2\!=\!(x_2\!=\!0)$, $D\!=\!D_1\cup D_2$, $A\!=\![1]\!\in\!H_2(\P^2,\Z)\cong \Z$, and $\mfs\!=\!((1,0),(0,1))$. Then, as we show below, the moduli space
\bEq{CaseStudy_e}
\ov\cM^{\log}_{0,\mfs}(X,D,[1])
\eEq
can be identified\footnote{the identification is homeomorphic with respect to the topology that we describe in Section~\ref{Cpt_s}.} with $B_{\pt_1,\pt_2} \P^2_{\tn{dual}}$ (two point blowup of $\P^2$), where $\P^2_{\tn{dual}}$ is the dual space of lines in $X\!=\!\P^2$, $\pt_1$ is the point corresponding to the line $D_1$, and $\pt_2$ is the point corresponding to the line $D_2$. 
Let $E_1$ and $E_2$ be the exceptional curves of $B_{\pt_1,\pt_2} \P^2_{\tn{dual}}$ and $L$ be the proper transform of the line connecting $\pt_1$ and $\pt_2$. Any line in $X$ not passing through $D_{12}$ intersects $D_1$ and $D_2$ at two disjoint points $z^1$ and $z^2$, respectively. 
By (\ref{Ordx_e1}),
$$
\ord(z^1)=(1,0)\qquad\tn{and}\qquad \ord(z^2)=(0,1).
$$
This gives an identification of 
$$
\cM_{0,\mfs}(X,D,[1]) \subset \ov\cM^{\log}_{0,\mfs}(X,D,[1])
$$
with $B_{\pt_1,\pt_2} \P^2\!-\! (E_1\cup E_2 \cup L)$. Every other log map $(u,[\ze])$ with smooth domain in (\ref{CaseStudy_e}) is either of depth $\{1\}$ or of depth $\{2\}$ with two marked points $z^1$ and $z^2$ of the corresponding orders. Those of depth $\{1\}$ are given by an isomorphism $u\colon \P^1\!\stackrel{\cong}{\lra}\!D_1$ and a holomorphic section $\ze$ of $\cN_XD_1\!\cong\!\cO_{\P^1}(1)$, such that $\ze$ has a simple zero at the marked point $z^1$ and $z^1\!\neq\!z^2\!=\!u^{-1}(D_2)$. Such $[\ze]$ is uniquely determined by $u(z^1)\!\in\!D_1\!\cong\! \P^1$. Therefore, via the identification 
$$
E_1\cong \P\big( H_0(\cN_X D_1)\big)\cong \P^1,
$$
such maps correspond to $E_1\!-\! \{E_1\cdot L\}\!\cong\!\C$. Similarly, the maps of depth $\{2\}$ with smooth domain correspond to $E_2\!-\! \{E_2\cdot L\}\!\cong \!\C$.
For other log maps $f$ in (\ref{CaseStudy_e}), $z^1$ and $z^2$ are mapped to  the point $D_{12}$ and thus live on a ``ghost bubble" $u_2\colon\!\P^1\!\lra\!X$, with $\tn{Im}(u_2)\!\equiv\! D_{12}$. This ghost bubble and the non-trivial map $u_1\colon\!\P^1\!\lra\!X$ are attached to each other at nodal points $z^3\in \tn{Dom}(u_2)$ and $z'\in \tn{Dom}(u_1)$. By definition, the meromorphic section $\ze\!=\!(\ze_1,\ze_2)$ defining the log map $(u_2,[\ze]\equiv([\ze_1],[\ze_2]))$ is a meromorphic section of the trivial bundle $u_2^*\cN_{X}D_{12}\!\cong\!\P^1\!\times\! \C^2$, such that 
$$
\ord_{z^1}(\ze)=(1,0)\qquad \tn{and}\quad \ord_{z^2}(\ze)=(0,1).
$$ 
Since $u_2^*\cN_{X}D_{12}$ is trivial, we should have $\ord_{z^3}(\ze)\!=\!(-1,-1)$ and these restrictions specify a unique $(\C^*)^2$-class $[\ze]$.
There are thus three possibilities for $f$:
\bEn
\item $u_1$ is of depth $\eset$: in this case, by Definition~\ref{PreLogMap_dfn}.\ref{MatchingOrders_l}, $u_1$ specifies an element of $\cM_{0,((1,1))}(X,D,[1])$ and we get an identification of such curves $f\!=\![u_1, (u_2,[\ze])]$ in (\ref{CaseStudy_e}) with the points of $L\!-\! \{L\cdot E_1, L\cdot E_2\}$.
The associated decorated dual graph $\Gamma$ is made of two vertices $v_1$ and $v_2$ corresponding to $u_1$ and $u_2$, with $I_{v_1}\!=\!\eset$ and $I_{v_2}\!=\!\{1,2\}$, connected by an edge $e$ with $I_e\!=\!\{1,2\}$ and $s_{\uvec{e}}\!=\!\pm (1,1)$ (depending on the choice of the orientation).
The group $\mc{G}(\Gamma)$ in this case is trivial and the function $s\colon\!\V\!\lra\!\R^2$ in Definition~\ref{LogMap_dfn}.\ref{Tropical_l} can be taken to be $s_{v_1}\!=\!(0,0)$ and $s_{v_2}\!=\!(1,1)$.

\item $u_1$ is of depth $\{1\}$: in this case $u_1$ comes with a holomorphic section $\ze'$ of $\cO_{\P^1}(1)$ as before. Since $\ord(z')\!=\!(1,1)$, by Definition~\ref{PreLogMap_dfn}.\ref{MatchingOrders_l}, $\ze'$ should be zero at $z'$ and this uniquely determines $[\ze']$. This unique element $f\!=\![(u_1,[\ze']),(u_2,[\ze])]$ corresponds to the point $E_1\cdot L$. The associated decorated dual graph $\Gamma$ is made of two vertices $v_1$ and $v_2$ corresponding to $u_1$ and $u_2$, with $I_{v_1}\!=\!\{1\}$ and $I_{v_2}\!=\!\{1,2\}$, connected by an edge $e$ with $I_e\!=\!\{1,2\}$ and $s_{\uvec{e}}\!=\!\pm (1,1)$ (depending on the choice of orientation).
The group $\mc{G}(\Gamma)$ in this case is trivial and the function $s\colon\!\V\!\lra\!\R^2$ in Definition~\ref{LogMap_dfn}.\ref{Tropical_l} can be taken to be $s_{v_1}\!=\!(1,0)$ and $s_{v_2}\!=\!(2,1)$.
\item $u_1$ is of depth $\{2\}$: similarly, there is a unique such map which corresponds to the point $E_2\cdot L$.
\eEn
\eEx

%%%%%%%%%%%%%%%
\subsection{Forgetful maps}\label{forgetful_ss}
In this section, we show that the process of forgetting some of the smooth components of an SNC divisor $D\!=\!\bigcup_{i\in [N]} D_i$ gives us a forgetful map between the corresponding log moduli spaces. The results are not used in the rest of the paper. While (\ref{FogetLog_e}) is not always an embedding, the map (\ref{ioSprod_e}) below is an embedding. This embedding can be used to reduce certain arguments to the case of smooth divisors.\\

\noindent
Let $D\!= \!\bigcup_{i\in [N]} D_i \!\subset\! X$ be an SNC symplectic divisor, $(\om,J)\!\in\!\cJ(X,D)$, $g,k\!\in\!\N$,  
\bEq{mfs_e}
\mfs\!=\!\big(s_a=(s_{ai})_{i\in [N]}\big)_{a=1}^k\!\in\!(\Z^N)^k,
\eEq
and $\Gamma\!\in\!\tn{DG}(g,\mfs,A)$.
Given $\mc{I}\!\subset\![N]$, 
let 
$$
\mfs|_{\mc{I}}\!=\!(s_a=(s_{ai})_{i\in \mc{I}})_{a=1}^k\!\in\!(\Z^{\mc{I}})^k,\quad D|_{\mc{I}}\!= \!\bigcup_{i\in \mc{I}} D_i,
$$
and $\Gamma|_{\mc{I}}\!\in\!\tn{DG}(g,\mfs|_{\mc{I}},A)$ be the decorated dual graph with the same set of vertices and edges, but with the reduced set of decorations 
$$
I'_v\!=\!I_v\cap \mc{I}\quad \forall~v\!\in\V, \quad I'_e\!=\!I_v\cap \mc{I}\quad \forall~e\!\in\E,\quad s'_{\uvec{e}}\!=\!(s_{\uvec{e},i})_{i\in \mc{I}}\!\in\!\Z^{ \mc{I}}\quad \forall~\uvec{e}\!\in\uvec{\E}.
$$
Define 
\bEq{prSS_e}
\iota_{[N],\mc{I}}\colon \cM^{\tn{plog}}_{g,\mfs}(X,D,A)_\Gamma\lra \cM^{\tn{plog}}_{g,\mfs|_{\mc{I}}}(X,D|_{\mc{I}},A)_{\Gamma|_{\mc{I}}}
\eEq
to be the (well-defined) forgetful map obtained by removing the meromorphic sections 
$$
(\ze_{v,i})_{i \in I_v -I'_v\subset [N]-\mc{I}}
$$ 
in (\ref{fplogSetUp_e}) for all $v\!\in\!\V$. 

\bLm{ForgetfulLog_lm}
With notation as above, the map $\iota_{[N],\mc{I}}$ in (\ref{prSS_e}) sends $\cM_{g,\mfs}(X,D,A)_\Gamma\!\subset\!\cM^{\tn{plog}}_{g,\mfs}(X,D,A)_\Gamma$ to $\cM_{g,\mfs|_{\mc{I}}}(X,D|_{\mc{I}},A)_{\Gamma|_{\mc{I}}}\!\subset\!\cM^{\tn{plog}}_{g,\mfs|_{\mc{I}}}(X,D|_{\mc{I}},A)_{\Gamma|_{\mc{I}}}$.
\eLm

\bPf
Fix an orientation $O$ on $\E$. With notation as in (\ref{DtoT_e}), the commutative diagram
$$
\xymatrix{
& \Z^\E\oplus \bigoplus_{v\in \V} \Z^{I_v}  \ar[rr]^{\vr} \ar[d]^{\tn{pr}_{\D}}&& \bigoplus_{e\in \E} \Z^{I_e}\ar[d]^{\tn{pr}_{\T}}\\
& \Z^\E\oplus \bigoplus_{v\in \V} \Z^{I'_v}  \ar[rr]^{\vr'} && \bigoplus_{e\in \E} \Z^{I'_e}\;,
}
$$
where $\tn{pr}_{\D}$ and $\tn{pr}_{\T}$ are the obvious projection maps and $\vr$ and $\vr'$ are defined via $O$, 
induces a group homomorphism $\tn{pr}_{[N],\mc{I}}\colon\!\mc{G}(\Gamma)\!\lra\!\mc{G}(\Gamma|_{\mc{I}})$ such that 
$$
\tn{pr}_{[N],\mc{I}}\big(\ob_\Gamma(f)\big)=\ob_{\Gamma|_{\mc{I}}}(\iota_{[N],\mc{I}}\big(f)\big)\qquad \forall~f\!\in\! \cM^{\tn{plog}}_{g,\mfs}(X,D,A)_\Gamma.
$$
Therefore, $\ob_\Gamma(f)\!=\!1$ implies $\ob_{\Gamma|_{\mc{I}}}(\iota_{[N],\mc{I}}\big(f)\big)\!=\!1$.
\ePf
\noindent
Taking union over all $\Gamma$, we obtain the stratified forgetful map 
$$
\iota_{[N],\mc{I}}\colon\ov\cM^{\tn{log}}_{g,\mfs}(X,D,A)\lra \ov\cM^{\tn{log}}_{g,\mfs|_{\mc{I}}}(X,D|_{\mc{I}},A).
$$
\noindent
For example, the $\mc{I}\!=\!\eset$ case of (\ref{prSS_e}) is the map $(\ref{FogetLog_e})$ into the underlying moduli space of stable maps; moreover, 
\bEq{SSSpr_e}
\iota_{[N],\mc{I}'}=\iota_{\mc{I},\mc{I}'}\circ \iota_{[N],\mc{I}}\colon \ov\cM^{\tn{log}}_{g,\mfs}(X,D,A)\lra \ov\cM^{\tn{log}}_{g,\mfs|_{\mc{I}'}}(X,D|_{\mc{I}'},A)\qquad \forall~\mc{I}'\!\subset\!\mc{I}\!\subset\![N].
\eEq

\noindent
For $\mfs$ as in (\ref{mfs_e}), let $\mfs_i\!=\mfs|_{\{i\}}\!=\!(s_{ai})_{a=1}^k\!\in\!(\Z)^k$, for all $i\!\in\![N]$, and define
\bEq{ioSprod_e}
\iota_{[N],1} \!=\prod_{i\in [N]}\iota_{[N],\{i\}} \colon \ov\cM^{\tn{log}}_{g,\mfs}(X,D,A)\lra \times_{i\in [N]}~ \ov\cM^{\tn{log}}_{g,\mfs_i}(X,D_i,A),
\eEq
where the right-hand side is the fiber product of 
$$
\big\{\iota_{\{i\},\eset}\colon \ov\cM^{\tn{log}}_{g,\mfs_i}(X,D_i,A)\!\lra\!\ov\cM_{g,k}(X,A)\big\}_{i\in [N]}.
$$
The map $\iota_{[N],1}$ is well-defined by (\ref{SSSpr_e}) and it is an embedding\footnote{by the results of Section~\ref{Cpt_s}, the maps $\iota_{[N],\mc{I}}$ and thus $\iota_{[N],1}$ are continuous.} by Remark~\ref{UniqueLog_rmk2}. As the following example shows, this embedding can be proper (i.e. not an equality).

\bEx{2dinP3toP2P2pr_ex}
In Example~\ref{2dinP3toP2P2_ex}, the obstruction groups $\mc{G}(\Gamma|_{\{1\}})$ and $\mc{G}(\Gamma|_{\{2\}})$ associated to $\Gamma|_{\{1\}}$ and $\Gamma|_{\{2\}}$ are trivial. Therefore, for an element of the right-hand side in (\ref{ioSprod_e}), the corresponding sections $\ze_{v_1,1}$ and $\ze_{v_2,2}$ can be arbitrary (modulo the combinatorial conditions imposed by Definitions~\ref{PreLogMap_dfn} and~\ref{LogMap_dfn}). On the other hand, for such a pair $(\ze_{v_1,1},\ze_{v_2,2})$ to define an element of the left-hand side, the corresponding group element in the non-trivial group $\mc{G}$ has be the identity.
Therefore, the restriction
$$
\iota_{\{1,2\},1}\colon 
\cM_{g,\mfs}(X,D,A)_{\Gamma}\lra \times_{i=1,2}~ \cM_{g,\mfs_i}(X,D_i,A)_{\Gamma|_{\{i\}}}
$$
of (\ref{ioSprod_e}) to $\ov\cM^{\tn{log}}_{g,\mfs}(X,D,A)_{\Gamma}$ is not an isomorphism.

\eEx

%%%%%%%%%%%%%%%

%------------------------------------------------------------------------------------------------------
\section{Compactness}\label{Cpt_s}
In this section, after a quick review of  the convergence problem for the Deligne-Mumford space and for the classical moduli spaces of $J$-holomorphic curves, we slightly rephrase and prove Theorem~\ref{Compactness_th} in several steps. The main step of the proof is Proposition~\ref{VertexOrder_prp}, that relates the sequence of ``gluing" and ``rescaling" parameters, when a sequence of $J$-holomorphic curves breaks into two pieces with at least one of them mapped into $D$. 

%----------------------------------------------------
\subsection{Classical Gromov convergence}\label{CptGromov_ss}

\bDf{def:cuttingset}
Given a $k$-marked genus $g$ (possibly not stable) nodal surface $C\equiv(\Si,\vec{z})$ with dual graph $\Gamma$, a 
\textbf{cutting configuration} with dual graph $\Gamma'$ is a set of disjoint embedded circles 
$$
\gamma\equiv\{\gamma_e\}_{e\in \E(\Gamma'/\Gamma)}\!\subset\! \Si,
$$ 
away from the nodes and marked points, such that  the nodal marked surface $(\Si',\vec{z}{\,'})$ obtained by pinching every $\gamma_e$ into a node $q_e$ has dual graph $\Gamma'$.
\eDf

\noindent
Thus, a cutting configuration corresponds to a continuous map
$$
\varphi_\gamma\colon C\lra C',
$$ 
denoted by a $\gamma$-\textbf{degeneration}\footnote{It is called ``deformation'' in \cite{RS}.} in what follows, onto a $k$-marked genus $g$ nodal surface $C'$ with dual graph $\Gamma'$
such that $\vec{z}{\,'}\!=\! \varphi_\gamma(\vec{z})$, the preimage of every node of $\Si$ is either a node in $\Si'$ or a circle in $\gamma$, and the restriction 
$$
\varphi_\gamma\colon \Si\!\setminus\! \gamma \lra \Si'\! \setminus \!(\varphi_\gamma (\gamma )\!\equiv \! \{q_e\}_{e\in \E(\Gamma'/\Gamma)})
$$
is a diffeomorphism. Let
\bEq{equ:deg-graphs}
\gamma^*\colon \Gamma'\lra \Gamma
\eEq
be the map corresponding to $\varphi_\gamma$ between the dual graphs. 
We have
$$
\E(\Gamma')\!\approx \!\E(\Gamma)\cup \E(\Gamma'/\Gamma), \qquad \L(\Gamma')\!\approx \!\L(\Gamma),
$$
such that $\gamma^*|_{\E(\Gamma)\subset \E(\Gamma')}$ and $\gamma^*|_{\L(\Gamma')}$ are isomorphisms and 
$$
\gamma^*\colon \E(\Gamma'/\Gamma)\lra \V(\Gamma)
$$
sends the edge $e$ corresponding to $\gamma_e$ to $v$, if $\gamma_e\!\subset\!\Sigma_v$.
For every ${v}'\!\in\!\V(\Gamma')$ there exists a unique $v\!\in\!\V(\Gamma)$ and a connected component $U_{v'}$ of $\Sigma_v\!\setminus\! \{\gamma_e\}_{e\in \E(\Gamma'/\Gamma)}$ such that $\Sigma'_{v'}\!\subset\!\Sigma'$ is obtained by collapsing the boundaries of $\tn{cl}({U}_{v'})$ ($\tn{cl}$ means closure). This identification determines the surjective map 
\bEq{gammaV_e}
\gamma^*\colon \V(\Gamma')\lra \V(\Gamma),\quad {v}'\lra v.
\eEq
From another perspective, a cutting configuration corresponds to expanding each vertex $v\!\in\!\V(\Gamma)$ into a sub-graph $\Gamma'_v\!\subset\!\Gamma'$ (some times, this involves just adding more loops to the existing graph) with the set of vertices and edges
$$
\V(\Gamma'_v)=(\gamma^*)^{-1}(v)\quad\tn{and}\quad \E(\Gamma'_v)=(\gamma^*)^{-1}(v)\!\cap \!\E(\Gamma'/\Gamma).
$$
Moreover, $g_v\!=\!g_{\Gamma'_v}$, the ordering of marked points are as before, and 
\bEq{equ:a-division}
A_v\!=\!\sum_{v'\in \V(\Gamma'_v)} A_{v'}.
\eEq
Figure~\ref{fig:cut-conf} illustrates a cutting configuration over a $1$-nodal curve of genus $3$ and the corresponding dual graphs. \\

\begin{figure}
\begin{pspicture}(7.8,1)(11,2.5)
\psset{unit=.26cm}

\pscircle*(35,9){.25}\pscircle*(39,9){.25}
\psline[linewidth=.07](35,9)(39,9)
\rput(35,10){$1$}\rput(39,10){$2$}
\psellipse[linewidth=.07](35,5)(5,2)
\psarc[linewidth=.07](35,7.8){3}{240}{300}
\psarc[linewidth=.07](35,2.2){3}{70}{110}
\psellipse[linestyle=dashed,dash=1pt,linewidth=.07](31.5,5)(0.5,1.46)
\pscircle*(40,5){.25}
\psellipse[linewidth=.07](46,5)(6,2)
\psarc[linewidth=.07](43,7.8){3}{240}{300}
\psarc[linewidth=.07](43,2.2){3}{70}{110}
\psarc[linewidth=.07](48,7.8){3}{240}{300}
\psarc[linewidth=.07](48,2.2){3}{70}{110}
\psellipse[linestyle=dashed,dash=1pt,linewidth=.07](45.5,5)(1.5,.5)

\psline[linewidth=1pt]{->}(53,5)(57,5)

\pscircle*(63,9){.25}\pscircle*(65,9){.25}\pscircle*(69,9){.25}
\psline[linewidth=.07](65,9)(69,9)\psline[linewidth=.07](63,9)(65,9)\pscircle(70,9){1}
\rput(63,10){$0$}\rput(65,10){$1$}\rput(68.7,10){$1$}
\psellipse[linewidth=.07](65,5)(5,2)
\psarc[linewidth=.07](65,7.8){3}{240}{300}
\psarc[linewidth=.07](65,2.2){3}{70}{110}
\pscircle[linewidth=.07](59,5){1}\pscircle*(60,5){.25}
\pscircle*(70,5){.25}
\psellipse[linewidth=.07](76,5)(6,2)\pscircle*(76,5){.25}
\psarc[linewidth=.07](73.3,7.4){3}{226}{300}
\psarc[linewidth=.07](75.9,1.8){3.2}{60}{120}
\psarc[linewidth=.07](78.5,7.4){3}{240}{313}
\psarc[linewidth=.07](73.3,2.6){3}{60}{125}
\psarc[linewidth=.07](75.9,8.2){3.2}{240}{300}
\psarc[linewidth=.07](78.5,2.6){3}{55}{120}

\end{pspicture}
\caption{On left, a $1$-nodal curve of genus $3$ and a cutting set made of two circles. On right, the resulting pinched curve.}
\label{fig:cut-conf}
\end{figure}

\noindent 
A sequence $\{\varphi_{\gamma_a}\colon C_a\lra C'\}_{a\in \N}$ of degenerations of marked nodal curves is called \textbf{monotonic} if $\Gamma(C_a)\!\cong\!\Gamma$ for some fixed $\Gamma$ and the induced maps $\gamma_a^*\colon \Gamma'\!\lra\! \Gamma$ are all the same. In this situation, the underlying marked nodal surfaces are isomorphic; i.e.,
\bEq{identified_e}
(C_a,\gamma_a)\cong \big((\Sigma,\mfj_a,\vec{z}),\gamma\big)\quad \forall~a\!\in\! \N,
\eEq
for some fixed marked surface $(\Si,\vec{z})$ with dual graph $\Gamma$ and cutting configuration $\gamma$. In the following,  we denote the complement of the set of nodes 
$$
 \{q_e\}_{e\in \E(\Gamma'/\Gamma)}\subset\!\Sigma'
$$ 
by $\Si'_*$.

\begin{definition}[{\cite[Dfn 13.3]{RS}}]\label{def:DM-convergence}
A sequence $\lrc{ C_a\!\equiv\!(\Si_a,\mfj_a,\vec{z}_a) }_{a\in \N}$ of genus $g$ $k$-marked nodal curves  \textbf{monotonically} converges to $C'\!\equiv \!(\Si', \mfj',\vec{z}{\,'} )$, if there exist a sequence of cutting configurations $\gamma_a$ on $C_a$ of type $\Gamma'$ and a monotonic sequence $\varphi_{\gamma_a}\colon C_a\!\lra\! C'$  of $\gamma_a$-degenerations such that the sequence $(\varphi_{\gamma_a}|_{\Si_a\setminus \gamma_a})_*\,\mfj_a $ converges to 
$\mfj'|_{\Si'_*}$ in the $C^\infty$-topology\footnote{uniform convergence on compact sets with all derivatives.}. 
\eDf

\noindent
By \cite[Sec~13]{RS}, the topology underlying the holomorphic orbifold structure of $\ov\cM_{g,k}$ is equivalent to the sequential $\tn{DM}$-convergence topology: a sequence $\lrc{ C_a}_{a\in \N}$ of genus $g$ $k$-marked stable nodal curves $\tn{DM}$-converges to $C'$ if a subsequence of that monotonically converges to $C'$.
The following result, known as Gromov's Compactness Thm \cite[Thm~1.5.B]{Gr}, describes a convergence topology on $\ov\cM_{g,k}(X,A,J)$  which is compact and metrizable; see \cite{Pan}, \cite[Thm~1.2]{Hu}, \cite[Thm 0.1]{Ye},  and \cite[Ch~5]{MS2} for further details.
In the special case of Deligne-Mumford space, Gromov convergence is equal to the DM-convergence discussed above.

\bTh{Gromov_Th}
Let $(X,\om)$ be a compact symplectic manifold,
$\{J_a\}_{a\in \N}$ be a sequence of $\om$-tame almost complex structures on $X$ converging in $C^\infty$-topology to $J$, and  
$$
\lrc{f_a \equiv \lrp{ u_a,C_a \equiv \lrp{\Si_a,\mfj_a,\vec{z}_a }}}_{a\in \N}
$$ 
be a sequence of stable $J_a$-holomorphic maps of bounded (symplectic) area into $X$. After passing to a subsequence, still denoted by $\{f_a\}_{a\in \N}$, there exists a unique (up to automorphism) stable $J$-holomorphic map 
$$
f'\!\equiv\!(u',C'\equiv(\Si',\mfj',\vec{z}{\,'}))
$$ 
such that $\{C_a\}_{a\in \N}$ monotonically converges to $C'$,
and 
\bEn
\item\label{l:smooth-conv} we can choose the $\gamma_a$-degeneration maps $\varphi_{\gamma_a}\colon \Si_a \!\lra\! \Si'$ of the monotonic convergence such that the restriction 
$$
u_a|_{\Sigma_a\setminus \gamma_a}\!\circ\! \varphi_{\gamma_a}^{-1}|_{\Si'_*}
$$ 
converges uniformly with all the derivatives to $u|_{\Si'_*}$ over compact sets;
\item with the dual graphs $\Gamma\!\cong\!\Gamma(C_a)$ and $\Gamma'\!=\!\Gamma(C')$ as in the definition of monotonic sequences, 
$$
\lim_{a\lra \infty} u_a(\gamma_{a,e})\!=\!u'(q_e)\qquad \forall~e\!\in\!\E(\Gamma'/\Gamma);
$$
\item\label{l:symp-area} symplectic area of $f'$ coincides with the symplectic area of $f_a$, for all $a\!\in\!\N$.
\eEn
\eTh
\noindent
It follows from the properties \ref{l:smooth-conv} and \ref{l:symp-area} that for every $v'\!\in\! \Gamma'$, with $U_{a,v'}\!\subset\! \Sigma_a$ as in the definition of $\gamma_a^*(v')$,
$$
\lim_{a\lra\infty} \int_{\tn{cl}(U_{a;v'})} {u_a}^*\om = \int_{\Si'_{v'}} (u')^*\om.
$$
Moreover, the stronger identity (\ref{equ:a-division}) holds. 
With respect to the identification of the domains and degeneration maps 
$$
(\varphi_{\gamma_a}\colon \Si_a\lra \Si')\cong (\varphi_{\gamma}\colon\Si\lra \Si')
$$ 
as in (\ref{identified_e}), the second property implies that the sequence $(u_a\!\colon\! \Si\!\lra\! X)_{a\in \N}$
$C^0$-converge to $u\circ \varphi_{\gamma}$. \\

\noindent
Assume $D\!\subset\! X$ is  an SNC symplectic divisor, $(\om,J)\!\in\!\cJ(X,D)$, and 
\bEq{OriginalSeq_e0}
\lrc{f_{a} \equiv \lrp{ u_a, [\ze_a] ,C_a \equiv \lrp{\Si_a,\mfj_a,\vec{z}_a}}}_{a\in \N}
\eEq
is a sequence of stable log maps in $\ov\cM^{\tn{log}}_{g,\mfs}(X,D,A)$. After passing to a subsequence, we may assume that all the maps in (\ref{OriginalSeq_e0}) have  the same decorated dual graph $\Gamma(\V,\E,\L)$, and that the underlying sequence of stable maps
\bEq{ClassicalPart_e}
\lrc{h_a \equiv \lrp{u_a,C_a \equiv \lrp{\Si_a,\mfj_a,\vec{z}_a}}}_{a\in \N}
\eEq
in $\ov\cM_{g,k}(X,A)$ (with the same domain) Gromov convergences to the stable map
$$
h\!\equiv\!(u,C\equiv(\Si,\mfj,\vec{z}))\in \ov\cM_{g,k}(X,A)
$$
as in Theorem~\ref{Gromov_Th}. Then, in order to prove Theorem~\ref{Compactness_th}, (for $J$ as in the statement of the theorem) after passing to a further subsequence, we prove  that $h$ lifts to a unique log map $f\!\in\!\ov\cM^{\tn{log}}_{g,\mfs}(X,D,A)$. The meromorphic sections that lift $h$ to the log map $f$ are specified in Section~\ref{CptLog_ss}. We first prove that $f$ is a pre-log map in Lemma~\ref{PreLog_lmm}; the proof works for arbitrary $(\om,J)\!\in\!\cJ(X,D)$. Then, in Proposition~\ref{LGSmoothV_prp}, we prove that $f$ satisfies the conditions of Definition~\ref{LogMap_dfn}. Since there are only finitely many possible log lifts of a stable map $f$ with different decorations on the dual graph, it follows with little effort that (\ref{FogetLog_e}) is a continuous local embedding.

%
%\begin{remark}
%A log map should be thought of as an stable map plus some partial tangent vector (to the moduli space) data. A sequence of stable maps $\{h_a\}_{a\in \N}$ converging to $h$ may approach it from various directions. Therefore, we need to further pass to a subsequence to pick one particular direction. 
%\end{remark}

%------------------------------------

\subsection{Log-Gromov convergence}\label{CptLog_ss}
In this section, first, we recall some basic structures associated to smooth/SNC symplectic divisors. Then we state the definition of log-Gromov convergence and a convergence result from which Theorem~\ref{Compactness_th} will be deduced.\vskip.1in

\noindent
Let $D \!\subset\! (X,\om)$ be  a smooth symplectic divisor, $J\!\in\!\cJ(X,D,\om)$, and $\mfi_{\cN_XD}$ be the induced complex structure on $\cN_XD$. 
Let $J_{X,D}$ be the almost complex structure on $\cN_XD$ induced by the $\dbar$-operator $\dbar_{\cN_XD}$ associated to $(\cN_XD,\mfi_{\cN_XD})$ as in the end of Section~\ref{dbar_ss}. 
Fix a compatible pair of a  Hermitian metric $\rho$ and a Hermitian connection $\nabla$ on $\cN_XD$. Such a  connection $\nabla$ defines a $1$-form $\al_\nabla$ on $\cN_XD-D$ whose restriction to each fiber $\cN_XD|_{p}-\{p\}\cong \C^*$ is the $1$-form $\nd\theta$ with respect to the polar coordinates $(r,\theta)$ determined by $\rho= r^2$ and the complex structure $\mfi_{\cN_XD}$. Recall from Section~\ref{dbar_ss} that the connection $\nabla$ gives a splitting 
\bEq{HorOVer_e}
T\cN_XD\cong \pi^*TD \oplus \pi^*\cN_XD,
\eEq
such that $J_{X,D}$ is equal to $\pi^*J_D$ on the first summand and  to $\pi^*\mfi_{\cN_XD}$ on the second one.
By the Symplectic Neighborhood Theorem~\cite[Thm~3.30]{MS1}, for $\cN'_XD$ sufficiently small, there exists 
a diffeomorphism
\bEq{RegPsi_e}
\Psi\colon\!\cN'_XD\!\lra\!X
\eEq
from a neighborhood of $D$ in $\cN_X D$ onto a neighborhood of $D$ in $X$ such that $\Psi(x)\!=\!x$, the isomorphism
\bEq{dPsi_e}
\cN_XD|_x\!=\!T_x^{\tn{ver}}\cN_XD \hookrightarrow T_x \cN_XD \stackrel{\nd_x\Psi}{\lra} T_xX \lra \frac{T_xX}{T_xD} \equiv \cN_XD|x
\eEq
is the identity map for every $x\!\in\!D$, and 
\bEq{om-Reg_e}
\Psi^*\om= \om_{X,D}=\pi^*(\om|_D) + \frac{1}{2}\nd(\rho\al_\nabla).
\eEq
The last property is not needed for many of the arguments in Section~\ref{CptLog_ss}.
In the language of \cite[Dfn~2.9]{FMZ1}, the tuple $\cR\!=\!(\rho,\nabla,\Psi)$ is called an $\om$-regularization. If $\Psi^*J=J_{X,D}$, then the tuple $(\om,\cR,J)$ is an element of $\tn{AK}(X,D)$ mentioned in (\ref{AKtoSymp_e}). In general, if $D\!=\!\bigcup_{i\in [N]} D_i$ is an SNC symplectic divisor in $(X,\om)$, a \textbf{system of regularizations for} $D$ in~$X$ is
a collection of smooth embeddings
$$
\Psi_I\!: \cN'_X D_I \lra X, \quad I \subset [N],
$$
from open neighborhoods $\cN'_XD_I\!\subset\!\cN_X D_I$ of~$D_I$ so that $\Psi_I|_{D_I}\!=\!\id_{D_I}$,
$\nd\Psi_I$ induces the identity map on~$\cN_XD_I\cong \bigoplus_{i\in I} \cN_XD_i|_{D_I}$, and
$$
\Psi_I\big(\cN_{I;I'}\!\cap\!\tn{Dom}(\Psi_I)\big)=D_{I'}\!\cap\!\tn{Im}(\Psi_I)
\qquad\forall~I'\!\subset\!I\!\subset\![N]\,.
$$
Here, 
$$
\pi_{I;I'}\colon \cN_{I;I'} \cong \bigoplus_{i\in I-I'} \cN_XD_i|_{D_I}\lra D_I
$$
is the normal bundle of $D_{I}$ in $D_{I'}$.
The last identity implies that the derivative $\nd\Psi_I$ induces an isomorphism of split vector bundles
\bEq{fDIIdfn_e}
\mf{D}\Psi_{I;I'}\!:  \pi_{I;I'}^*\cN_{I;I-I'}\big|_{\cN_{I;I'}\cap\tn{Dom}(\Psi_I)}
\lra \cN_XD_{I'}\big|_{D_{I'}\cap\tn{Im}(\Psi_I)};
\eEq
see \cite[Sec~2.2]{FMZ2}. A \textbf{regularization} for $D$ in $X$ is a
system of regularizations for $D$ in~$X$ as above
satisfying the compatibility conditions
$$
\tn{Dom}(\Psi_I)=\mf{D}\Psi_{I;I'}^{-1}\big(\tn{Dom}(\Psi_{I'})\big), 
\quad  \Psi_I\!=\!\Psi_{I'}\!\circ\!\mf{D}\Psi_{I;I'}|_{\tn{Dom}(\Psi_I)}
\qquad \forall~I'\!\subset\!I\!\subset\![N]\,.
$$
\vskip.05in

\begin{definition}[{\cite[Dfn~2.9]{FMZ2}}]\label{omRegV_dfn}
An \textbf{$\om$-regularization} for $D$ in $X$ consists of a choice of Hermitian structure
$(\mfi_{I;i},\rho_{I;i},\na^{(I;i)})$ on $\cN_X D_i|_{D_I}$
for all $i\!\in\!I\!\subset\![N]$ together with
a regularization for~$D$ in~$X$ as above so that
$$
\Psi_I^*\om\!=\!\pi^*(\om|_{D_I})+\frac{1}{2} \sum_{i\in I} \nd(\rho_{I;i} \al_{\na^{(I;i)}})
\quad \forall~I\subset [N],
$$
and (\ref{fDIIdfn_e}) is an isomorphism of split Hermitian vector bundles
for all $I'\!\subset\!I\!\subset\![N]$. 
\end{definition}
\noindent
Finally, an element of $\tn{AK}(X,D)$ is a tuple of $(\om,\cR,J)$ where $\cR$ is an $\om$-regularization as in Definition~\ref{omRegV_dfn} and 
$$
\Psi_I^* J= \pi_I^*(J|_{TD_I}) \oplus  \bigoplus_{i\in I} \pi_I^* \,\mfi_{I;i}
$$
with respect to the decomposition (\ref{HorOVer_e}).
The main reason for restricting to $\tn{AK}(X,D)$ or the integrable almost complex structures in Theorem~\ref{Compactness_th} is that in the proof of Proposition~\ref{VertexOrder_prp}, for any $p\!\in\! D_I$, we need $J$ to be $(\C^*)^I$-equivariant in a neighborhood of $p$ with respect to a (local) $(\C^*)^I$-action that preserves $D$ and fixes $D_I$.\\

\noindent
For any $c\!\in\!\R_{>0}$, define 
$$
\cN_XD(c)\!=\!\{v\!\in\! \cN_XD\colon \rho(v)\!<\!c\}.
$$  
For any $t\!\in\!\C^*$, define
\bEq{Psit_e}
\aligned
&R_t\colon \cN_XD\lra\cN_XD, \qquad R_t(v)=tv\quad \forall v\!\in\!\cN_XD,\\
&\Psi_t=\Psi\circ R_t \colon R_t^{-1}(\cN'_XD)\lra X,\qquad J_t=\Psi_t^*J.
\endaligned
\eEq
Note that if $\Psi^*J=J_{X,D}$ then $J_t\equiv J_{X,D}$ is independent of $t$. The following lemma is an expansion of the sentence after \cite[(6.5)]{IP1}.

\bLm{JLemma_e}
For $J$ satisfying (\ref{intInnormal_e0}), we have 
$$
\lim_{t\lra 0} J_t|_{\ov{\cN_XD(c)}}= J_0:= J_{X,D}|_{\ov{\cN_XD(c)}}\qquad \forall~c\!\in\!\R_{>0}
$$
uniformly with all derivatives.
\eLm
\bPf
In order to simplify the notation, let us forget about $\Psi$ and think of $J$ as an almost complex structure on $\cN'_XD$ itself; then $J|_{D}\!=\!J_{X,D}|_D$ and $J_t\!=\!R_t^*J$, for every $t\!\in\!\C^*$. 
Via (\ref{HorOVer_e}), we decompose $J$ into $4$ components
$$
\aligned
&J_v(\al)= \lrp{J_v^{11}(\al_1)+J_v^{2 1}(\al_2)}\oplus 
 \lrp{J_v^{12}(\al_1)+J_v^{2 2}(\al_2)},\\ 
& \forall~x\!\in\!D,~v\!\in\!\cN_XD|_{x},~\al=\al_1\oplus\al_2\!\in\!\big(\pi^*TD \oplus \pi^*\cN_XD\big)|_{v},
\endaligned
$$
where, for example, $J^{11}$ is the component which maps the horizontal subspace $\pi^*TD$ to itself. Identifying $\al_1$ and $\al_2$ with the corresponding vectors in $T_xD$ and $\cN_XD|_x$, respectively, we get 
$$
(J_t)_v(\al)= \lrp{J_{tv}^{11}(\al_1)+J_{tv}^{2 1}(t\al_2)}\oplus 
 \lrp{\frac{1}{t}J_{tv}^{12}(\al_1)+J_{tv}^{2 2}(\al_2)}.
 $$
On each compact set $\ov{\cN_XD(c)}$, the first summand uniformly converges to $J_D(\al_1)$, and $J_{tv}^{2 2}(\al_2)$ uniformly converges to $\mfi_{\cN_XD}(\al_2)$ (with all the derivatives). Finally, the term 
 $$
 \frac{1}{t}J_{tv}^{12}(\al_1)
 $$
 $C^\infty$-converges to the normal part of (a multiple of) $N_J(v,J\al_1)$, which is zero by (\ref{intInnormal_e0}); see Remark~\ref{dbarvsD_rmk}.
\ePf

\noindent
For any (continuous) map $u\colon\!\Si\!\lra\!\cN_{X}D$, let
$$
\ov{u}=\pi\circ u\colon \Si\!\lra\!D
$$
denote its projection to $D$. Then $u$ is equivalent to a section $\ze\!\in\!\Gamma(\Si,\ov{u}^*\cN_XD)$ in the sense that 
\bEq{uzeCor_e}
u(x)\!=\!\ze(x)\!\in\!\cN_XD|_{\ov{u}(x)}\qquad \forall~x\!\in\!\Si.
\eEq 
We  will use this correspondence repeatedly in the following arguments. 
In particular, by \ref{HoloProj_l}-\ref{HoloSec_l} in Page~\pageref{HoloProj_l}, $u$ is $J_{X,D}$-holomorphic if and only if $\ov{u}$ is $J_D$-holomorphic and $\dbar_{\cN_XD}\ze\!=\!0$.

\bDf{STNS1_dfn} 
With $(X,D,\om,J,\Psi)$ as above (i.e. $D$ is smooth),
let 
$$
\bigg(f_{a} \equiv \big( u_a, C_a \equiv \lrp{\Si_a,\mfj_a,\vec{z}_a}\big)\bigg)_{a\in \N}
$$
be a sequence of stable maps with smooth domain in $\cM_{g,\mfs}(X,D,A)$ that Gromov converges, considered as a sequence in $\ov\cM_{g,k}(X,A)$, to the marked nodal map
$$
f\!\equiv\!(u_v,C_v\equiv(\Si_v,\mfj_v,\vec{z}_v))_{v\in \V}\in \ov\cM_{g,k}(X,A)
$$
with dual graph $\Gamma\!=\!\Gamma(\V,\E,\L)$ and nodal domain $\Si\!=\!\bigcup_{v\in \V}\Si_v$. 
With notation as in (\ref{nodalcurve_e}), (\ref{Divide_e}), and Theorem~\ref{Gromov_Th}, for each $v\!\in\!\V_1$, we say \textbf{$(u_a)_{a \in \N}$ is asymptotic to
 $$
 \ze_{v}\in \Om_{\tn{mero}}(\Si_v,u_v^* \cN_XD)
 $$  
on $\Si_v$ in the normal direction to $D$} if there exists a sequence of non-zero complex numbers $(t_{a,v})_{a\in \N}$ satisfying
\bEq{Convtoze_e}
\tn{(uniformly)  }\lim_{a\lra \infty} \Psi_{t_{a,v}}^{-1}\circ u_a \circ \varphi_{\gamma_a}^{-1}|_{K} = \ze_v|_K
\eEq
 in the sense of (\ref{uzeCor_e}) for every compact set $K\!\subset\!\Si_{v}\!-\!q_v$ .
\eDf

\noindent
Proposition~\ref{STNS1_prp} below shows that, after passing to a subsequence, the limiting $J$-holomorphic map $f$ always admits such meromorphic sections $\ze_v$, and that they are unique up to multiplication by a constant in~$\C^*$. Since $\nd\Psi$ in (\ref{dPsi_e}) is supposed to be the identity map on $\cN_XD$, (\ref{Convtoze_e}) does not depend on the particular choice of $\Psi$ in~(\ref{RegPsi_e}).

\bDf{LogConv_dfn}
Let $D\!\subset\! X$ be an SNC symplectic divisor, $(\om,J)\!\in\!\cJ(X,D)$, and 
$$
\bigg(f_a\equiv \big(u_{a,v},[\ze_{a,v}]=([\ze_{a,v,i}])_{i\in I_v},C_{a,v}=(\Si_v,\mfj_{a,v},\vec{z}_v))_{v\in \V}\bigg)_{a\in \N}\!\in\!\ov\cM^{\tn{log}}_{g,\mfs}(X,D,A)
$$
be a sequence of stable log maps in\footnote{More precisely, they represent equivalence classes of elements in $\ov\cM^{\tn{log}}_{g,\mfs}(X,D,A)$. }   $\ov\cM^{\tn{log}}_{g,\mfs}(X,D,A)$ with a fixed decorated dual graph $\Gamma\!=\!\Gamma(\V,\E,\L)$.
We say this sequence \textbf{log-Gromov converges} to the log (resp. pre-log) map
$$
f'=\big(u_{v'},[\ze_{v'}]=([\ze_{v',i}])_{i\in I_{v'}},C_{v'}\big)_{v'\in \V'}
$$
in $\ov\cM^{\tn{log}}_{g,\mfs}(X,D,A)$ (resp. $\cM^{\tn{plog}}_{g,\mfs}(X,D,A)_{\Gamma'}$) with the decorated dual graph $\Gamma'\!=\!\Gamma(\V',\E',\L')$  if the underlying sequence of stable maps in $\ov\cM_{g,k}(X,A)$ Gromov converges to the underlying marked nodal map
\bEq{Stable-limit_e}
\iota(f')= \big(u_{v'},C_{v'}\big)_{v'\in \V'}\!\in\!\ov\cM_{g,k}(X,A)
\eEq
with the nodal domain $\Si'\!=\!\bigcup_{v'\in \V'}\Si_{v'}$, 
and the followings hold. With $\gamma^*\colon\!\V'\!\lra\!\V$ as in (\ref{gammaV_e}) and notation as in Theorem~\ref{Gromov_Th}, for each $v\!\in\!\V$ and $v'\!\in\!\V'$, with $\gamma^*(v')\!=\!v$, \begin{itemize} 
\item if $i\!\in\!I_{v'}\!-\!I_v$,  $(u_{a,v})_{a \in \N}$ is asymptotic to $\ze_{v',i}$ on $\Si_{v'}$ in the normal direction to $D_i$ in the sense of Definition~\ref{STNS1_dfn}; \item 
 if $i\!\in\!I_v$, there exists a sequence 
$(t_{a,v',i})_{a\in \N}\!\in\!\C^*$ such that for every compact set $K\!\in\!\Si_{v'}\!-\!(z_{v'} \cup q_{v'})$, $t_{a,v',i}^{-1} ~\ze_{a,v,i} \circ \varphi_{\gamma_a}^{-1}|_{K}$ uniformly converges to $\ze_{v',i}|_{K}$.
 \end{itemize}
\eDf

\bTh{Log_thm11} 
Assume $D\!\subset\! X$ is an SNC symplectic divisor, $(\om,\cR,J)\!\in\!\tn{AK}(X,D)$ for some regularization $\cR$ or $J$ is integrable, and 
\bEq{faSequence_e11}
\big(f_a\equiv \big(u_{a,v},[\ze_{a,v}]=([\ze_{a,v,i}])_{i\in I_v},C_{a,v}=(\Si_v,\mfj_{a,v},\vec{z}_v))_{v\in \V}\big)_{a\in \N}\!\in\!\ov\cM^{\tn{log}}_{g,\mfs}(X,D,A)
\eEq
is a sequence of log maps in $\ov\cM^{\tn{log}}_{g,\mfs}(X,D,A)$. After passing to a subsequence, there exists a unique log map (up to reparametrization)
 \bEq{LogLimit_e11}
f'=\big(u_{v'},[\ze_{v'}]=([\ze_{v',i}])_{i\in I_{v'}},C_{v'}\big)_{v'\in \V'}
\eEq
such that (\ref{faSequence_e11}) log-Gromov converges to (\ref{LogLimit_e11}) in the sense of Definition~\ref{LogConv_dfn}.
\eTh

\noindent
We break the proof Theorem~\ref{Log_thm11} into smaller steps. The main steps are proved in the subsequent sections. \\

\noindent
For two sequences of non-zero complex numbers $(t_a)_{a\in \N}$ and $(t'_a)_{a\in \N}$, we write 
\bEq{ESeq_e}
(t_a)_{a\in \N}\sim (t'_a)_{a\in \N}\qquad \tn{if}\qquad \lim_{a\lra \infty} t_a/t'_a\!=\!1.
\eEq 
The right-hand side of (\ref{ESeq_e}) defines an equivalence relation on the set of such sequences and we denote the equivalence class of a sequence $(t_a)_{a\in \N}$ by  $[(t_a)_{a\in \N}]$. For  an equivalence class $[(t_a)_{a\in \N}]$ and $t\!\in\!\C^*$, the equation 
$$
t [(t_a)_{a\in \N}]:=[(tt_a)_{a\in \N}]
$$
is well-defined and defines an action of $\C^*$ on the set of equivalence classes. Moreover, the operation of point-wise multiplication/divison  between such sequences 
$$
(t_a)_{a\in \N}\cdot  (t'_a)_{a\in \N} = (t_at'_a)_{a\in \N}
$$
descends to a well-defined multiplication/division operation between the equivalence classes.\\

\noindent
The following proposition corresponds to \cite[Proposition 6.6]{IP1}.

\begin{remark}
There is a minor issue in the proof of  \cite[Proposition 6.6]{IP1}. In  \cite[(6.13)]{IP1}, the authors use Intermediate Value Theorem to find the right rescaling parameter $t\!=\!t_m$. However, the energy function used there is not necessarily continuous in $t$. For example, apply their argument to the example where $X=C_1\times C_2$ is the product of two curves, the divisor $V$ is $\{p\}\times C_2$ for some point $p\in C_1$, and the sequence of curves is $\{p_i\times C_2\}_{i=1}^\infty$, with $\lim_{i\to \infty } p_i=p$. \end{remark}

\bPr{STNS1_prp} As in Definition~\ref{STNS1_dfn} (i.e. $D$ is smooth), let 
\bEq{OriginalSeq_e}
\lrc{f_{a} \equiv \lrp{ u_a, C_a \equiv \lrp{\Si_a,\mfj_a,\vec{z}_a}}}_{a\in \N}
\eEq
be a sequence of stable maps with smooth domain in $\cM_{g,\mfs}(X,D,A)$ that Gromov converges, considered as a sequence in $\ov\cM_{g,k}(X,A)$, to the marked nodal map
$$
f\!\equiv\!\big(u_v,C_v\equiv(\Si_v,\mfj_v,\vec{z}_v)\big)_{v\in \V}\in \ov\cM_{g,k}(X,A).
$$ 
After passing to a subsequence (which we still denote by $\N$), for each $v\!\in\!\V_1$, there exists a unique 
 $$
 [\ze_{v}]\in \Om_{\tn{mero}}(\Si_v,u_v^* \cN_XD)/\C^*
 $$  
such that $(u_a)_{a \in \N}$  is asymptotic to $\ze_{v}$ on $\Si_v$ in the normal direction to $D$ in the sense of Definition~\ref{STNS1_dfn}. Furthermore, $\ze_v$ has no pole/zero in $\Si_{v}\!-\!(q_v\cup z_v)$, and  it has a zero of order $s_i$ at $z^i$, for all $z^i\!\in\!\vec{z}_v$.
\ePr

\bPf
For every fixed $K\!\in\!\Si_v-q_v$, by Theorem~\ref{Gromov_Th}, the sequence
$$
\ov{u}_{a,K}= \pi \circ u_{a,K} \colon K \lra D,\quad\tn{with}\quad u_{a,K}\equiv \Psi^{-1}\circ u_a \circ \varphi_{\gamma_a}^{-1}|_{K}\colon K\lra \cN'_XD\qquad \forall~a\!>\!\!>\!1,
$$
converges uniformly with all the derivative to $u_v|_K$, and 
$$
u_{a,K}(z)\!=\!\ze_{a,K}(z) 
$$
for some non-trivial smooth section $\ze_{a,K}\!\in\!\Gamma(K, \ov{u}_{a,K}^*\cN_XD)$ in the sense of (\ref{uzeCor_e}), such that the sequence $\ze_{a,K}$ converges uniformly with all the derivatives (with respect to a connection $\nabla$) to $0$. Choose $(t_{a,v,K})_{a>\!\!>1}$ so that 
\bEq{VNorm_e}
||t^{-1}_{a,v,K}\ze_{a,K}||_{L^\infty(K)}= c_K\qquad  \forall~a\!>\!\!>\!1
\eEq
for some arbitrary non-zero constant $c_K$.
Then, by \cite[Thm~4.1.1]{MS2} (after passing to a subsequence), the sequence  
$$
(\Psi_{t_{a,v,K}}^{-1}\circ u_v \circ \varphi_{\gamma_a}^{-1}|_{K})_{a>\!\!>1}
$$ 
of $J_{t_{a,v,K}}$-holomorphic maps in $\cN_XD(c_K)$  converges uniformly with all the derivatives to a unique $J_{X,D}$-holomorphic map 
$$
u_{\infty,K}\colon\!K\!\lra\!\cN_XD(c_K).
$$
By (\ref{VNorm_e}), \ref{HoloSec_l} on Page~\pageref{HoloSec_l}, and since $\ov{u}_{a,v}$ converges to $u_v|K$, we have 
$$
\ov{u}_{\infty,K}\!=\!u_v|K\qquad\tn{and}\qquad u_{\infty,K}\!=\!\ze_{v,K}
$$ 
for some unique non-trivial $\dbar_{\cN_XD}$-holomorphic section $\ze_{v,K}$ of $u_v^* \cN_XD|_{K}$. Since $\ze_{a,K}$ is non-zero away from $\vec{z}_a\cap\varphi_{\gamma_a}^{-1}(K)$, $\ze_{v,K}$ is non-zero away from $\vec{z}_v \cap K$. \\

\noindent
Let 
$$
K_1\subset K_2 \subset \cdots
$$
be a sequence exhausting $\Si_{v}\!-\!q_v$. For each $K_i$ let $\ze_{v,K_i}$ and $c_{K_i}$ be the section and constant corresponding to $K_i$, respectively, in the argument above.  Choose a reference point $p\!\in\!K_1$ and  fix a non-zero vector $v_p\!\in\!\cN_XD|_{u_v(p)}$. For each $i$, we can equally rescale $c_{K_i}$ and $(t_{a,v,K_i})_{a>\!\!>1}$ by a constant number in $\C^*$ so that $\ze_{v,K_i}(p)$ becomes equal to $v_p$. Then, by the uniqueness of the limiting section, we get 
$$
\ze_{v,K_i}\!=\!\ze_{v,K_{i+1}}|_{K_i}\qquad \forall~i\!\in\!\N.
$$ 
Therefore, the equation
$$
\ze_v(x):=\ze_{v,K_i}(x)\qquad \forall~x\!\in\!\Si_{v}\!-\!q_v,~i\!\in\!\N~~\tn{s.t.}~~x\!\in\!K_i,
$$
defines a holomorphic section of $u_v^*\cN_XD|_{\Si_{v}-q_v}$ such that (\ref{Convtoze_e}) holds. Moreover, 
$$
(t_{a,v,K_i})_{a\in \N}\sim (t_{a,v,K_j})_{a\in \N}\qquad \forall i,j \in \N.
$$
It remains to show that $\ze_v$ has at most finite order poles at the nodes and $\tn{ord}_{z^i}(\ze_v)\!=\!s_i$ for all $z^i\!\in\!\vec{z}_v$.\\

\noindent
For any marked point $z^i\!\in\!\vec{z}_v$, let $\De_i\!\subset\!\Si_v$ be a sufficiently small disk around $z^i$ that contains no other marked point or nodal point. For $a$ sufficiently large, the order of vanishing of $u_a$ at $z_a^i$ is equal to the winding number of 
$$
 \Psi_{t_{a,v,\partial\ov\De_i}}^{-1}\circ u_a \circ \varphi_{\gamma_a}^{-1}|_{\partial\ov\De_i}\qquad  \forall~a\!>\!\!>\!1
 $$ 
 around $D$. With $K\!=\!\ov\De_i$ in (\ref{Convtoze_e}), these numbers are the same for  $a\!>\!\!>\!1$ and they are equal to the winding number of ${u}_{\infty,\ov\De_i}|_{\partial\ov\De_i}$ around $D$. The latter is equal to the order of $\ze_v$ at $z^i$. We conclude that the contact orders stay the same at the marked points. \\

\noindent
Similarly, for any nodal point $q_{\uvec{e}}\!\in\!\Si_v$, with $\uvec{e}\!\in\!\uvec{\E}$ and $v_1(\uvec{e})\!=\!v$, let $\De\!\subset\!\Si_v$ be a sufficiently small disk around $q_{\uvec{e}}$ that contains no other marked point or nodal point. Choose a compact set 
$K\!\subset\!\Si_{v}\!-\!q_v$ so that one of whose boundary circles coincides with $\partial \ov\De$. Since the convergence in (\ref{Convtoze_e}) is uniform, the winding numbers of
$$
 \Psi_{t_{a,v,K}}^{-1}\circ u_a \circ \varphi_{\gamma_a}^{-1}|_{\partial\ov\De}\qquad  \forall~a\!>\!\!>\!1
 $$ 
 around $D$ are the same as the winding number of ${u}_{\infty,K}|_{\partial\ov\De}$ around $D$. The latter is equal to the order of $\ze_v$ at $q_{\uvec{e}}$. We conclude that $\ze_v$ extends to a meromorphic section at $q_{\uvec{e}}$.
\ePf

\begin{remark}
Note that that the sections $\ze_v$ and the equivalence class of the rescaling sequence $[(t_{a,v})_{a\in \N}]$ are independent of the choice of $\Psi$. It is also clear from (\ref{Convtoze_e}) that if $(t_{v,a})_{a\in \N}$ is a rescaling sequence associated to $\ze_v$ and  $(t'_{v,a})_{a\in \N}$ is a rescaling sequence associated to $c\ze_v$, for any $c\!\in\!\C^*$, then
\bEq{Rescalingt_e}
[(t'_{v,a})_{a\in \N}]=c^{-1}[(t_{v,a})_{a\in \N}].
\eEq
\end{remark}

\noindent
The following is the analogue of Proposition~\ref{STNS1_prp} for a sequence of stable log maps with smooth domain and image in $D$.

\bCr{STNS1_cr} If $D$ is smooth, suppose 
\bEq{OriginalSeq_e2}
\lrc{f_{a} \equiv \lrp{ u_a, \ze_a,C_a \equiv \lrp{\Si_a,\mfj_a,\vec{z}_a}}}_{a\in \N}
\eEq
is a sequence of representatives of stable log maps with smooth domain in $\cM_{g,\mfs}(X,D,A)_{\{1\}}$ such that the underlying sequence of stable $J_D$-holomorphic maps
\bEq{OriginalSeq_e3}
\lrc{f_{a} \equiv \lrp{ u_a, C_a \equiv \lrp{\Si_a,\mfj_a,\vec{z}_a}}}_{a\in \N}
\eEq
Gromov converges, as a sequence in $\ov\cM_{g,k}(D,A)$, to the nodal map
$$
f\!\equiv\!\big(u_v,C_v\equiv(\Si_v,\mfj_v,\vec{z}_v)\big)_{v\in \V}\in \ov\cM_{g,k}(D,A).
$$
With notation as in (\ref{nodalcurve_e}), (\ref{Divide_e}), and Theorem~\ref{Gromov_Th}, after passing to a subsequence (which we still denote by $\N$), for every $v\!\in\!\V$, there exists a unique  
$$
 [\ze_{v}]\in \Om_{\tn{mero}}(\Si_v,u_v^* \cN_XD)/\C^*
 $$ 
 and a unique equivalence class of sequences of non-zero complex numbers $[(t_{a,v})_{a\in \N}]$ such that 
\bEq{LinSTNS1_e2}
  \lim_{a\lra \infty} t_{a,v}^{-1} ~\ze_a \circ \varphi_{\gamma_a}^{-1}|_{K} = \ze_v|_K,
\eEq
for any compact set $K\!\subset\!\Si_{v}\!-\!q_v$. 
Furthermore, $[\ze_v]$ only depends on the sequence of equivalence classes $\big([\ze_a]\big)_{a\in \N}$, it has no pole/zero in $\Si_{v}- (q_v\cup z_v)$, and it has a zero/pole of the same order $s_i$ at $z^i$, for all $z^i\!\in\!\vec{z}_v$. 
\eCr

\bPf
If (\ref{LinSTNS1_e2}) holds for a sequence $(\ze_a,t_{a,v})_{a\in \N}$, then it also holds for any other simultaneous reparametrization $(t_a\ze_a,t_at_{a,v})_{a\in \N}$. Therefore, (\ref{LinSTNS1_e2}) only depends on the sequence of equivalence classes $\big([\ze_a]\big)_{a\in \N}$.
Every map in the sequence (\ref{OriginalSeq_e2}) corresponds to a $J_{X,D}$-holomorphic map in 
$$
\cM_{g,\mfs}(\cN_XD,D,A).
$$
We can choose the representatives $\ze_a$ so that their image in $\cN_XD$ lie in an arbitrary small compact neighborhood\footnote{so that we can still apply the Gromov convergence theorem. We can also use the compact manifold $\P_XD$ in (\ref{PXV_e}) instead of $\cN_XD$ with the symplectic form $
\om_{X,D}= \pi^*(\om|D) + \ep \nd \big(\frac{\rho\al_{\na}}{1+\rho}\big)
$, where $\ep\!>\!0$ is a sufficiently small constant. Then, for $t$ sufficiently small, by interpolating between 
$$J_t|_{R_t^{-1}(\cN'_XD)}\quad\tn{and}\quad J_{X,D}|_{\P_XD},$$
we can construct a family of almost complex structures $\wt{J}_t$ on $\P_XD$ such that  $\wt{J}_t$ converges to $J_{X,D}$; see \cite[Prp~6.6]{IP1}.} of $D$.
Replacing $(X,D,\om,J)$ with $(\cN_XD,D,\om_{X,D},J_{X,D})$ and $\Psi$ with the identity map in Proposition~\ref{STNS1_prp}, we get the desired result.
\ePf

\noindent
From Proposition~\ref{STNS1_prp} and Corollary~\ref{STNS1_cr} we derive the following conclusion.

\bLm{PreLog_lmm} 
Let $D\!\subset\! X$ be an SNC symplectic divisor, $(\om,J)\!\in\!\cJ(X,D)$, and 
\bEq{faSequence_e}
\big(f_a\equiv \big(u_{a,v},[\ze_{a,v}]=([\ze_{a,v,i}])_{i\in I_v},C_{a,v}=(\Si_v,\mfj_{a,v},\vec{z}_v))_{v\in \V}\big)_{a\in \N}\!\in\!\ov\cM^{\tn{log}}_{g,\mfs}(X,D,A)
\eEq
be a sequence of stable log maps in $\ov\cM^{\tn{log}}_{g,\mfs}(X,D,A)$. After passing to a subsequence, there exists a unique pre-log map
 \bEq{LogLimit_e}
f'=\big(u_{v'},[\ze_{v'}]=([\ze_{v',i}])_{i\in I_{v'}},C_{v'}\big)_{v'\in \V'}
\eEq
such that (\ref{faSequence_e}) log-Gromov converges to (\ref{LogLimit_e}) in the sense of Definition~\ref{LogConv_dfn}.
\eLm 

\bPf
First, we apply the Gromov convergence to the underlying sequence of stable maps. Then, running through all $D_i$ and $v\!\in\! \V$, one at a time, by applying Proposition~\ref{STNS1_prp} (with $D=D_i$) to the sequence  
$$
\big(u_{a,v},C_{a,v}\big)_{a\in \N},
$$
whenever $i\!\not\in\! I_v$, and Corollary~\ref{STNS1_cr}  to the sequence  
$$
\big(u_{a,v},\ze_{a,v,i},C_{a,v}\big)_{a\in \N},
$$
whenever $i\!\in\! I_v$, we obtain $f'$. We need to show that $f'$ satisfies the conditions of Definition~\ref{PreLogMap_dfn}. The first condition is obviously satisfied.\\

\noindent
\textbf{Continuity.} The matching condition~\ref{MatchingPoints_l} of Definition~\ref{PreLogMap_dfn} is about the continuity of the underlying stable map $f'$ and already holds by the Gromov compactness.\\

\noindent
\textbf{Contact orders at the nodes.} In order to show that the condition~\ref{MatchingOrders_l} of Definition~\ref{PreLogMap_dfn} is satisfied, let us first fix some notation first. Since  
$$
s_{\uvec{e}}\!=\!-s_{\uvec{e}}\quad \Leftrightarrow\quad  s_{\uvec{e},i}\!=\!-s_{\uvec{e},i} \!\in\! \Z \qquad \forall i\in [N],
$$
it is enough to show that the condition~\ref{MatchingOrders_l} is satisfied relative to each smooth component $D_i$; i.e. we may assume $D$ is smooth. In the context/notation\footnote{Note that the notation used for the limiting map in Proposition~\ref{STNS1_prp} is different than the one in the statement of Lemma~\ref{PreLog_lmm}.}  of Proposition~\ref{STNS1_prp}, for every $v,v'\!\in\!\V$ and any node $q_{e}\!=\!(q_{\uvec{e}}\sim q_{\scz\ucev{e}})$, with $\uvec{e}\!\in\!\uvec{\E}_{v,v'}$, connecting $\Si_{v}$ and $\Si_{v'}$, let 
$\De_{\uvec{e}}\!\subset\!\Si_{v}$ be a sufficiently small disk around $q_{\uvec{e}}$ (not containing any other marked point or nodal point), $\De_{\scz\ucev{e}}\!\subset\!\Si_{v'}$ be a sufficiently small disk around $q_{\scz\ucev{e}}$, and $A_e\!=\! \De_{\uvec{e}} \cup \De_{\scz\ucev{e}}$ be a the resulting neighborhood of $q_e$ in $\Si$. We orient each circle $\partial\De_{\uvec{e}}$ in the direction of the counter-clock wise rotation in $\De_{\uvec{e}}\!\subset\!\C$.
For each $e\!\in\!\E$, $A_{a,e}\!=\!\varphi_{\gamma_a}^{-1}( {A}_e)$ is a cylinder in $\Si_a$ with two (oppositely oriented) boundaries 
\bEq{partialAe_e}
\partial A_{a,\uvec{e}}=\varphi_{\gamma_a}^{-1}(\partial\De_{\uvec{e}}) \qquad \tn{and}\quad \partial A_{a,\scz{\ucev{e}}}=\varphi_{\gamma_a}^{-1}(\partial\De_{\scz\ucev{e}})
\eEq
such that $u_a|_{A_{a,e}}$ does not intersect $D$ for $a\!>\!\!>\!1$. Since $u_a|_{A_{a,e}}$ is continuous and does not intersect $D$, the winding numbers of $u_a$ around $D$ on the two boundary circles of the annulus $A_{a,e}$ (if oriented compatibly) are the same. 
But $\partial A_{a,\uvec{e}}$ and $\partial A_{a,\scz{\ucev{e}}}$ are the boundary circles of the annulus $A_{a,e}$ with opposite orientations, therefore, the winding numbers of 
$$
u_a|_{\partial A_{a,\uvec{e}}}\quad \tn{and}\quad  u_a|_{\partial A_{a,\scz{\ucev{e}}}}
$$
are opposite of each other. If $v\!\in\!\V_1$, by the proof of Proposition~\ref{STNS1_prp}, 
$$
s_{\uvec{e}}:=\tn{ord}_{q_{\uvec{e}}} \ze_v= \tn{winding number of~}(u_a|_{\partial A_{a,\uvec{e}}})\qquad \forall~a\!>\!\!>\!1.
$$
Similarly, if $v\!\in\!\V_0$, then 
$$
s_{\uvec{e}}:=\tn{ord}_{q_{\uvec{e}}} (u_v,D)= \tn{winding number of~}(u_a|_{\partial A_{a,\uvec{e}}})\qquad \forall~a\!>\!\!>\!1.
$$
Therefore, 
\bEq{NodeContactCondition_e}
s_{\uvec{e}}=-s_{\scz\ucev{e}}\qquad \forall~v,v'\!\in\!\V ,~e\!\in\!\E_{v,v'}.
\eEq
The same conclusion holds in the case of Corollary~\ref{STNS1_cr} (since it is a corollary of Proposition~\ref{STNS1_prp}).
\vskip.1in

\noindent
The contact order condition~\ref{MatchingOrders_l} in Definition~\ref{PreLogMap_dfn}, for every $e\!\in\!\E(\Gamma'/\Gamma)$, follows from (\ref{NodeContactCondition_e}). For each $\uvec{e}\!\in\!\uvec{\E}\!=\!\uvec{\E}(\Gamma)\!\subset\!\uvec{\E'}\!=\!\uvec{\E}(\Gamma')$, with $v\!=\!v_1(\uvec{e})\!\in\!\V$, the nodal point $q_{\uvec{e}}$ is a marked point for $(u_v,C_v)$. For such $e$, by the last statements in Proposition~\ref{STNS1_prp} and Corollary~\ref{STNS1_cr}, the contact order $s_{\uvec{e}}$ remains unchanged in the limiting process. Therefore, the contact order condition~\ref{MatchingOrders_l} in Definition~\ref{PreLogMap_dfn}, for every $e\!\in\!\E\!\subset\! \E'$, follows from the corresponding condition on $(f_{a})_{a\in \N}$.\\

\noindent 
\textbf{Contact orders at the marked points.} Finally, Condition~\ref{OatMkd_l} in Definition~\ref{PreLogMap_dfn} follows from the corresponding statements in Proposition~\ref{STNS1_prp} and Corollary~\ref{STNS1_cr}.
\ePf

\noindent
 In order to prove Theorem~\ref{Log_thm11} (and thus Theorem~\ref{Compactness_th}), it just remains to prove the following proposition.

\noindent
\bPr{LGSmoothV_prp}
If further $(\om,\cR,J)\!\in\!\tn{AK}(X,D)$ for some regularization $\cR$ or if $J$ is integrable, then the pre-log $J$-holomorphic map $f'$ in (\ref{LogLimit_e}) satisfies Conditions \ref{Tropical_l} and \ref{GObs_e} of Definition~\ref{LogMap_dfn}.
\ePr

\noindent
We prove Proposition~\ref{LGSmoothV_prp} in Section~\ref{CptLogSmooth_ss}. The proof uses a fine comparison result between the rescaling parameters $(t_{a,v',i})_{a\in \N}$ corresponding to the sections  $\ze_{v',i}$, for all $v'\!\in\!\V'$ and $i\!\in\!I_{v'}$, and the ``gluing parameters" of the nodes. 
We expect Proposition~\ref{VertexOrder_prp} and thus Proposition~\ref{LGSmoothV_prp} to be true for a larger 
class of almost K\"ahler structures containing $\tn{AK}(X,D)$ and the space of K\"ahler structures.

%----------------------------------------------------
\subsection{Local behavior of convergence}\label{LocalTh_ss}

Proposition~\ref{LGSmoothV_prp} is essentially a consequence of  Proposition~\ref{VertexOrder_prp} below that relates the sequence of rescaling parameters $(t_{a,v',i})_{a\in \N}$ corresponding to the sections  $\{\ze_{v',i}\}_{i\in I_{v'}, v'\in \V'}$ in Lemma~\ref{PreLog_lmm} to the ``gluing parameters" at the nodes and the ratios of leading order coefficients $0\!\neq\!\eta_{\uvec{e}',i}\!\in \!\cN_XD_i|_{u'(q_{e'})}$ in (\ref{etae_e}). We use the natural log of these parameters to cook up the map required in Condition \ref{Tropical_l} of Definition~\ref{LogMap_dfn}.\\

\noindent
Let us start with a local picture of what is happening in Lemma~\ref{PreLog_lmm} with respect to any smooth component of $D$. Suppose $D$ is a smooth symplectic divisor in $(X,\om)$ and $J\!\in\!\cJ(X,D,\om)$. Fix a regularization $\Psi\colon\!\cN'_XD\!\lra\!X$ as in (\ref{RegPsi_e}). Let $\De_1$ and $\De_2$ be compact discs of some fixed sufficiently small  radius $\de$ around $0\!\in\!\C$ with coordinates $z_1$ and $z_2$. For $i\!=\!1,2$, let $\{z_{i,a}\}_{a\in \N}$ be a sequence of complex-coordinates\footnote{more precisely, $z_{i,a}\colon\!\De_i\!\lra\!\C$ is a sequence of smooth functions converging to the function $z_i\colon\!\De_i\!\lra\!\C$ in $C^\infty$-topology.} on $\De_i$ converging to $z_i$ uniformly with all the derivatives. \\

\noindent
\textbf{Local case 1}. For a sequence of complex numbers $(\ve_a)_{a\in \N}$ converging to zero, suppose $u_a\colon\!A_a\lra \tn{Im}(\Psi)\!\subset \!X$, where
\bEq{LocalSeq_e}
A_{a}=\{(z_{1,a},z_{2,a})\colon z_{1,a}z_{2,a}=\ve_{a}\colon z_{1,a}\in \De_1,z_{2,a}\in \De_2 \}\subset \De_1\times\De_2, \quad \forall~a\!\in\!\N,
\eEq
is a sequence of $J$-holomorphic maps that Gromov converges to the nodal map 
$$
\big(u_1(z_1)\colon\!\De_1 \!\lra\! X,~u_2(z_2)\colon\!\De_2\lra D),\quad x\!=\!u_1(0)\!=\!u_2(0)\!\in\!D.
$$ 
In other words, for any $\ep\!>\!0$, 
\bEnalph
\item\label{Glimit1_it} the sequence of $J$-holomorphic maps 
$$
u_a(z_{1,a})\!\equiv \!u_a(z_{1,a},\ve_a/z_{1,a})\colon A_a\approx \{z_{1,a}\in\C\colon  \ve_{a}/\de\leq |z_{1,a}|\!\leq\! \de\} \lra X
$$
converges uniformly with all the derivatives on the compact set 
$$
\{(z_{1,a},z_{2,a})\!\in\!A_a\colon \ep\!\leq\! |z_{1,a}| \}\approx \{z_{1,a}\!\in\!\C\colon \ep\!\leq\! |z_{1,a}|\!\leq\! \de \}
$$
to $u_1|_{\{z_{1}\in \C\colon \ep \leq  |z_{1}| \leq \de \}}$, 
\item\label{Glimit2_it} the reparametrization 
$$
u_a(z_{2,a})\!\equiv \!u_a(\ve_a/z_{2,a},z_{2,a})\colon A_a\approx \{z_{2,a}\in\C\colon  \ve_{a}/\de\leq |z_{2,a}|\!\leq\! \de\} \lra X
$$
converges uniformly with all the derivatives on the compact set 
$$
\{(z_{1,a},z_{2,a})\!\in\!A_a\colon \ep\!\leq\! |z_{2,a}| \}\approx \{z_{2,a}\!\in\!\C\colon \ep\!\leq\! |z_{2,a}|\!\leq\! \de \}
$$
to $u_2|_{\{z_{2}\in \C\colon \ep \leq  |z_{2}| \leq \de \}}$, and
\item\label{energy_it} we do not get any bubbling in between the two maps (i.e. the energy in between shrinks to zero with $\ep$).
\eEnalph
Furthermore, suppose 
\bEn
\item $u_1$ has a tangency order of $s\!>\!0$ with $D$ at $z_1\!=\!0$, and
\item\label{zeCondition_it} there exists a meromorphic section $\ze$ of $u_2^*\cN_XD$ with (only) a pole of order $s$ at the origin and a sequence of complex numbers 
$(t_a)_{a\in \N}$ converging to zero such that $t_a^{-1}\Psi^{-1}(u_a(z_{2,a}))$ converges to $\ze(z_2)$ uniformly with all the derivatives on any compact set $\{ z_2\!\in\!\C\colon \ep\leq  |z_2| \leq \de\} \subset \De_2$. 

\eEn
Let $0\!\neq\!\eta_2\!\in\!\cN_XD|_{x}$ be the leading coefficient of $\ze$ with respect to the coordinate $z_2$ as in (\ref{LocalCoord_e}),
and  $0\!\neq\!\eta_1\!\in\!\cN_XD|_{x}$ be the $s$-th derivative of $u_1$  in the normal direction to $D$ at $0$ with respect to the coordinate $z_1$ as in (\ref{normalDerivative}). Proposition~\ref{VertexOrder_prp} below shows that there is an explicit relation between the sequence of gluing parameters $(\ve_a)_{a\in N}$, the sequence of rescaling parameters $(t_a)_{a\in \N}$, and the ratio $\eta_2/\eta_1\!\in\!\C^*$.\\

\noindent
\textbf{Local case 2}. 
Similarly, consider the situation where  the sequence of $J$-holomorphic maps  $\{u_a\}_{a\in \N}$ in (\ref{LocalSeq_e}) Gromov converges to the nodal map 
$$
\big(u_1\colon\!\De_1 \!\lra\! D, u_2\colon\!\De_2\lra D),\quad x\!=\!u_1(0)\!=\!u_2(0)\!\in\!D,
$$ 
with the following property: there exist meromorphic sections $\ze_1(z_1)$ and $\ze_2(z_2)$ of $u_1^*\cN_XD$ and $u_2^*\cN_XD$, respectively, such that
$$
\tn{ord}_0(\ze_1)\!=\!s, \qquad \tn{ord}_0(\ze_2)\!=\!-s,
$$ 
and, for $i\!=\!1,2$, there exists a sequence of complex numbers  $(t_{i,a})_{a\in \N}$ converging to zero such that $t_{i,a}^{-1}\Psi^{-1}(u_a(z_{i,a}))$ converges to $\ze_i(z_i)$ uniformly with all the derivatives on any compact set $ \{ z_i\colon \ep\leq  |z_i| \leq \de\}\! \subset\! \De_i$.
With $\eta_1$ and $\eta_2$ as before, the following theorem also shows that there is a similar relation between the sequence of gluing parameters $(\ve_a)_{a\in \N}$, rescaling parameters $(t_{i,a})_{a\in \N}$, and the ratio $\eta_2/\eta_1\!\in\!\C^*$.

\bPr{VertexOrder_prp}
With notation as above, if further $(\om,\cR,J)\!\in\!\tn{AK}(X,D)$ for some regularization $\cR$ or if $J$ is integrable, in the local case 1, we have 
\bEq{ratio_e1Local}  
\lim_{a\lra \infty}\frac{\ve_{a}^s}{t_a} = \frac{\eta_{2}}{ \eta_{1}};
\eEq
in the local case 2, we have 
\bEq{ratio_e2Local}  
\lim_{a\lra \infty}\frac{t_{1,a}~\ve_{a}^s} {t_{2,a}} = \frac{\eta_{2}}{\eta_{1}} .
\eEq
\ePr

\noindent
Note that the situation in (\ref{ratio_e2Local}) reduces to the situation in (\ref{ratio_e1Local}) after a rescaling of the sequence $\{u_a\}_{a\in \N}$ via $(t_{1,a})_{a\in \N}$.  For the rescaled sequence we will have $\big(t_a=\frac{t_{2,a}}{t_{1,a}}\big)_{a\in \N}$. 
We prove Proposition~\ref{VertexOrder_prp} in the next section. The proof of this proposition is the only place where we use  the extra assumption on $J$ in the statement of Theorem~\ref{Compactness_th}, but we expect this proposition, and thus Theorem~\ref{Compactness_th}, to be true for a larger class of almost complex structures that contains $\cJ(X,D)$ and holomorphic structures.

\bRm{onSTNS1_rmk} 
It is easy to see that the limit conditions in (\ref{ratio_e1Local}) and (\ref{ratio_e2Local})  are independent of $\Psi$, the representatives $\ze_1$, $\ze_{2}$, and the local coordinates $z_1$ and $z_{2}$. For example, in (\ref{ratio_e2Local}), changing $z_{2}$ with $\al z_{2}$ and $\ze_2$ with $\beta\ze_2$, for some $\al,\beta\!\in\!\C^*$, changes $\eta_{2}$ on the right-hand side of (\ref{ratio_e2Local}) to $\al^{s}\beta\eta_{2}$, changes $\ve_{a}$ and $t_{2,a}$ on the left-hand side of (\ref{ratio_e2Local}) 
to $\al \ve_{a}$ and $\beta^{-1} t_{2,a}$, respectively, and has no effect on the other terms. Thus it affects both sides of (\ref{ratio_e2Local}) equally. It is also clear that (\ref{ratio_e1Local}) and (\ref{ratio_e2Local}) only depend on the equivalence classes $[(\ve_a)_{a\in \N}]$, $[(t_a)_{a\in \N}]$, $[(t_{1,a})_{a\in \N}]$, and $[(t_{2,a})_{a\in \N}]$.
\eRm

\begin{remark}
In the case of smooth divisors, a significantly simpler version of (\ref{ratio_e2Local}) suffices for proving Proposition~\ref{LGSmoothV_prp}. 
Instead of (\ref{ratio_e2Local}), in order to get the partial order in Lemma~\ref{Partial-Order_lmm}, we only need to prove that 
\bEq{Weaker-limits_e}
\lim_{a\lra \infty}\frac{t_{1,a}} {t_{2,a}} \!=\! \frac{\eta_{2}}{\eta_{1}}\qquad \tn{if}~~s\!=\!0, \quad\tn{and}\quad 
\lim_{a\lra \infty}\frac{t_{1,a}} {t_{2,a}} \!=\! \infty\qquad \tn{if}~~s\!>\!0.
\eEq
The equalities in (\ref{Weaker-limits_e}) can be proved without the extra restriction on $J$. Thus, if $D$ is smooth, Theorem~\ref{Compactness_th} holds for arbitrary $(\om,J)\!\in\!\cJ(X,D)$.
\end{remark}

\newtheorem*{proofof-VertexOrder_prp}{Proof of Proposition~\ref{VertexOrder_prp}}
\begin{proofof-VertexOrder_prp}
The proof below is by constructing a modified sequence of $J$-holomorphic maps in $\cN_XD$.\vskip.1in

\noindent
Let $(\cR\!=\!(\rho,\nabla,\Psi),\mfi_{\cN_XD},\dbar_{\cN_XD},J_{X,D})$ be as in the beginning of Section~\ref{CptLog_ss}. If $(\om,\cR,J)\!\in\!\tn{AK}(X,D)$, then $\Psi^*J\!=\!J_{X,D}$. If $J$ is holomorphic, we consider a holomorphic chart $(z_1,\ldots,z_n)$ around $x\!=\!u_1(0)\!=\!u_2(0)\!\in\!D$ such that $D\!=\!(z_1\!=\!0)$. Then, replacing the rescaling procedure in the proof below with holomorphic rescaling of $z_1$, the same proof works for the holomorphic case.  
\\

\noindent
Assume $\Psi^*J\!=\!J_{X,D}$. Note that $J_{X,D}$ is $\C^*$-invariant. Since the argument is local, in order to simplify the notation, let us forget about $\Psi$ and think of $\{u_a\}_{a\in \N}$ as a sequence of $J_{X,D}$-holomorphic maps into $\cN'_XD$ itself. \\

\noindent
Assume that we are in the situation of local case 1. For each $a\!\in\!\N$, let 
$$
\la_a\!=\!\frac{\ve_a^s}{t_a}.
$$

\noindent
(\textbf{Claim 1}) First, we show that there is NO subsequence $(a_1,a_2,\ldots)$ of $\N$ such that 
$$
\lim_{i \lra \infty} \la_{a} = 0~\tn{or}~\infty.
$$
Thus, we conclude that there is $M\!>\!0$ such that $M^{-1}\!<\! \abs{\la_a}\!<\!M$ for all $a\!\in\!\N$.\\

\noindent
(\textbf{Claim 2}) Then, for any subsequence $(a_1,a_2,\ldots)$ such that the limit 
$$
\lim_{i \lra \infty} \la_{a}=\la
$$
exists, we show that $\la\!=\!\eta_2/\eta_1$. This implies that (\ref{ratio_e1Local}) holds over entire $\N$.\\

\noindent
In order to prove these claims, first, we construct two new sequences of $J$-holomorphic maps.
For $a\!\in\!\N$, define 
\begin{gather}
\label{u1a_e}
u_{1,a}\colon A_a\lra \cN_XD,\qquad u_{1,a}(z_{1,a},z_{2,a})=z_{1,a}^{-s}\,u_a(z_{1,a},z_{2,a}),\\
\label{u2a_e}
u_{2,a}\colon A_a\lra \cN_XD,\qquad u_{2,a}(z_{1,a},z_{2,a})=z_{1,a}^{-s}\, \la_a\,  u_a(z_{1,a},z_{2,a}),
\end{gather}
where the multiplications on the right-hand sides are with respect to the complex structure $\mfi_{\cN_XD}$ on $\cN_XD$. 
By \ref{HoloProj_l}-\ref{HoloSec_l} in Page~\pageref{HoloProj_l}, both (\ref{u1a_e}) and (\ref{u2a_e}) are sequences of $J_{X,D}$-holomorphic maps in $\cN_XD$.\\

\noindent
We will also use the following fact.
For any $c\!>\!0$, there exists a sufficiently small $\ep_c\!>\!0$ such that 
$$
\om_{\ep_c}= \pi^*(\om|_D) + \frac{\ep_c}{2}\nd(\rho\al_\nabla)
$$
tames $J_{X,D}$ on $\ov{\cN_XD(c)}$.
For any compact $2$-dimensional domain $\Si$ and a smooth map $u\colon\!\Si\!\lra\! \ov{\cN_XD(c)}$, 
let 
$$
\om_{\ep_c}(u)\!=\!\int_{\Si} u^*\om_{\ep_c}
$$
denote the symplectic area of $u$.\\

\noindent
In order to prove Claim 1, we separate the problem into two cases. In the first and second parts below, we consider the cases where the limit is $\infty$ and zero, respectively. In each case, we apply the Gromov convergence to the auxiliary sequences in (\ref{u1a_e}) and (\ref{u2a_e}) to get a contradiction if  the limit is $\infty$ or $0$.  \\

\noindent
\textbf{Proof of Claim 1 part 1}. 
After passing to a subsequence, suppose 
\bEq{limitinfty_case}
\lim_{a \lra \infty} \la_{a} = \infty.
\eEq
By assumption \ref{Glimit1_it} in Page~\pageref{Glimit1_it} and the previous paragraph, for any $0\!<\!r\!<\!\de$, restricted to  $r\!\leq\!|z_1|\!\leq \!\de$ (and its pre-images in $A_a$), the sequence $\{u_{1,a}(z_{1,a})\}_{a\in \N}$ converges uniformly with all the derivatives to the $J_{X,D}$-holomorphic map
$$
u_{1,\infty,1}(z_1)=z_{1}^{-s}\, u_1(z_1).
$$
By definition of $\eta_1$, the function $u_{1,\infty,1}(z_1)$ extends to $z_1\!=\!0$ with $u_{1,\infty,1}(0)\!=\!\eta_1 \!\in\!\cN_XD|_x$, where $x\!=\!u_1(0)\!=\!u_2(0)\!\in\!D$. 
By assumptions  \ref{Glimit2_it} and \ref{zeCondition_it} in Page~\pageref{zeCondition_it}, (\ref{limitinfty_case}), and since 
$$
z_{1,a}^{-s}= \ve_a^s \,z_{2,a}^s,
$$
restricted to  $r\!\leq\!|z_2|\!\leq \!\de$ (and its pre-images in $A_a$), the sequence $\{u_{1,a}(z_{2,a})\}_{a\in \N}$ converges uniformly with all the derivatives to the $J_{X,D}$-holomorphic map
$$
u_{1,\infty,2}(z_2)=u_2(z_2) \!\subset\!D.
$$
This obviously extends to entire $\De_2$ with $u_{1,\infty,2}(0)\!=\!x$. The following sub-claim shows that the sequence $\{u_{1,a}\}_{a\in \N}$ is bounded in between, so that Gromov's convergence applies.\\

\noindent
\textbf{Sub-Claim.} There exists a sufficiently large $c\!>\!0$ such that  
\bEq{Ctl-Im-En_e}
\tn{Im}\big(u_{1,a}\big) \subset \ov{\cN_XD(c)}\quad \tn{and}\quad
\om_{c_\ep}\big(u_{1,a}\big)\leq c\qquad \forall~a\!\in\!\N.
\eEq

\noindent
\textbf{Proof of Sub-Claim.} Suppose (\ref{Ctl-Im-En_e}) does not hold. Then (after passing to a subsequence), by assumptions  \ref{Glimit1_it}-\ref{energy_it} in Page~\pageref{zeCondition_it}, for any $c\!>\!1$, there exists a sequence $\{r_a\}_{a\in \N}$, with 
\bEq{r_alimit_e}
\lim_{a\lra \infty} r_a\!=\!\infty\qquad \tn{and}
\eEq
\bEq{max_e}
\max_{(z_{1,a},z_{2,a})\in A_a}\bigg(\abs{r_a^{-1}z_{1,a}^{-s}\, u_a(z_{1,a},z_{2,a})}, ~
\om_{\ep_c}\!\lrp{(z_{1,a},r_a^{-1}z_{1,a}^{-s}\, u_a(z_{1,a},z_{2,a}))}
\bigg)\!=\!c\quad \forall~a\!>\!\!>\!1.
\eEq
Let  
$$
\wt{u}_{1,a}\colon A_{a}\lra \cZ,\qquad u_{1,a}(z_{1,a},z_{2,a})=r_a^{-1}z_{1,a}^{-s}\, u_{a}(z_{1,a},z_{2,a}).
$$
Then 
\bIt
\item by assumption \ref{Glimit1_it} in Page~\pageref{Glimit1_it} and (\ref{r_alimit_e}), for any $0\!<\!r\!<\!\de$, restricted to  $r\!\leq\!|z_1|\!\leq \!\de$, the rescaled sequence $\{\wt{u}_{1,a}(z_{1,a})\}_{a\in \N}$ converges uniformly with all the derivatives to the $J_{X,D}$-holomorphic map
$$
\wt{u}_{1,\infty,1}(z_1)=\ov{u}_1(z_1)\!\subset\!D,
$$
where $\ov{u}_1$ is the image of $u_1$ in $D$;
\item by assumptions  \ref{Glimit2_it} and \ref{zeCondition_it} in Page~\pageref{zeCondition_it}, restricted to  $r\!\leq\!|z_2|\!\leq \!\de$, the sequence $\{\wt{u}_{1,a}(z_{2,a})\}_{a\in \N}$ still converges uniformly with all the derivatives to the $J_{X,D}$-holomorphic map
$$
\wt{u}_{1,\infty,2}(z_2)=u_2(z_2) \!\subset\! D; 
$$
\item and, by (\ref{max_e}), (the proof of) Gromov convergence theorem in \cite{MS2} applies\footnote{the Gromov convergence applies, because on the open ends of $A_a$ we already know that $\{\wt{u}_{1,a}\}_{a\in \N}$ uniformly converges to $\wt{u}_{1,\infty,1}$ and $\wt{u}_{1,\infty,2}$, and in the middle the sequence is bounded with bounded energy.} to the sequence $\{\wt{u}_{1,a}\}_{a\in \N}$. In the limit we get a bubble domain $\Si_\infty$ with $\De_1$ and $\De_2$ at the two ends and at least one  closed bubble in between (because of (\ref{max_e})), and a continuous $J_{X,D}$-holomorphic map 
$$
\wt{u}_{1,\infty}\colon \Si_\infty \lra  \cZ 
$$ 
such that 
$$
\wt{u}_{1,\infty}|_{\De_1}=\wt{u}_{1,\infty,1}\qquad\tn{and}\qquad \wt{u}_{1,\infty}|_{\De_2}=\wt{u}_{1,\infty,2}.
$$
\eIt
Any non-trivial bubble would have  trivial image in $D$, thus its image lives in 
$\ov{\cN_{X}D(c)}|_{x}$.
This is impossible since the latter is open and there are no-marked points to stabilize such a bubble.\qed\\

\noindent
Going back to the proof of Claim 1-Part 1, by (\ref{Ctl-Im-En_e}), (the proof of) Gromov convergence theorem in \cite{MS2} applies to the sequence $\{u_{1,a}\}_{a\in \N}$. In the limit we get a bubble domain $\Si_\infty$ with $\De_1$ and $\De_2$ at the two ends and possibly some closed bubbles in between, and a continuous $J_{X,D}$-holomorphic map 
$$
{u}_{1,\infty}\colon \Si_\infty \lra  \cZ 
$$ 
such that 
$$
{u}_{1,\infty}|_{\De_1}={u}_{1,\infty,1}\qquad\tn{and}\qquad {u}_{1,\infty}|_{\De_2}={u}_{1,\infty,2}.
$$
Since 
$$
{u}_{1,\infty,1}(0)\!\neq\!{u}_{1,\infty,2}(0),
$$
$\Si_\infty$ should include at least one non-trivial bubble. Such a non-trivial bubble would have  trivial image in $D$, thus its image lives in 
$\ov{\cN_{X}D(c)}|_{x}$.
This is impossible since the latter is a domain in $\C$ and there are no marked points to stabilize such a bubble.\\\qed

\noindent
\textbf{Proof of Claim 1 part 2}. 
After passing to a subsequence, suppose 
\bEq{iftozero_e}
\lim_{a\lra \infty} \la_{a} = 0.
\eEq
By assumptions  \ref{Glimit2_it} and \ref{zeCondition_it} in Page~\pageref{zeCondition_it}, since 
$$
u_{2,a}(z_{1,a},z_{2,a})=z_{1,a}^{-s}\, \la_a\,  u_a(z_{1,a},z_{2,a})=z_{2,a}^{s}\, t_a^{-1}\, u_a(z_{1,a},z_{2,a}),
$$
for any $0\!<\!r\!<\!\de$, restricted to  $r\!\leq\!|z_2|\!\leq \!\de$ (and its pre-images in $A_a$), the sequence $\{u_{2,a}(z_{2,a})\}_{a\in \N}$ converges uniformly with all the derivatives to the $J_{X,D}$-holomorphic map
$$
u_{2,\infty,2}(z_2)=z_{2}^{s}\,\ze(z_2).
$$
By definition of $\eta_2$, the function $u_{2,\infty,2}(z_2)$ extends to $z_2\!=\!0$ with $u_{2,\infty,2}(0)\!=\!\eta_2$. 
On the other hand,  by assumption  \ref{Glimit1_it} in Page~\pageref{Glimit1_it} and (\ref{iftozero_e}),
restricted to  $r\!\leq\!|z_1|\!\leq \!\de$ (and its pre-images in $A_a$), the sequence $\{u_{2,a}(z_{1,a})\}_{a\in \N}$ converges uniformly with all the derivatives to the $J_{X,D}$-holomorphic map
$$
u_{2,\infty,1}(z_1)=\ov{u_1}(z_1)\!\subset\!D.
$$
This obviously extends to the entire $\De_1$ with $u_{2,\infty,1}(0)\!=\!x$. By a similar argument as in the previous case, the inequality 
$$
u_{2,\infty,1}(0)\!\neq\! u_{2,\infty,2}(0),
$$
leads to a contradiction. This finishes the proof of Claim 1.\qed\vskip.2in 

\noindent
\textbf{Proof of Claim 2}. After passing to a subsequence, suppose 
$$
\lim_{a \lra \infty}\la_a = \la\neq 0.
$$
Then, going back to the proof of Claim 1 part 1, since
$$
z_{1,a}^{-s}= \ve_a^s \,z_{2,a}^s\qquad \forall~a\!\in\!\N,
$$
restricted to  $r\!\leq\!|z_2|\!\leq \!\de$, the sequence  $\{u_{1,a}(z_{2,a})\}_{a\in \N}$ converges uniformly with all the derivatives to the $J_{X,D}$-holomorphic map
$$
u_{1,\infty,2}(z_2)=\la^{-1}\, z_2^{s}\,\ze_2(z_2).
$$
This extends to the entire $\De_2$ with $u_{1,\infty,2}(0)\!=\!\la\eta_2$.  By a similar argument as in the proof of Claim~1-part~1, if
$$
u_{1,\infty,1}(0)\!\neq\! u_{1,\infty,2}(0),
$$
we get a contradiction. Therefore, 
$$
\eta_1=u_{1,\infty,1}(0)\!=\! u_{1,\infty,2}(0)=\la^{-1}\eta_2;
$$
in other words, $\la\!=\!\eta_2/\eta_1$.
This finishes the proof of Proposition~\ref{VertexOrder_prp} in the local case 1.\\

\noindent
For the local case 2, repeat the exact same proof with 
$$
\aligned
&u_{1,a}\colon A_a\lra \cN_XD,\qquad u_{1,a}(z_{1,a},z_{2,a})=z_{1,a}^{-s}\,t_{1,a}^{-1}\, u_a(z_{1,a},z_{2,a})
\quad\tn{and}\\
&u_{2,a}\colon A_a\lra  \cN_XD,\qquad u_{2,a}(z_{1,a},z_{2,a})=z_{1,a}^{-s}\, \la_a\, u_a(z_{1,a}),
\endaligned
$$
in place of (\ref{u1a_e}) and (\ref{u1a_e}), respectively, where 
$$
\la_a\!=\!\frac{t_{1,a}\,\ve_a^s}{t_{2,a}}\qquad \forall~a\!\in\!\N.
$$
This finishes the proof of Proposition~\ref{VertexOrder_prp} under the assumption that $\Psi^*J\!=\!J_{X,D}$. \qed
\end{proofof-VertexOrder_prp}

\noindent
\begin{remark}
For arbitrary $J$ on $\cN'_XD$, define
$$
\cZ\!=\!\lrc{(t,v)\!\in\!\C\!\times\!\cN_XD\colon ~t^sv\!\in\!\cN'_XD}, \quad \cZ_*\!=\!\lrc{(t,v)\!\in\!\cZ\colon~ t\!\in\!\C^*},
$$ 
and 
$$
F\colon\!\cZ_*\!\lra\!\C\!\times\!\cN'_XD, \qquad F(t,v)=(t, t^sv).
$$
Let $J_{\cZ}=F^*(\mfi\!\times\!J)$, where $\mfi$ is the standard almost complex structure on $\C$ and $\mfi\!\times\!J$ is the product almost complex structure on the target. 
By an argument similar to Lemma~\ref{JLemma_e}, the almost complex structure $J_{\cZ}$ on $\cZ_*$ extends to a (similarly denoted) almost complex structure on the entire $\cZ$ satisfying 
\bEq{limitJZ_e}
J_{\cZ}|_{\{0\}\times\cN_XD \,\cup\, \C\times D}  \cong \mfi\!\times\! J_{X,D}.
\eEq
Similarly, for every $a\!\in\!\N$, let
$$
\wt\cZ_{a}\!=\!\lrc{(t,v)\!\in\!\C\!\times\!\cN_XD\colon~t^{s}\la_a^{-1}v\!\in\!\cN'_XD},\quad \wt\cZ_{a,*}=\{(t,v)\!\in\!\wt\cZ_a\colon\!t\!\in\!\C^*\},
$$ 
and define
\bEq{F_anew}
F_a\colon\!\wt\cZ_{a,*}\!\lra\!\C\!\times\!\cN'_XD, \qquad F_a(t,v)=(t, \la_a^{-1}t^{s}v).
\eEq
For each $a\!\in\!\N$, let $J_{a}\!=\!F_a^*(\mfi\!\times\!J)$. By Lemma~\ref{JLemma_e} and the previous paragraph, for each $a\!\in\!\N$, the almost complex structure $J_{a}$ on $\wt\cZ_{a,*}$ extends to a (similarly denoted) almost complex structure on the entire $\wt\cZ_a$ satisfying (\ref{limitJZ_e}).\\

\noindent
For $a\!\in\!\N$, define 
\begin{gather}
\label{u1a_General}
u_{1,a}\colon A_a\lra \cZ,\qquad u_{1,a}(z_{1,a},z_{2,a})=\big(z_{1,a},z_{1,a}^{-s}\,u_a(z_{1,a},z_{2,a})\big),\\
\label{u2a_eGeneral}
u_{2,a}\colon A_a\lra \cZ_a,\qquad u_{2,a}(z_{1,a},z_{2,a})=\big(z_{1,a},z_{1,a}^{-s}\, \la_a\,  u_a(z_{1,a},z_{2,a})\big).
\end{gather}
By definition, (\ref{u1a_e}) is a sequence of $J_\cZ$-holomorphic maps in $\cZ$ and (\ref{u2a_e}) is a sequence of $J_a$-holomorphic maps in $\cZ_a$. In principle, one may try the proof above by replacing  (\ref{u1a_e}) and (\ref{u1a_e}) with (\ref{u1a_General}) and (\ref{u1a_General}), respectively.  However, multiplication by $\la_a^{-1}$ in (\ref{F_anew}) and by $r_a^{-1}$ in (\ref{max_e}) have adverse effects on the almost complex structure, making it hard to apply the Gromov convergence. 
\end{remark}

%----------------------------------------------------
\subsection{Proof of Proposition~\ref{LGSmoothV_prp} and~Theorem~\ref{Compactness_th}}\label{CptLogSmooth_ss}

Going back to the set up of Proposition~\ref{LGSmoothV_prp}, first, assume that the dual graph $\Gamma$ of $f_a$ in (\ref{faSequence_e}) is made of only one vertex $\V\!=\!\{v\}$ (in other words, restrict to the $v$-th component of the sequence $(f_a)_{a\in \N}$ in (\ref{faSequence_e})) and fix a set of representatives 
$$
(\ze_{a,v,i})_{i\in I_v}
$$
for $[\ze_{a,v}]$. 
For each $v'\!\in\!\V'$ and $i\!\in\!I_{v'}$ fix a representative $\ze_{v',i}$ of the $\C^*$-equivalence class $[\ze_{v',i}]$ in Lemma~\ref{PreLog_lmm}, and a sequence of rescaling parameters $(t_{a,v',i})_{a\in \N}$ satisfying Proposition~\ref{STNS1_prp} or Corollary~\ref{STNS1_cr}, depending on whether $i\!\notin\!I_v$ or $i\!\in\! I_v$, respectively. \\

\noindent
By the surjectivity of the classical gluing theorem of $J$-holomorphic maps (e.g. \cite[Sec~7]{FOOO}), for $a$ sufficiently large, the domain $\Si_a\cong \Si$ of (the stable map underlying) $f_a$ can be obtained from the nodal domain $\Si'$ of (the stable map underlying) $f$ in the following way.~There exist
\bIt
\item a sequence of complex structure $\mfj'_a=(\mfj_{v',a})_{v'\in \V'}$ on the nodal domain $\Si'\!=\!(\Si_{v'})_{v'\in\V'}$ of the stable nodal map $\iota(f)$ in (\ref{Stable-limit_e}),
\item a sequence of local $\mfj_{v',a}$-holomorphic coordinates $z_{\uvec{e}',a}\colon \De_{\uvec{e}'}\!\lra\! \C$ around $q_{\uvec{e}'}\!\in\!\Si_{v'}$, for all $v'\!\in\!\V'$ and $\uvec{e}'\!\in\!\uvec{\E}'_{v'}$, and
\item  a sequence of non-zero complex numbers $(\ve_{e',a})_{e'\in \E'}$ converging to zero,
\eIt
such that 
\bEn
\item\label{glue-domain_it} $(\Si_a,\mfj_a,\vec{z_a})$ is isomorphic to the smoothing of $(\Si', \mfj'_a=(\mfj_{v',a})_{v'\in \V'})$ defined by
\bEq{GluignRelation_e}
z_{\uvec{e}',a}z_{\scz\ucev{e}',a}=\ve_{e',a} \qquad \forall~e'\!\in\!\E', \quad \tn{and}
\eEq
\item the sequence $(\mfj_{v',a})_{a\in \N}$ $C^\infty$-converges to $\mfj_{v'}$ for all $v'\!\in\!\V'$,
\item\label{limitze_e}  the sequence  $(z_{\uvec{e}',a})_{a\in \N}$ $C^\infty$-converges to $z_{\uvec{e}'}$, where $z_{\uvec{e}'}\colon \De_{\uvec{e}'}\!\lra\! \C$ is some fixed local $\mfj_{v'}$-holomorphic coordinate  around $q_{\uvec{e}'}\!\in\!\Si_{v'}$, for all $v'\!\in\!\V'$ and $\uvec{e}'\!\in\!\uvec{\E}'_{v'}$.
\eEn
We will use this standard presentation of $(\Si_a,\mfj_a)$ in the proof of Proposition~\ref{LGSmoothV_prp} and Theorem~\ref{Compactness_th}.

\bRm{STDpres_rmk}
For $\de\!>\!0$ sufficiently small, let 
$$
\aligned
 &\De_{\uvec{e}',a}(\de)=\{x\!\in\!\De_{\uvec{e}'}\colon~ |z_{\uvec{e}',a}(x)|\!<\!\de\}\qquad \forall~\uvec{e}'\!\in\!\E',~a\!>\!\!>\!1,\quad \tn{and}\\
 &A_{e',a}=\{ z_{\uvec{e}',a}z_{\scz\ucev{e}',a}=\ve_{e',a}\colon z_{\uvec{e}',a}\in \De_{\uvec{e}'}(2\ve_{e',a}),z_{\scz\ucev{e}',a}\in \De_{\scz\ucev{e}'}(2\ve_{e',a})  \} \subset \Si_a\quad \forall~e'\!\in\!\E',~a\!>\!\!>\!1.
 \endaligned
$$
Then, with respect to the identification of the domains in \ref{glue-domain_it}, the $\gamma_a$-degeneration maps 
$$
\varphi_{\gamma_a}\colon\!\Si_a\!\lra\!\Si',
$$
can be taken to be identity on the complement of  $\cup_{e'\in \E'} A_{e',a}$ and some ``nice" degeneration map
$$
A_{e',a} \lra  \De_{\uvec{e}'}(2\ve_{e',a})\cup \De_{\scz\ucev{e}'}(2\ve_{e',a})
$$
on the neck region. 
\eRm

\noindent
For each $\uvec{e}'\!\in\!\uvec{\E}'$ and $i\!\in\!I_{e'}$, let 
$$
0\!\neq\!\eta_{\uvec{e}',i}\!\in \!\cN_XD_i|_{u'(q_{e'})}
$$ 
be the leading coefficient term in (\ref{etae_e}) with respect to $z_{\uvec{e}'}$ (and $\ze_{v',i}$, if $i\!\in\!I_{v'}$). By Proposition~\ref{VertexOrder_prp}, for every $v'_1,v'_2\!\in\!\V'$ and $\uvec{e}'\!\in\!\uvec{\E}'_{v'_1,v'_2}$ we have 
\begin{gather}
\label{ratio_e2}
\lim_{a\lra \infty}\frac{t_{a,v'_1,i}~\ve_{e',a}^{s_{\uvec{e}',i}}}{t_{a,v'_2,i}} = \eta_{\scz\ucev{e}',i}/ \eta_{\uvec{e}',i} \qquad \forall~i\!\in\!I_{v_1'}\!\cap\!I_{v_2'},\\
\label{ratio_e1}  \lim_{a\lra \infty}t_{a,v'_1,i}~\ve_{e',a}^{s_{\uvec{e}',i}} = \eta_{\scz\ucev{e}',i}/ \eta_{\uvec{e}',i}\qquad \forall~i\!\in\!I_{v_1'}\!-\!I_{v_2'}.
\end{gather}

\noindent
The following proposition shows that, for $a$ sufficiently large, we can adjust the choices involved to get equality at each $a$.

\bPr{VertexOrder_Pr}
There exists a choice of the coordinates $\{z_{\uvec{e}'}\}_{\uvec{e}'\in \uvec{\E}'}$ and $\{z_{\uvec{e}',a}\}_{\uvec{e}'\in \uvec{\E}',a\in \N}$ satisfying (\ref{GluignRelation_e}) and item~\ref{limitze_e} after that, and the representatives $\ze_{v',i}$ and $(t_{a,v',i})_{a\in \N}$ for $[\ze_{v',i}]$ and $[(t_{a,v',i})_{a\in \N}]$, respectively, such that 
\begin{gather}
\label{ratio-equal_e2}
t_{a,v'_1,i}~\ve_{e',a}^{s_{\uvec{e}',i}} = t_{a,v'_2,i} \qquad \forall~i\!\in\!I_{v_1'}\!\cap\!I_{v_2'},~a\!>\!\!>\!1,\\
\label{ratio-equal_e1} 
t_{a,v'_1,i}~\ve_{e',a}^{s_{\uvec{e}',i}} = 1 \qquad \forall~i\!\in\!I_{v_1'}\!-\!I_{v_2'},~a\!>\!\!>\!1.
\end{gather}
\ePr
\vskip.1in

\noindent
Proof of Proposition~\ref{VertexOrder_Pr}  uses the following lemma with the linear map 
$$
\vr_\C\colon\!\C^{\E'}\oplus \bigoplus_{v'\in \V'} \C^{I_{v'}}\lra \bigoplus_{e'\in \E'} \C^{I_{e'}}
$$
defined in (\ref{DtoT_e}). We will use Proposition~\ref{VertexOrder_Pr} to construct a map 
$$
s \colon\! \V'\!\lra\!\R^N, \quad v'\!\lra\!s_{v'},\qquad\tn{and}\qquad \la \colon\! \E'\!\lra\!\R_+, \quad e'\!\lra\!\la_{e'},
$$ 
satisfying Condition \ref{Tropical_l} and also to show that the limit satisfies  Condition~\ref{GObs_e} of Definition~\ref{LogMap_dfn}

\bLm{Adjust_lmm}
Assume $f\colon\!\C^n\!\lra\!\C^m$ is a complex-linear map and 
$
(\xi_a)_{a\in \N}\!\subset\!\C^n
$ 
is a sequence such that 
\bEq{Limitxi_e}
\lim_{a\lra \infty} f(\xi_a)\!=\!\eta.
\eEq
Then, there exists a convergent sequence $(\xi'_a)_{a\in \N}\!\subset\!\C^n$ (i.e. $\exists~\xi'\!\in\!\C^n$ such that $\lim_{a\lra \infty}\xi'_a\!=\!\xi'$)  such that $f(\xi_a-\xi'_a)\!=\!0$ for all $a\in \N$.
\eLm

\bPf
Since $\tn{Im}(f)\!\subset\!\C^m$ is closed, (\ref{Limitxi_e}) implies that $\eta\!\in\!\tn{Im}(f)$. Let $\eta\!=\!f(\xi)$. Fix an affine subspace\footnote{i.e. a shifted linear subspace.} $H\!\subset\!\C^n$ passing through $\xi$ and transverse\footnote{assuming $f$ in not trivial; otherwise, the Lemma is obvious.} to the hyperplane $f^{-1}(\eta)\!\subset\!\C^n$.
By (\ref{Limitxi_e}), there exists $M\!\in\!\N$ such that $H$ is transverse to $f^{-1}(f(\xi_a))$ for all $a\!>\!M$.
Then the sequence $\xi'_a\!=\! f^{-1}(f(\xi_a))\!\cap\! H$, if $a\!>\!M$, and $\xi'_a\!=\!\xi_a$, if $a\!\leq\!M$, has the desired properties.
\ePf

\newtheorem*{proofof-VertexOrder_Pr}{Proof of Proposition~\ref{VertexOrder_Pr}}
\begin{proofof-VertexOrder_Pr}
Throughout the proof we assume $I_v\!=\!\eset$; for $I_v\!\neq\!\eset$, the argument reduces to $I_v\!=\!\eset$ by considering the associated sequence of maps in $\cN_XD_{I_v}$. We modify a given set of representatives to another set satisfying (\ref{ratio-equal_e2}) and (\ref{ratio-equal_e1}). Assuming $I_v\!=\!\eset$, fix an orientation $O$ on $\E'$, and choose some branch 
$$
\eta=\bigoplus_{\uvec{e}'\in O} \eta_{e'}\!\in\! \bigoplus_{e'\in \E'} \C^{I_{e'}},\qquad \eta_{e'}=\big(-\log (\eta_{\scz\ucev{e}',i}/ \eta_{\uvec{e}',i})\big)_{i\in I_{e'}}\!\in\! \C^{I_{e'}}\quad \forall~\uvec{e}'\!\in \!O,
$$
of the multi-valued function $\log$. By (\ref{ratio_e2}) and (\ref{ratio_e1}) and definition of $\vr_\C$ in (\ref{DtoT_e}) (via the chosen orientation $O$), we can choose the branches
$$
\xi_{a}=\big((-\log(\ve_{e',a}))_{e'\in \E'},(-\log(t_{a,v',i}))_{v'\in \V', i\in I_{v'}}  \big)\in \C^{\E'} \oplus \bigoplus_{v'\in \V'} \C^{I_{v'}}\qquad \forall~a\!\in\!\N.
$$
so that 
$$
\lim_{a\lra \infty} \vr_\C(\xi_{a}) =\eta.
$$
By Lemma~\ref{Adjust_lmm} applied to $\vr_\C$, there exists a sequence 
$$
(\xi'_a)_{a\in \N}\!\subset\!\C^{\E'} \oplus \bigoplus_{v'\in \V'} \C^{I_{v'}}
$$ 
such that $\vr_\C(\xi_a-\xi'_a)\!=\!0$ for all $a\in \N$ and $\lim_{a\lra \infty}\xi'_a\!=\!\xi'$.
Taking the exponential of $\xi'_a$ and $\xi'$, we conclude that there exist
$$
\big((\al_{e'})_{e'\in \E'}, (\al_{v',i})_{v'\in\V',i\in I_{v'}}\big),\big((\al_{e',a})_{e'\in \E'}, (\al_{a,v',i})_{v'\in\V',i\in I_{v'}}\big)_{a\in \N}\in (\C^*)^{\E'}\times \prod_{v'\in \V'} (\C^*)^{I_{v'}}
$$
such that 
$$
\lim_{a\lra\infty} \big((\al_{e',a})_{e'\in \E'}, (\al_{a,v',i})_{v'\in\V',i\in I_{v'}}\big) =\big((\al_{e'})_{e'\in \E'}, (\al_{v',i})_{v'\in\V',i\in I_{v'}}\big)
$$
and
\begin{gather}
\label{ratio_e2adjusted}
\frac{\big(\al^{-1}_{v'_1,i} t_{a,v'_1,i}\big)~\big(\al_{e',a}^{-1}\ve_{e',a}\big)^{s_{\uvec{e}',i}}}
{\big(\al^{-1}_{v'_2,i}t_{a,v'_2,i}\big)} =1 \qquad \forall~i\!\in\!I_{v_1'}\!\cap\!I_{v_2'},~a\!\in\!\N
,\\
\label{ratio_e1adjusted}  
\big(\al^{-1}_{v'_1,i} t_{a,v'_1,i}\big)~\big(\al_{e',a}^{-1}\ve_{e',a}\big)^{s_{\uvec{e}',i}}= 1 \qquad \forall~i\!\in\!I_{v_1'}\!-\!I_{v_2'},~a\!\in\!\N.
\end{gather}
By (\ref{ratio_e2adjusted}) and (\ref{ratio_e1adjusted}), for $a$ sufficiently large, replacing 
\bIt
\item $\{z_{\uvec{e}'}\}_{\uvec{e}'\in O}$ with $\{\al_{e'}^{-1} z_{\uvec{e}'}\}_{\uvec{e}'\in O}$, 
\item  $\{z_{\uvec{e}',a}\}_{\uvec{e}'\in O}$ with $\{\al^{-1}_{e',a} z_{\uvec{e}',a}\}_{\uvec{e}'\in O}$,
\item  $\{\ve_{e',a}\}_{e'\in \E'}$ with $\{\al^{-1}_{e',a} \ve_{e',a}\}_{e'\in \E'}$,
\item $(t_{a,v',i})_{v'\in \V', i\in I_{v'}}$ with $(\al^{-1}_{a,v',i}t_{a,v',i})_{v'\in \V', i\in I_{v'}}$, and
\item  $(\ze_{v',i})_{v'\in \V', i\in I_{v'}}$ with $(\al_{v',i}\ze_{v',i})_{v'\in \V', i\in I_{v'}}$, 
\eIt
we get a new set of representatives satisfying (\ref{ratio-equal_e2}) and (\ref{ratio-equal_e1}). In particular, the limits in (\ref{ratio_e2}) and (\ref{ratio_e1}) can be set to be equal to $1$.\qed
\end{proofof-VertexOrder_Pr}

\newtheorem*{proofof-LGSmoothV_Th}{Proof of Proposition~\ref{LGSmoothV_prp}}
\begin{proofof-LGSmoothV_Th}
First, assume that the dual graph $\Gamma$ of $f_a$ in (\ref{faSequence_e}) is made of only one vertex $\V\!=\!\{v\}$ and fix a set of representatives 
$$
(\ze_{a,v,i})_{i\in I_v}
$$
for $[\ze_{a,v}]$.
By Propositions~\ref{VertexOrder_prp} and \ref{VertexOrder_Pr}, we can choose the coordinates $\{z_{\uvec{e}'}\}_{\uvec{e}'\in \uvec{\E}'}$ and $\{z_{\uvec{e}',a}\}_{\uvec{e}'\in \uvec{\E}',a\in \N}$, and the representatives $\ze_{v',i}$ and $(t_{a,v',i})_{a\in \N}$ so that
(\ref{ratio-equal_e2}) and (\ref{ratio-equal_e1}) hold. For each $v'\!\in\!\V'$ and $i\!\in\!I_{v'}\!-\!I_v$, note that  $(t_{a,v',i})_{a\in \N}$ converges to $0$; therefore,
$$
-\log(|t_{a,v',i}|)>0 \qquad \forall~v'\!\in\!\V,~i\!\in\!I_{v'}\!-\!I_v,~a\!>\!\!>1,
$$
and it converges to infinity.
Choose a sequence of positive vectors $s^{a}_{v}\!=\!(s^a_{v,i})_{i\in I_v}\!\in\!\R_+^{I_v}$  such that 
\bEq{LargeEnough_e}
s^a_{v,i}-\log(|t_{a,v',i}|)> 1 \qquad \forall~v'\!\in\!\V',~i\!\in\!I_v.
\eEq
With these choices, for $a\!>\!\!>1$, the functions $s^a\colon \V'\lra \R^N$ defined by 
\bEq{sv-from-log_e}
s^a_{v'}=\bigg(  \big( s^a_{v,i}-\log(|t_{a,v',i}|)  \big)_{i\in I_v},     \big(-\log(|t_{a,v',i}|)\big)_{i\in I_{v'}\!-\!I_v}\bigg)\!\in\!\R_+^{I_{v'}}\qquad \forall~v'\!\in\!\V',
\eEq
and $\la^a\colon \E'\lra \R_+$ defined by
$$
\la_{e'}^{a}=-\log (|\ve_{e',a}|) \qquad \forall~e'\!\in\!\E'
$$
satisfy Condition \ref{Tropical_l} of Definition~\ref{LogMap_dfn}. By (\ref{ratio_e2}) and (\ref{ratio_e1}), $f'$ also satisfies 
Condition \ref{GObs_e} of Definition~\ref{LogMap_dfn}.\\

\noindent
For general $\Gamma$, by the definition of $\eta$ in (\ref{Totaleta_e}), we can choose a set of representatives 
$$
(\ze_{a,v,i})_{a\in \N, v\in \V, i\in I_v}
$$
and coordinates $(z_{\uvec{e},a}=\de_{a,\uvec{e}}z_{\uvec{e}})_{a\in \N, \uvec{e}\in\uvec{\E}}$
such that the leading coefficients $\eta_{\uvec{e},i,a}$ in (\ref{etae_e}) satisfy 
\bEq{assumption}
\eta_{\uvec{e},i,a}= \eta_{\scz\ucev{e},i,a} \qquad \forall~e\!\in\!\E,~i\!\in\!I_e,~a\!\in\!\N.
\eEq
Let $\de_{a,e}=\de_{a,\uvec{e}}\de_{a,\scz\ucev{e}}$, for all $e\!\in\!\E$.
For each $v\!\in\!\V$, choose the representatives 
$$
(\ze_{v',i})_{v'\in \gamma^*(v), i\in I_{v'}}\quad\tn{and}\quad (t_{a,v',i})_{v'\in \gamma^*(v), i\in I_{v'},a\in \N}
$$ 
so that
(\ref{ratio-equal_e2}) and (\ref{ratio-equal_e1}) hold. 
By (\ref{assumption}), we have 
\begin{gather}
\label{ratio_e2new}
\lim_{a\lra \infty}\frac{t_{a,v'_1,i}~\de_{e,a}^{s_{\uvec{e},i}}}{t_{a,v'_2,i}} = 1 \qquad \forall~i\!\in\!I_{v_1'}\!\cap\!I_{v_2'},\\
\label{ratio_e1new}  \lim_{a\lra \infty}t_{a,v'_1,i}~\de_{e,a}^{s_{\uvec{e},i}} = 1\qquad \forall~i\!\in\!I_{v_1'}\!-\!I_{v_2'}.
\end{gather}
With an argument similar to the proof of Proposition~\ref{VertexOrder_Pr}, we can choose these representatives so that further
\begin{gather}
\label{ratio-equal_e2new}
t_{a,v'_1,i}~\de_{e,a}^{s_{\uvec{e},i}} = t_{a,v'_2,i} \qquad \forall~i\!\in\!I_{v_1'}\!\cap\!I_{v_2'},~a\!\in\!\N,\\
\label{ratio-equal_e1new} 
t_{a,v'_1,i}~\de_{e,a}^{s_{\uvec{e},i}} = 1 \qquad \forall~i\!\in\!I_{v_1'}\!-\!I_{v_2'},~a\!\in\!\N.
\end{gather}
Also choose the functions $s^a\colon \V\lra \R^N$ and $\la^a\colon \V\lra \R_+$ satisfying Definition~\ref{LogMap_dfn}.\ref{Tropical_l} so that (\ref{LargeEnough_e}) holds and 
$$
\la^a_e-\log (\de_{e,a})>1\qquad \forall~e\!\in\!\E,~a\!>\!\!>1.
$$
Then, similarly to (\ref{sv-from-log_e}), for $a\!>\!\!>1$, the extended functions $s_{\tn{new}}^a\colon \V'\lra \R^N$ given by 
\bEq{sv-from-log_e2}
s^a_{\tn{new},v'}=\bigg(  \big( s^a_{v,i}-\log(|t_{a,v',i}|)  \big)_{i\in I_v},     \big(-\log(|t_{a,v',i}|)\big)_{i\in I_{v'}\!-\!I_v}\bigg)\!\in\!\R_+^{I_{v'}}\qquad \forall~v'\!\in\!\V',~v\!=\!\gamma^*(v'),
\eEq
and $\la^a_{\tn{new}}\colon \E'\lra \R_+$ given by $\la^a_e-\log (\de_{e,a})$, if  $e\!\in\!\E\!\subset\!\E'$, and 
$$
\la_{e'}^{a}=-\log (|\ve_{e',a}|) \qquad \tn{if}~e'\in\!\E'/\E, 
$$
satisfy Condition \ref{Tropical_l} of Definition~\ref{LogMap_dfn}. By (\ref{ratio_e2}) and (\ref{ratio_e1}) applied to $\E'/\E$, the assumption (\ref{assumption}), and (\ref{ratio-equal_e2new}) and (\ref{ratio-equal_e1new}), $f$ also satisfies 
Condition \ref{GObs_e} of Definition~\ref{LogMap_dfn}.\qed
\end{proofof-LGSmoothV_Th}

\newtheorem*{proofof-Compactness_th}
{Proof of Theorem~\ref{Compactness_th}}
\begin{proofof-Compactness_th}
Similarly to the classical case, consider the sequential convergence topology on $\ov\cM_{g,\mfs}^{\log}(X,D,A)$ given by Definition~\ref{LogConv_dfn}: a subset $W$ of  $\ov\cM_{g,\mfs}^{\log}(X,D,A)$ is closed if every sequence in $W$ has a subsequence with a log-Gromov limit in $W$.  Note that as in \cite[Section~5.1]{MS2}, we must show that convergence with respect to the topology defined above is equivalent to log-Gromov convergence. Since the forgetful map $\ov\cM_{g,\mfs}^{\log}(X,D,A)\!\to\! \ov\cM_{g,k}(X,A)$ is finite-to-one and log-Gromov convergence is a lift of the classical Gromov convergence, this property follows from the the corresponding statement for the Gromov convergence topology on $\ov\cM_{g,k}^{\log}(X,A)$. In other words, the five axioms\footnote{Even though \cite[Section 5.1]{MS2} is about the genus $0$ moduli spaces, the statements used here are valid in all genus.} in \cite[Lemma~5.6.4]{MS2} lift to sequences in $\ov\cM_{g,\mfs}^{\log}(X,D,A)$.\\

\noindent
Suppose $W\!\in\!\ov\cM_{g,k}(X,A)$ is closed and let $W'\!=\!\iota^{-1}(W)$. Let $(f_{a})_{a\in \N}$ be any sequence in $W'$. Its image $(h_a\!=\!\iota(f_{a}))_{a\in \N}$ in $W$ has a subsequence, still denoted by $(h_a)_{a\in \N}$, that Gromov converges to some $h\!\in\!W$. On the other hand, by Theorem~\ref{Log_thm11}, $(f_{a})_{a\in \N}$ has a subsequence that log-Gromov converges to some $f\!\in\!\ov\cM_{g,\mfs}^{\log}(X,D,A)$. By Definition~\ref{LogConv_dfn}, we have $\iota(f)=h$, i.e. $f\!\in\!W'$. Therefore, $W'$ is closed. We conclude that $\iota$ is \textbf{continuous}.\\

\noindent
Let $f$ be an arbitrary log map in $\ov\cM_{g,\mfs}^{\log}(X,D,A)$ with the decorated dual graph $\Gamma$ and $h\!=\!\iota(f)$ be the underlying stable map in $\cM_{g,k}(X,A)$.  Let $(U_a)_{a\in \N}$ be a shrinking basis for the (metrizable) topology of  $\ov\cM_{g,k}(X,A)$ around $h$. Recall from Lemma~\ref{UniqueLog_lm} that every stable map $h$ admits at most finitely many log lifts $f$, each of which is uniquely specified by the vector decorations on the nodes of its dual graph (i.e. the contact data $s_{\uvec{e}}$ at the nodes $q_{\uvec{e}}$). Furthermore, by Lemma~\ref{g0Embedding}, such a lift is unique if the genus is zero. As we explained before Remark~\ref{STDpres_rmk}, for $a$ sufficiently large, by the classical gluing theorem, the domain of every map $h'$ in $U_a$ is obtained from the nodal domain $\Si$ of $h$ by gluing the nodes in a standard way. Furthermore, the image of $h'$ is $C^0$-close to the image of $h$. The dual graph $\Gamma'$ of $h'$ is a contraction of $\Gamma$ in the sense\footnote{Their roles are reversed here.} of (\ref{equ:deg-graphs}). With these identifications, if $f'$ is a log lift of $h'$ in $U_a$, by its \textbf{decoration type}, we mean 
(1) the vector decorations $s_{\uvec{e}}$  at its nodes $q_{\uvec{e}}$, together with (2) the winding\footnote{Contact points with $D_i$ are among the marked/nodal points and are away from the neck region.} number of $h'$ around $D_i$ along the circles $\partial A_{\uvec{e}}$ (see (\ref{partialAe_e})) on every neck $A_e$ obtained from gluing the node $q_e$ of the domain of $h$; see proof of Lemma~\ref{PreLog_lmm}. 
Thus, $f'$ has \textbf{the same decoration type as} $f$ if (1) at every node of the domain of $f'$ the vector decoration $s_{\uvec{e}}$ is the same as the vector decoration at the corresponding node of $f$, (2) on every neck $A_e$ the winding number of $h'$ around $D_i$ along the circle $\partial A_{\uvec{e}}$ is the same as the tangency order $s_{\uvec{e},i}$ for $f$.\\

\noindent 
For $a$ sufficiently large, define $U'_a$ be the set of elements $f'$ in $\ov\cM_{g,\mfs}^{\log}(X,D,A)$ whose image $h'$ under $\iota$ lies in $U_a$ and $f'$ has the same decoration type as $f$.
By Remark~\ref{UniqueLog_rmk2}, the restriction of $\iota$ to $U'_a$ is one-to-one. We show that $U'_a$ is open. Let $(f_{b})_{b\in \N}$ be a sequence in the complement of $U'_a$ that log-Gromov converges to $f'$. After possibly passing to a subsequence, we can assume that the underlying sequence of stable maps $(h_{b})_{b\in \N}$ lies either in $U_a$ or its complement $U^c_a$. In the latter case, by Definition~\ref{LogConv_dfn}, $f'$ belongs to the complement of $U_a'$. In the former case, the decoration type of $f'$ (with respect to $f$) will be the same as the decoration type of $f_{b}$ which is, by definition, different from the decoration type of $f$. Therefore, $f'$ belongs to the complement of $U_a'$. We conclude that  $U_a'$ is open. Furthermore, it is easy to see that  $(U'_a)_{a\in \N}$ is a shrinking basis for the topology of $\ov\cM_{g,\mfs}(X,D,A)$ at $f$. Therefore, the log-Gromov topology of $\ov\cM_{g,\mfs}(X,D,A)$ is first-countable.  \\

\noindent
Hausdorffness is  the consequence of uniqueness of limit in Theorem~\ref{Log_thm11}. If $Y$ is a first-countable topological space and has the property that every convergent sequence has a unique limit then $Y$ is Hausdorff.
Finally, compactness of $\ov\cM_{g,\mfs}^{\log}(X,D,A)$ is the consequence of the existence of limit in Theorem~\ref{Log_thm11}.\qed

\end{proofof-Compactness_th}

%%%%%%%%%%%%%%%%%%%%%%
\section{Log vs. relative compactification; review}\label{LvsR_s}
In Section~\ref{RelComp_ss}, following the description in \cite{FZ}, we review the construction of the relative moduli spaces for smooth symplectic divisors in \cite{LR, IP1}. In Section~\ref{RelToRel_ss}, we show that  the natural forgetful map from the relative compactification to our log compactification is onto.\\

\noindent
First, let us recall some relevant facts from Section~\ref{dbar_ss}.
Suppose $D\!\subset\!(X,\om)$ is a smooth symplectic divisor, $J$ is an $\om$-tame almost complex structure on $X$ such that $J(TD)\!=\!TD$, and $\dbar_{\cN_XD}$ is the $\dbar$-operator in Lemma~\ref{ConnectionToDbar_lmm}. With notaion as in Section~\ref{dbar_ss}, choose a Hermitian connection~$\na^{\cN}$ on $(\cN_XD,\mfi_{\cN_XD})$ such that $\dbar_{\cN_XD}\!=\!\dbar_{\nabla^{\cN}}$. The connection $\na^{\cN}$ gives a splitting of the exact sequence
\bEq{NXVsplit_e}
\begin{split}
0&\lra \pi^*\cN_XD\lra T(\cN_XD)
\stackrel{\nd\pi}{\lra} \pi^*TD\lra0
\end{split}
\eEq
of vector bundles over~$\cN_XD$ which restricts to the canonical splitting over the zero section and is preserved by the multiplication by~$\C^*$; see \cite[Sec~4.1]{FZ}.
Let
\bEq{PXV_e}
\P_X D= \P(\cN_X D\oplus D\!\times\!\C),\qquad D_{0}= \P(0\oplus D\!\times\!\C)~~\tn{and}~~ D_{\infty}=\P(\cN_X D\oplus 0) \subset \P_XD.
\eEq
The splitting of (\ref{NXVsplit_e}) extends to a splitting of the exact sequence 
$$
0\lra T^{\tn{vrt}}(\P_XD) \lra T(\P_XD)
\stackrel{\nd\pi}{\lra} \pi^*TD\lra0,
$$
where $\pi\!:\P_XD\!\lra\!D$ is the bundle projection map induced by (\ref{Firstpi_e}); this splitting restricts to the canonical splittings over $D_{0}\!\cong\!D_{\infty}\!\cong\! D$
and is preserved by the multiplication by~$\C^*$.
Via this splitting, the almost complex structure $J_D$ and 
the complex structure $\mfi_{\cN_XD}$ in the fibers of~$\pi$ induce
an almost complex structure~$J_{X,D}$ on~$\P_XD$ which restricts to $J_D$ on $D_{0}$ and~$D_{\infty}$ and is preserved by the~$\C^*$-action. 
In fact, $J_{X,D}|_{\cN_XD}$ is the almost complex structure $J_{\dbar_{\cN_XD}}$ associated to $\dbar_{\cN_XD}$ described in items~\ref{HoloProj_l}-\ref{HoloSec_l} of Page~\pageref{HoloProj_l}  and is independent of choice of $\na^{\cN}$.
By  property \ref{HoloProj_l}, the projection $\pi\!:\P_XD\!\lra\!D$ is $(J_D,J_{X,D})$-holomorphic. By~\ref{HoloSec_l}, there is a one-to-one correspondence between 
the space of $J_{X,D}$-holomorphic maps $u\colon\!(\Si,\mfj)\!\lra\!(\P_X D,J_{X,D})$ (not mapped into $D_{X,0}$ and $D_{X,\infty}$) and tuples $(u_D,\ze)$ where $u_D\colon (\Si,\mfj)\!\lra\!(D,J_D)$ is a  $J_D$-holomorphic map into $D$ and $\ze$ is a non-trivial meromorphic section of $u_D^*\cN_XD$ with respect to the holomorphic structure defined by $u^*\dbar_{\cN_XD}$.

%---------------------------------------------------
\subsection{Relative compactification}\label{RelComp_ss}
Let $(X,\om)$ be a smooth symplectic manifold, $D\!\subset\!X$ be a smooth symplectic divisor, and $J\!\in\!\cJ(X,D,\om)$.
With notation as in (\ref{PXV_e}), for each $m\!\in\!\N$, let 
$$
X[m]=\big(X\sqcup\{1\}\!\times\!\P_XD\sqcup\ldots\sqcup
\{m\}\!\times\!\P_XD\big)/\!\!\sim\,,
$$
where
$$
D\sim \{1\}\!\times\!D_{\infty}\,,~~~
\{r\}\!\times\!D_{0}\sim \{r\!+\!1\}\!\times\!\P_{\infty}D
\quad  \forall~r\!=\!1,\ldots,m\!-\!1;
$$
see Figure~\ref{relcurve_fig}. This is a \textit{basic} (i.e. there are no triple or higher intersections) SNC variety which is smoothable to (a symplectic manifold deformation equivalent to) $X$ itself.
There exists a continuous projection map $\pi_m\colon\!X[m]\!\lra\!X$ which is identity on $X$ and $\pi$ on each $\P_{X}D$.
We denote by $J_m$ the almost complex structure on~$X[m]$ so that 
$$
J_m|_X=J_X \qquad\hbox{and}\qquad 
J_m|_{\{r\}\times\P_XD}=J_{X,D} \quad \forall~r=1,\ldots,m.$$
For each $(c_1,\ldots,c_m)\!\in\!\C^*$, define
$$
\Theta_{c_1,\ldots,c_m}\!:X[m]\lra X[m] \qquad\hbox{by}\quad
\Theta_{c_1,\ldots,c_m}(x)=\begin{cases}(r,[c_rv,w]),&\hbox{if}~x\!=\!(r,[v,w])\!\in\!\{r\}\!\times\!\P_XD;\\x,&\hbox{if}~x\!\in\!X.
\end{cases}
$$
This diffeomorphism is biholomorphic with respect to~$J_m$ and
preserves the fibers of the projection $\P_XD\!\lra\!D$
and the sections~$D_{0}$ and~$D_{\infty}$.\\

\begin{figure}
\begin{pspicture}(38,-1)(11,3.2)
\psset{unit=.3cm}
\psline[linewidth=.1](48,-2)(80,-2)
\psline[linewidth=.1](48,2)(80,2)
\psline[linewidth=.1](48,6)(80,6)

\rput(51,6.8){\small{$D$}}\rput(51,5.2){\small{$D_{\infty}$}}
\rput(51,2.8){\small{$D_{0}$}}\rput(51,1.2){\small{$D_{\infty}$}}
\rput(51,-1.2){\small{$D_{0}$}}

\pscircle*(62,-2){.2}\rput(62,-2.8){\small{$z^2$}}
\pscircle*(67.8,-2){.2}\rput(67.8,-2.8){\small{$z^3$}}

\rput(46,8){$X$}\rput(46,4){$1\!\times\!\P_XD$}\rput(46,0){$2\!\times\!\P_XD$}

\psccurve(63,10)(62,10.25)(60,12)(57,8)(60,6)(63,8)(66,6)(69,8)(66,12)(64,10.25)\pscircle*(60,6){.2}\pscircle*(66,6){.2}

\pscircle*(60,2){.2}\pscircle*(64.2,2){.2}\pscircle*(67.8,2){.2}\pscircle*(60,4){.2}
\pscircle(60,4){2}
\rput(61,4){\small{$z^1$}}
\pscircle(67.8,0){2}
\psccurve(64,2.1)(66,3)(68,2.1)(66,6)
\psccurve(59.8,1.9)(62,1)(64.4,1.8)(62,-2)

\psarc(60,6){3}{67}{113}
\psarc(60,11.5){3}{240}{300}

\end{pspicture}
\caption{A relative map with $k\!=\!3$ and $\mfs\!=\!(0,2,2)$
into the expanded degeneration $X[2]$.}
\label{relcurve_fig}
\end{figure}

\noindent
The moduli space of relative stable curves for $(X,D)$ in~\cite[Sec 7]{IP1} is defined  in the following way.
With slight modification, we follow the description in \cite{FZ}.
Suppose $k\!\in\!\N$, $A\!\in\!H_2(X,\Z)$, and $\mfs\!=\!(s_1,\ldots,s_{k})\!\in\!\N^{k}$ 
is a tuple satisfying
\bEq{bsumcond_e}
\sum_{a=1}^k s_a=A\cdot D.
\eEq

\noindent
A  level zero genus $g$ $k$-marked degree $A$ \textbf{relative} $J$-holomorphic map into~$X$ of contact type~$\mfs$ with $D$ is
simply a stable $J$-holomorphic map in $\ov\cM_{g,k}(X,A)$
such~that 
\bEq{RelativeOrder_e}
u^{-1}(D)\subset \big\{z^{1},\ldots,z^{k}\big\}, \quad
\tn{ord}_{z^{a}}\big(u,D\big)=s_a \quad \forall a\!=\!1,\ldots,k.
\eEq
For $m\!\in\!\Z_+$, a level $m$ $k$-marked \textbf{relative} $J$-holomorphic map of contact type~$\mfs$ is
a continuous map $u\!:\Si\!\lra\!X[m]$ 
from a marked connected nodal curve $\big(\Si,\mfj,\vec{z}\!=\!(z^1,\ldots,z^k)\big)$ 
such~that 
$$
u^{-1}\big(\{m\}\!\times\!D_0\big)
\subset\big\{z^{1},\ldots,z^{k}\big\}, \quad
\tn{ord}_{z^a}(u,\{m\}\!\times\!D_0)\!=\!s_a\quad \forall~z^a\!\in\!u^{-1}\big(\{m\}\!\times\!D_0\big),
$$
 $s_a\!=\!0$ if and only if $z^a\!\notin\!u^{-1}\big(\{m\}\!\times\!D_0\big)$, and the restriction of~$u$ to each irreducible component $\Si_j$ of~$\Si$ is either 
\bEn
\item a $J$-holomorphic map to $X$ such that the set $u|_{\Si_j}^{\,-1}(D)$ 
consists~of the nodes joining~$\Si_j$ to irreducible components of~$\Si$ mapped to
$\{1\}\!\times\!\P_XD$, or
\item a $J_{X,D}$-holomorphic map to $\{r\}\!\times\!\P_XD$ for some $r\!=\!1,\ldots,m$ such~that 
\bEnalph
\item the set $u|_{\Si_j}^{\,-1}(\{r\}\!\times\!D_{\infty})$ 
consists~of the nodes~$q_{j,i}$ joining~$\Si_j$ to irreducible components of~$\Si$ mapped to
$\{r\!-\!1\}\!\times\!\P_XD$ if $r\!>\!1$ and to~$X$ if $r\!=\!1$
and 
$$
\tn{ord}_{q_{j,i}}\big(u,D_{\infty}\big)=
\begin{cases}
\tn{ord}_{q_{i,j}}(u,D_{0}),&\hbox{if}~r\!>\!1,\\
\tn{ord}_{q_{i,j}}(u,D),&\hbox{if}~r\!=\!1,
\end{cases}
$$ 
where $q_{i,j}\!\in\!\Si_{i,j}$ is the point identified with~$q_{j,i}$, 
\item  if $r\!<\!m$, the set $u|_{\Si_j}^{\,-1}(\{r\}\!\times\!D_{0})$ 
consists~of the nodes joining~$\Si_j$ to irreducible components of~$\Si$ mapped to
$\{r\!+\!1\}\!\times\!\P_XD$;
\eEnalph
\eEn
see Figure~\ref{relcurve_fig}.
The \textbf{genus} and the \textbf{degree} of such a map $u\!:\!\Si\!\lra\!X[m]$
are the arithmetic genus of~$\Si$ and the homology class
$$
A=\big[\pi_m\!\circ\!u\big]\in H_2(X,\Z).
$$

\noindent
Two tuples $(u_\al,\Si_\al,\mfj_\al,\vec{z}_\al)$ and $(u_\beta,\Si_\beta,\mfj_\beta,\vec{z}_\beta)$  as above are \textbf{equivalent}
if there exist a biholomorphic map $\varphi\colon(\Si_\al,\mfj_\al)\!\lra\!(\Si_\beta,\mfj_\beta)$ 
and \hbox{$c_1,\ldots,c_m\!\in\!\C^*$} so~that 
$$
\varphi(z^a_\al)=z^a_\beta \quad \forall~a\!=\!1,\ldots,k \quad\hbox{and}\quad
u_\beta=\Theta_{c_1,\ldots,c_m}\circ u_\al\circ\varphi.
$$
A tuple as above is \textbf{stable} if it has finitely many automorphisms (self-equivalences).\\

\noindent
If $A\!\in\!H_2(X,\Z)$, $g,k\!\in\!\N$, and $\mfs\!=\!(s_1,\ldots,s_{k})\!\in\!\N^{k}$ 
is a tuple satisfying~(\ref{bsumcond_e}), then the \textbf{relative moduli space} 
\bEq{relmoddfn_e}
\ov\cM_{g,\mfs}^{\tn{rel}}(X,D,A)
\eEq
is the set of equivalence classes of such connected stable $k$-marked genus~$g$ degree~$A$ $J$-holomorphic maps
into  $X[m]$ for any $m\!\in\!\N$. If $X$ is compact, the latter space has a natural compact Hausdorff topology. 

\begin{remark}
In (\ref{bsumcond_e}), we are allowing $s_a$ to be zero for some $a\!\in\!\{1,\ldots,k\}$. A marked point $z$ with contact order $0$ has image away from $D$ (or $D_0$, $D_\infty$). Therefore, such points are ordinary marked points as in the classical moduli spaces of $J$-holomorphic curves. In the literature, marked points are usually divided into the classical part $(z^1,\ldots,z^k)$ and the relative part $(z^{k+1},\ldots,z^{k+\ell})$ such that $s_a\!=\!\tn{ord}_{z^{k+a}}(u,D)\!>\!0$ and $\sum_{a=1}^\ell s_a\!=\!A\cdot D$. 
Then the moduli space (\ref{relmoddfn_e}) is denoted by $\ov\cM^{\tn{rel}}_{g,k,\mfs}(X,D,A)$ with $\mfs\!\in\!(\Z_+)^\ell$. This sort of separation works fine in the relative case,  because there are only two types of points: in $D$ or away from $D$. In the general SNC case $D\!=\!\bigcup_{i\in [N]} D_i$, however, there are $2^N$ types of points and it is notationally cumbersome (and useless) to divide points into separate groups based on their type.
\end{remark}

\begin{remark}\label{dbarvsD_rmk}
Let $(X,\om)$ be a smooth symplectic manifold, $D\!\subset\!X$ be a smooth symplectic divisor, and $J$ be an $\om$-tame almost complex structure on $X$ such that $J(TD)\!=\!TD$. If $u\colon\!(\Si,\mfj)\!\lra\!(X,J)$ is $J$-holomorphic, the linearization of the Cauchy-Riemann operator (\ref{CR_e}) at~$u$ 
is given~by
\begin{gather}
D_u\dbar\colon\! \Gamma\big(\Si;u^*TX\big)\lra \Gamma\big(\Si,\Om_{\Si,\mfj}^{0,1}\!\otimes_{\C}\!u^*TX\big),\notag\\
\label{Dudfn_e}
D_u\dbar(\xi)= u^*\dbar_{\wh\na}
+\frac{1}{4} N_{J}(\xi,\nd u),
\end{gather}
where $\wh\na$ is the $\C$-linear connection in (\ref{CLConn_e2}) and $\dbar_{\wh\na}$ is the associated $\dbar$-operator on $\Gamma(X,TX)$ in Lemma~\ref{ConnectionToDbar_lmm}; see  \cite[Ch~3.1]{MS2}. The kernel of $D_u\dbar$ corresponds to infinitesimal deformations of $u$ (over the fixed domain $(\Si,\mfj)$) and the cokernel of that is the obstruction space for integrating infinitesimal deformations  to actual deformations.\\

\noindent
If furthermore $\tn{Im}(u)\!\subset\!D$, then the linearization map $D_u\dbar$, defined in (\ref{Dudfn_e}), satisfies
$$
D_u\bar\partial \big(\Ga(\Si,u^*TD)\big)\subset \Gamma\big(\Si,\Om_{\Si,\mfj}^{0,1}\!\otimes_{\C}\!u^*TD\big),
$$
because the restriction of  $D_u\bar\partial$ to $\Ga(\Si,u^*TD)$ is the linearization\footnote{the linearization of (\ref{CR_e}) is independent of the choice of the connection at every $J$-holomorphic map.}
of the $\dbar$-operator at~$u$  for the space of maps into~$D$.
Thus, $D_u\bar\partial$ descends to a first-order differential operator
\bEq{DuNXV_e}
D_u^{\cN_XD}\bar\partial\!: \Ga(\Si,u^*\cN_XD)\lra \Gamma\big(\Si,\Om_{\Si,\mfj}^{0,1}\!\otimes_{\C}\!u^*\cN_XD\big).
\eEq
If $J\!\in\!\cJ(X,D,\om)$, i.e. (\ref{intInnormal_e0}) holds,  then the normal part of $N_{J}(*,\nd u)$ vanishes.  From (\ref{Dudfn_e}) and Lemma~\ref{ConnectionToDbar_lmm} we conclude that
$$
D_u^{\cN_XD}\bar\partial\!=\! u^* \dbar_{\cN_XD}
$$
is a complex linear operator. 
From another point of view, we can use (\ref{intInnormal_e0}) to show that certain sequence of almost complex structures on the normal bundle $\cN_XD$ converges to $J_{X,D}$, see Lemma~\ref{JLemma_e}. 
\end{remark}

%-------------------------------------------------
\subsection{Relative vs. Log compactification}\label{RelToRel_ss}
In this section, for the case where $D$ is smooth (i.e. $N\!=\!1$ in Definition~\ref{LogMap_dfn}), we compare $\ov\cM^{\tn{rel}}_{g,\mfs}(X,D,A)$  and $\ov\cM^{\tn{log}}_{g,\mfs}(X,D,A)$. Proposition~\ref{ReltoLog_prp} shows that the latter is smaller and there is a projection map from the relative compactification onto the log compactification. This is expected, since the notion of nodal log curve involves more $\C^*$-quotients on the set of meromorphic sections. 
In the algebraic case, \cite[Thm~1.1]{AMW} shows that an algebraic analogue of the projection map (\ref{RelToLogPi_e}) induces an equivalence of virtual fundamental classes.
We expect the same to hold for the invariants/VFC arising from our log compactification.\\

\noindent
First, we start with a simple lemma that highlights the relation between Definition~\ref{LogMap_dfn}.\ref{Tropical_l} and the layer structure in the relative compactification. In the following, when $D$ is smooth ($N\!=\!1$), for a (pre-)log map with the decorated dual graph $\Gamma(\V,\E,\L)$ we define
\bEq{Divide_e}
\aligned
&\V_i=\{v\in \V\colon |I_v|\!=\! i \},\quad  \E_i=\{e\in \E\colon |I_e|\!=\!i \},\qquad \tn{with } i\!=\!0,1,\\
&\E_{1,0}=\{e\in \E\colon |I_e|\!=\!1,~s_{\uvec{e}}\!=\!0 \},\quad \E_{1,\star}=\{e\in \E\colon |I_e|\!=\!1,~s_{\uvec{e}}\!\neq \!0 \}. 
\endaligned
\eEq  

\bLm{Partial-Order_lmm}
Let $D\!\subset\! (X,\om)$ be a smooth symplectic divisor, $J\!\in\!\cJ(X,D,\om)$, and 
\bEq{SmoothVPreLogMap_e}
\big[f\equiv \big((u_v,[\ze_v],C_v)_{v\in \V_1}, (u_v,C_v)_{v\in \V_0} \big)\big]\!\in\!\cM^{\tn{plog}}_{g,\mfs}(X,D,A)_\Gamma
\eEq
be a pre-log $J$-holomorphic curve with dual graph $\Gamma(\V,\E,\L)$.
The  there exists a function $s\colon\!\V\!\lra\!\R_{\geq 0}$ satisfying Definition~\ref{LogMap_dfn}.\ref{Tropical_l}, if and only if 
the relations  
\bEnalph
\item $v_1\!\approx_\Gamma \!v_2$ if $v_1$ and $v_2$ are connected and $s_{\uvec{e}}\!=\!0$ for any $\uvec{e}\!\in\!\uvec{\E}_{v_1,v_2}$, and 
\item $v_1\!\prec_\Gamma\!v_2$ if $v_1$ and $v_2$ are connected and $s_{\uvec{e}}\!>\!0$ for any $\uvec{e}\!\in\!\uvec{\E}_{v_1,v_2}$, 
\eEnalph
are independent of the choice of $\uvec{e}\!\in\!\uvec{\E}_{v_1,v_2}$ (i.e. they are well-defined), and generate a partial order $\preceq_{\Gamma}$ on $\V$.
\eLm
\noindent
Note that for a classical edges $e$ connecting  $v_1,v_2\!\in\!\V_0$, since $I_e\!=\!\eset$ by (\ref{Ieunion_e}), we always have 
$$
s_{\uvec{e}}\!=\!0\!\in\!\{0\}=\R^{\eset}\!\subset\! \R^{N=1}\!=\!\R.
$$

\bPf 
If $(a)$ and $(b)$ define a partial order $(\V,\preceq_{\Gamma})$, we construct $s\colon\!\V\!\lra\! \R$ satisfying Definition~\ref{LogMap_dfn}.\ref{Tropical_l} in the following way. For every $v\!\in\!\V_0$ define $s_v\!=\!0$. 
Let $\V_{\tn{min}}^{(1)}$ be the subset of minimal vertices in $ \V_1$. For every $v\!\in\!\V_{\tn{min}}^{(1)}$ define $s_v\!=\!1$. Having constructed $\V_{\tn{min}}^{(1)},\ldots,\V_{\tn{min}}^{(\ell)}$, let $\V_{\tn{min}}^{(\ell+1)}$ be the subset of minimal vertices in 
$$
\V_1\!-\!\big(\V_{\tn{min}}^{(1)}\cup\!\cdots \!\V_{\tn{min}}^{(\ell)}\big).
$$
For every $v\!\in\!\V_{\tn{min}}^{(\ell+1)}$ define $s_v\!=\!\ell+1$. This function clearly satisfies Definition~\ref{LogMap_dfn}.\ref{Tropical_l}. Conversely, given such a function $s\colon\!\V\!\lra\!\R$ satisfying Definition~\ref{LogMap_dfn}.\ref{Tropical_l}, define $v_1\!\approx_\Gamma\!v_2$ (resp. $v_1\!\prec_{\Gamma}\!v_2$)  if they are connected by a path and $s_{v_1}\!=\!s_{v_2}$ (resp. $s_{v_1}\!<\!s_{v_2}$). This is a partial order whose defining conditions match with (a) and (b).
\ePf

\bLm{evGaSmooth_lmm}
With notation as in Lemma~\ref{Partial-Order_lmm}, the pre-log curve $f$ satisfies Definition~\ref{LogMap_dfn}.\ref{GObs_e} if and only if there exists a set of representatives $\{\ze_v\}_{v\in \V_1}$ such that 
\bEq{MatchAtEsetNodes_e}
\ze_v(q_{\uvec{e}})=\ze_{v'}(q_{\scz\ucev{e}})\qquad \forall~v,v'\!\in\!\V_1,~e\!\in\! \E_{v,v'}\quad \tn{s.t.}\quad s_{\uvec{e}}\!=\!0.
\eEq
\eLm

\bPf
The last equation is well-defined by Definition~\ref{PreLogMap_dfn}.\ref{MatchingOrders_l}.
Then the homomorphism (\ref{DtoT_e}) (corresponding to some fixed orientation $O$ on $\E$) takes the form
\bEq{DtoT_e2}
 \Z^{\E_0} \oplus \Z^{\E_1}\oplus \Z^{\V_1}  \stackrel{\vr}{\xrightarrow{\hspace*{1.5cm}}} \Z^{\E_1}
\eEq
where $\vr|_{ \Z^{\E_0}}\!\equiv\!0$, $\vr(1_{e})\!=\!s_{\uvec{e}}\!\in\!\Z$ for all $e\!\in\!\E_1$, and 
$$
\vr(1_{v}\equiv 1_{v,1})_{e}=\begin{cases} 1_{e} &\mbox{if } v_1(\uvec{e})=v; \\ 
-1_e & \mbox{if } v_2(\uvec{e})=v;\\
0 & \mbox{if $e$ is a loop or otherwise.} \end{cases} 
$$
Therefore, tensoring (\ref{DtoT_e2}) with $\C$, the cokernel $\C\K_\C$ of $\vr_\C$ is equal to the cokernel of the induced map
$$
\C^{\V_1}  \stackrel{\ov\vr_\C}{\xrightarrow{\hspace*{1.5cm}}} 
\C^{\E_{1,0}}.
$$
Fix an arbitrary set of representatives  
\bEq{zeRep_e2}
\big(\ze_v\in \Om_{\tn{mero}}(\Si_v, u_v^*\cN_{X}D)\big)_{v\in \V_1}.
\eEq
By (\ref{LocalCoord_e}) and (\ref{etae_e}), for every $e\!\in\!\E_{1,0}$ with $v\!=\!v_1(\uvec{e})$ and $v'\!=\!v_2(\uvec{e})$, we have
$$
\eta_e= \ze_v(q_{\uvec{e}})/\ze_{v'}(q_{\scz\ucev{e}})\in \C^*.
$$
Therefore, 
$$
\ov\eta\!\equiv\!(\eta_e)_{e\in \E_{1,0}}\!\in\! \prod_{e\in \E_{1,0}} (\C^*)^{\E_{1,0}}
$$
is equal to $(1)_{e\in \E_{1,0}}$ if and only if (\ref{MatchAtEsetNodes_e}) holds. Since cokernel of $\ov\vr_\C$ coincides with cokernel of $\vr_\C$, the element 
$$
[\eta]\!\in\!(\C^*)^{\E}/\tn{exp}(\tn{im}(\vr_\C))
$$ 
in (\ref{[eta]_e}) is the identity element if and only if 
$$
[\ov\eta]\!\in\!(\C^*)^{\E_{1,0}}/\tn{exp}(\ov\vr_\C(\C^{\V_1}))
$$ 
is the identity element. The latter holds if and only if is there exists a rescaling of the sections $(\ze_v)_{v\in \V_1}$ for which (\ref{MatchAtEsetNodes_e}) holds.
\ePf

\bPr{ReltoLog_prp}
Let $D\!\subset\! (X,\om)$ be a smooth symplectic divisor, $J\!\in\!\cJ(X,D,\om)$, and $\mfs\!\in\!\N^k$. Then there exists a natural surjective map
\bEq{RelToLogPi_e}
\pi\colon\ov\cM^{\tn{rel}}_{g,\mfs}(X,D,A)\lra \ov\cM_{g,\mfs}^{\log}(X,D,A).
\eEq
\ePr

\bPf
For each relative curve $f$, $\pi(f)$ is the log curve obtained by forgetting those unstable $\P^1$-components of the domain which are isomorphically mapped to the trivial fibers of $\P_XD$, and by restricting the equivalence class of each section defining a map into a $\P_XD$ to the equivalence classes of its restrictions to each connected component. The required function $s\colon\!\V(\Gamma)\!\lra\!\R_{\geq 0}$ in Definition~\ref{LogMap_dfn}\ref{Tropical_l} can be taken to be the one given by the layer structure of the relative moduli space. Moreover, by Lemma~\ref{evGaSmooth_lmm}, $\pi(f)$ satisfies (\ref{GammaLogStratum_e}) because a set of sections representing $f$ have equal values at the nodes $q_e$ with $I_e\!=\!\{1\}$ and $s_{\uvec{e}}\!=\!0$.\\

\noindent
Conversely, let $f$ be any log map with dual graph $\Gamma$. By Corollary~\ref{m_0integral}, we can assume that the function $s\colon \V(\Gamma)\!\lra\!\R_{\geq 0}$ in Definition~\ref{LogMap_dfn}.\ref{Tropical_l} is integral. Furthermore, we take $s$ so that $\max(s)$ is the smallest among all such $s$.
For each connected component $\Si_v$ of $\Si$ in $f$ with $I_v\!=\!\{1\}$, choose an arbitrary section $\ze_v$ representing the equivalence class $[\ze_v]$ in $f$. By Lemma~\ref{evGaSmooth_lmm}, we can choose these sections to have equal values at the nodes $q_e$ with $I_e\!=\!\{1\}$ and $s_{\uvec{e}}\!=\!0$. Define a relative map $\wt{f}$ whose restriction to $\Si_v$ is the map  corresponding to $\ze_v$ into the $s_v$-th $\P_XD$ and such that disconnected nodes are connected by adding extra $\P^1$-components to the domain and by mapping them bijectively to the $\P^1$-fibers of $\P_{X}D$.
Since $\max(s)$ is the smallest among all such $s$, there is at least one non-trivial component in each $\P_XD$ of the expanded degeneration $X[\max(s)]$; i.e. $\wt{f}$ defines a stable map into $X[\max(s)]$. It is clear from the construction that $\pi(\wt{f})\!=\!f$.
\ePf

\noindent
Next, we give an example where the projection map  (\ref{RelToLogPi_e}) is non-trivial and both the relative and the log moduli spaces are smooth. The relative moduli space in this example is some blowup of the log moduli space.
\bEx{P1Deg1TwoPt_ex}
Let $X\!=\!\P^1$, $D\!=D_1=\!\pt_1\!\sqcup \pt_2$ (so $N\!=\!1$) be the disjoint union of two points, $g\!=\!0$, $k\!=\!4$, and $A\!=\![1]\!\in\!H_2(\P^1,\Z)\!\cong\! \Z$. Therefore $\mfs=(0,0,1,1)\in \N^4$ (or a permutation of that) is the only option for the contact pattern. Then the relative moduli space
$\ov\cM^{\tn{rel}}_{0,\mfs}(X,D,[1])$
can be identified with a blowup of $\P^1\!\times\!\P^1$ at $4$ points, while $\ov\cM_{0,\mfs}^{\log}(X,D,[1])$ can be identified with a  blowup of $\P^1\!\times\!\P^1$ at $2$ (of those) points. The projection map in (\ref{RelToLogPi_e}) corresponds to the blowdown of the two extra exceptional  curves.
\eEx

\section{Comments on deformation theory and gluing}\label{DG_s}

\subsection{Deformation theory and the expected dimension}\label{EXDim_ss}
In this section, we outline a Fredholm set up for studying the deformation  theory of log $J$-holomorphic maps and make some conclusions. This setup is discussed in details in \cite{FT2}, where it is also extended to log $(J,\nu)$-holomorphic maps. \\

\noindent
In the case of the classical moduli space of stable $J$-holomorphic curves $\ov\cM_{g,k}(X,A)$, for a $J$-holomorphic map $u\colon\!(\Si,\mfj)\!\lra\!(X,J)$ with smooth domain, the linearization $D_u\dbar$ of the Cauchy-Riemann equation in (\ref{Dudfn_e}) is Fredholm. 
Therefore, the real vector spaces
$$
\Def(u)\!=\!\ker(\tn{D}_u\dbar)\quad \tn{and} \quad \Obs(u)\!=\!\coker(\tn{D}_u\dbar)
$$
are finite dimensional. The first space corresponds to infinitesimal deformations of $u$ (over the fixed domain $C$) and the second one is the obstruction space for integrating the elements of $\Def(u)$ to actual deformations. In the nodal case, the kernel $\Def(u)$ of the similarly defined linearization map in \cite[Sec~6.3]{FF} corresponds to infinitesimal deformations of $u$ in the stratum $\ov\cM_{g,k}(X,A)_\Gamma$. Deformations into $\ov\cM_{g,k}(X,A)$ correspond to gluing the nodes of the domain with gluing parameters from $\C^{\E}$ and the gluing is virtually un-obstructed, i.e.  if $\Obs(u)\!=\!0$, for every sufficiently small smoothing $(\Si',\mfj')$ of the nodes of the domain $(\Si,\mfj)$, there exists a $J$-holomorphic map $u'\colon\! (\Si',\mfj')\!\lra\!(X,J)$ close $u$; see \cite[Thm~6.3.5]{FF} for $\Obs(u)\!\neq \!0$. In other words, moduli spaces $\ov\cM_{g,k}(X,A)$ are virtually smooth (orbifolds) and the ``virtual normal cone" of the stratum $\ov\cM_{g,k}(X,A)_\Gamma$ is an (orbi-) bundle of rank $|\E|$. 
For the log moduli spaces defined in this paper, as (\ref{GammaLogStratum_e}) indicates, there are new obstructions for smoothability of nodal pre-log curves. The claim is that, in addition to a logarithmic version of $\tn{D}_u\dbar$, the deformation-obstruction is encoded in the combinatorial linear map (\ref{DtoT_e}).\\

\noindent
In the setting of Theorem~\ref{Compactness_th}, suppose $(u,\Si,z^1,\ldots,z^k)$ is an element of 
$$
\cM_{g,\mfs}(X,D,A)\!\subset\! \cM_{g,k}(X,A),\qquad\tn{for some}\quad\mfs\!=\!\big(s_a\!=\!(s_{ai})_{i\in [N]}\big)_{a\in [k]}\!\in\!(\N^N)^k;
$$ 
i.e. $\Si$ is smooth, $u^{-1}(D) \!\subset\! \{z^1,\ldots,z^k\}$ and $\ord_{z^a}(u,D_i)\!=\!s_{ai}$, for all $i\!\in\![N]$ and $a\!\in\![k]$. If $(\om,\cR,J)\!\in\!\AK(X,D)$, let  $TX(-\log D)$ be log tangent bundle in  \cite[(8)]{FMZ2}, and if $J$ is integrable, let  $TX(-\log D)$ be the usual holomorphic logarithmic tangent bundle.
There is a natural complex linear homomorphism 
$$
\iota\colon TX(-\log D)\lra TX
$$ 
(covering $\tn{id}_X$) that is an isomorphism away from $D$. This homomorphism induces similarly denoted maps 
$$
\aligned
\iota_1&\colon \Gamma(\Si,u^*TX(-\log D))\lra  \Gamma(\Si,u^*TX),\\
 \iota_2&\colon \Gamma(\Si,\Om^{0,1}_{\Si,\mfj}\otimes_\C u^*TX(-\log D))\lra  \Gamma(\Si,\Om^{0,1}_{\Si,\mfj}\otimes_\C u^*TX).
\endaligned
$$   
The following is one of the key steps in understanding the deformation theory of $J$-holomorphic maps relative to an SNC divisor; see \cite[Sec~5.1]{FT2}. 

\begin{theorem}[\cite{FT2}]\label{CRLLog_thm}
With notation as above, the linearization $\tn{D}_u\dbar$ naturally lifts to a Fredholm linear map 
\bEq{CRLog_e}
\tn{D}^{\log}_u\dbar\colon W^{\ell,p}(\Si,u^*TX(-\log D))\lra W^{\ell-1,p}(\Si,\Om^{0,1}_{\Si,\mfj}\otimes_\C u^*TX(-\log D))
\eEq
such that $\iota_2\circ \tn{D}^{\log}_u\dbar = \tn{D}_u\dbar \circ \iota_1$ over the space of smooth sections. Furthermore, if $\tn{coker}(\tn{D}^{\log}_u\dbar)\!=\! 0$, the set of $J$-holomorphic maps (the marked domain is fixed) of contact type $\mfs$ close to $u$ (in a suitable Banach manifold) form an oriented smooth manifold of real dimension
\bEq{RRLog_e}
2\Big( \tn{deg}(u^*TX(-\log D) )\!+\! \dim_\C\! X (1\!-\!g)\Big).
\eEq
\eTh
\noindent
Note that (\ref{RRLog_e}) follows from Riemann-Roch and (\ref{CRLog_e}). Considering the deformations of the marked domain $(\Si,\mfj,z^1,\ldots,z^k)$, it follows from (\ref{RRLog_e}) that the expected dimension of $\cM_{g,\mfs}(X,D,A)$ is equal to the naive dimension count (\ref{dlog_e}).\\

\noindent
Next, consider a log map $f\!=\!(u, [\ze_i]_{i\in I}, \Si,z^1,\ldots,z^k)$ in the stratum $\cM_{g,\mfs}(X,D,A)_I$, i.e. $\Si$ is smooth, $u(\Si)\!\subset\! D_I$ for a non-trivial maximal subset $I\!\subset\![N]$, $\ord_{z^a}(u,D_i)\!=\!s_{ai}$ for all $i\!\notin\! I$, and $\ord_{z^a}(\ze_i)\!=\!s_{ai}$ for all $i\!\in\! I$.  Forgetting the meromorphic sections, by Remark~\ref{UniqueLog_rmk}, we get an inclusion map
$$
\cM_{g,\mfs}(X,D,A)_I\hookrightarrow \cM_{g,\ov{\mfs}}(D_I,\ov{D},A), \qquad (u, [\ze_i]_{i\in I}, \Si,z^1,\ldots,z^k)\lra (u, \Si,z^1,\ldots,z^k)
$$
where 
$$
\ov{D}=\bigcup_{i\in S-I} D_{I\cup i}\subset D_I\qquad \tn{and}\qquad \ov\mfs\!=\!\big(s_a\!=\!(s_{ai})_{i\in [N]-I}\big)_{a=1}^k\!\in\!(\N^{[N]-I})^k.
$$
With $(D_I,\ov{D})$ in place of $(X,D)$ in (\ref{CRLog_e}), deformation theory of $\cM_{g,\ov\mfs}(D_I,\ov{D},A)$ is given by the restriction of $\tn{D}^{\log}_u\dbar $ to $TD_I(-\log \ov{D})$. It is worth mentioning that restricted to $D_I$, there is a natural isomorphism
\bEq{Tlogdec_e}
 TX(-\log D)|_{D_I} \cong TD_I(-\log \ov{D})\oplus D_I\times\C^I.
\eEq

\bLm{HodgeObs_lmm}
There exists a map $P_I=(P_{I,i})_{i\in I}\colon\cM_{g,\ov{\mfs}}(D_I,\ov{D},A)\!\lra\! (\tn{Pic}^0(\Si))^I$ such that 
$$
\cM_{g,\mfs}(X,D,A)_I=P_I^{-1}\big(\cO^I\big).
$$ 
In particular, 
$$\cM_{0,\mfs}(X,D,A)_I=\cM_{0,\ov{\mfs}}(D_I,\ov{D},A).
$$ 
Here $\tn{Pic}^0(\Si)$ is the group of degree $0$ holomorphic lines bundles on $(\Si,\mfj)$ and $\cO\!\in\!\tn{Pic}^0(\Si)$ is the trivial holomorphic line bundle. 
\eLm 

\bPf
For each $i\!\in\!I$, define 
$$
P_{I,i}\big([u, \Si,z^1,\ldots,z^k]\big)=u^*\cN_{X}D_i \otimes \cO\big(-\sum_{a=1}^k s_{ai} z^a\big)\in\!\tn{Pic}^0(\Si),
$$
where $\cO(-\sum_{a=1}^k s_{ai} z^a)$ is the line bundle corresponding to the divisor $-\sum_{a=1}^k s_{ai} z^a$. Therefore, 
$$
P_{I,i}\big([u, \Si,z^1,\ldots,z^k]\big)\!=\!\cO
$$ 
if and only if there exists a meromorphic section $\ze_{I,i}$ of $u^*\cN_{X}D_i$ with zeros/poles of order $s_a$ and $z^a$ (and nowhere else).
\ePf

\noindent
The conclusion is that the deformation/obstruction theory of the stratum $\cM_{g,\mfs}(X,D,A)_I$ is given by $\tn{D}^{\log}_u\dbar$  on $\cM_{g,\ov{\mfs}}(D_I,\ov{D},A)$ and the linearization of  $P_I$. By (\ref{dlog_e}) and Lemma~\ref{HodgeObs_lmm}, the expected real dimension of $\cM_{g,\mfs}(X,D,A)_I$ is  
\bEq{VirDimRel_e}
2\Big(c_1^{TX(-\log D)}(A)+(\dim_\C X-3)(1-g)+k-|I|\Big).
\eEq
Via the identification (\ref{Tlogdec_e}), the maps $\tn{D}^{\log}_u\dbar$ on $\cM_{g,\ov{\mfs}}(D_I,\ov{D},A)$ and $P_I$ can be combined into a single Fredholm operator as in (\ref{CRLog_e}); see \cite[Sec~5.2]{FT2}. \\

\noindent
Moving to the nodal case, with notation as in (\ref{FVersion_e}), let 
\bEq{SiCone_e}
\si=\si(\Gamma)=\K_\R\cap \big(\R_{\geq 0}^\E \oplus \bigoplus_{v\in \V} \R_{\geq 0}^{I_v}\big)\!\subset\!\K_\R
\eEq
be the cone of non-negative elements in the kernel of $\vr_\R\colon\! \D_\R\!\lra\!\T_\R$. This cone is independent of the choice of the orientation $O$ used to define (\ref{DtoT_e}); in fact, by (\ref{kernel_e}), 
\bEq{SiCone_e2}
\si=\big\{\big((\la_e)_{e\in \E},(s_v)_{v\in \V}\big)\!\in\!\ \R_{\geq 0}^\E\oplus \bigoplus_{v\in \V} \R_{\geq 0}^{I_v} \colon~~ s_{v}\!-\!s_{v'}\!=\!\la_es_{\uvec{e}}\quad\forall~v,v'\!\in\!\V,~\uvec{e}\!\in\!\uvec{\E}_{v',v}\big\}.
\eEq
The integral lattice underlying $\si$ coincides with the monoid $Q^{\vee}$ in \cite[Sec~2.3.9]{ACGS}.

\bLm{Convexity_lmm}
For every $\Gamma\!\in\!\tn{DG}(g,\mfs,A)$, $\si(\Gamma)$ is a top-dimensional strictly convex rational polyhedral cone in $\K_\R(\Gamma)$. 
\eLm

\bPf
The functions $s$ and $\la$ in Definition~\ref{LogMap_dfn}.\ref{Tropical_l}  define an element $m_+$ of 
\bEq{KcapPos_e}
\K_\R\cap \big(\R_{+}^\E \oplus \bigoplus_{v\in \V} \R_{+}^{I_v}\big).
\eEq
Since all of the coefficients in $m_+$ are positive, for any arbitrary $m\!\in\!\K_\R$ there exists a sufficiently large $r\!>\!0$ such that $m\!+\!rm_+\!\in\!\si$. We conclude that $\si$ is top-dimensional.
Since $\R_{\geq 0}^\E \oplus \bigoplus_{v\in \V} \R_{\geq 0}^{I_v}$ is a strictly convex rational polyhedral cone  and $\K_\R$ is an integrally defined sub-vector space, the intersection (\ref{KcapPos_e}) is  a strictly convex rational polyhedral cone.
\ePf

\bCr{m_0integral}
By Lemma~\ref{Convexity_lmm}, the functions $s$ and $\la$ in  Definition~\ref{LogMap_dfn}.\ref{Tropical_l} can be chosen to be 
integral-valued.
\eCr

\noindent
In conclusion, with a set up similar to \cite[Sec~6.3]{FF}, the deformation/obstruction theory of any startum 
$$
\cM_{g,\mfs}(X,D,A)_\Gamma
$$ 
around $f\!=\!(u, [\ze], \Si,z^1,\ldots,z^k)$ is given by (1) $\tn{D}^{\log}_u\dbar$ and $P_I$ for each smooth component $\Si_v$ of $\Si$, and (2) the obstruction map (\ref{PLtoG_e}).

\bLm{ExpectedGamma_lmm}
For any decorated dual graph $\Gamma\!\in\!\tn{DG}(g,\mfs,A)$, the expected complex dimension of $\cM_{g,\mfs}(X,D,A)_\Gamma$ is 
\bEq{ExpectedGamma_e}
c_1^{TX(-\log D)}(A) + (n-3)(1-g) + k - \dim_\R \K_\R(\Gamma)
\eEq
The only stratum with $\dim \K_\R(\Gamma)\!=\!0$ is $\cM_{g,\mfs}(X,D,A)$.
\eLm

\bPf
By (\ref{VirDimRel_e}), the expected complex dimension of each component $\cM_{g_v,\mfs_v}(X,D,A_v)_{I_v}$ is equal to
$$
c_1^{TX(-\log D)}(A_v)+(n-3)(1-g_v)+k_v+\ell_v-|I_v|
$$
where $k_v=|\vec{z}_v|$, $\ell_v=|q_v|$, and $\mfs_v$ is the set of contact order vectors at $\vec{z}_v \cup q_v$.
The pre-log space $\cM^{\tn{plog}}_{g,\mfs}(X,D,A)_\Gamma$ in (\ref{PLtoG_e}) is the fiber product of  $\{\cM_{g_v,\mfs_v}(X,D,A_v)_{I_v}\}_{v\in \V}$
over the evaluation maps at the nodal points
$$
\prod_{v\in \V}\cM_{g_v,\mfs_v}(X,D,A_v)_{I_v} \lra \prod_{e\in \E} (D_{I_e}\times D_{I_e}).
$$
Therefore, using (\ref{genus_e}), the expected complex dimension of $\cM^{\tn{plog}}_{g,\mfs}(X,D,A)_\Gamma$   is 
\bEq{dimBG_e}
\aligned
&\sum_{v\in \V}\big(c_1^{TX(-\log D)}(A_v)+(n-3)(1-g_v)+k_v+\ell_v-|I_v|\big)-\sum_{e\in \E} (n-|I_e|)=\\
&c_1^{TX(-\log D)}(A)+(n-3)(1-g)+k-|\E|-\sum_{v\in \V}|I_v|+\sum_{e\in \E} |I_e|.
\endaligned
\eEq
By (\ref{DtoT_e}), 
$$
\dim_\R \K_\R(\Gamma)-\tn{dim}_\C(\mc{G})=|\E|+\sum_{v\in \V}|I_v|-\sum_{e\in \E}|I_e|.
$$
By (\ref{PLtoG_e}), the stratum $\cM_{g,\mfs}(X,D,A)_\Gamma$ is the pre-image of the identity element under the map
$$
\tn{ob}_{\Gamma}\colon\cM^{\tn{plog}}_{g,\mfs}(X,D,A)_\Gamma\lra \mc{G}.
$$
Therefore, the expected complex dimension of $\cM_{g,\mfs}(X,D,A)_\Gamma$ is equal to the difference of (\ref{dimBG_e}) and 
$$
\tn{dim}_\C(\mc{G})=\dim_\R \K_\R(\Gamma)- \bigg(|\E|+\sum_{v\in \V}|I_v|-\sum_{e\in \E}|I_e|\bigg)
$$
which is equal to (\ref{ExpectedGamma_e}).\\

\noindent
By Definition~\ref{LogMap_dfn}.\ref{Tropical_l} and (\ref{kernel_e}), a function $(s,\la)$ as in Definition~\ref{LogMap_dfn}.\ref{Tropical_l} gives us an element of $\K_\R(\Gamma)$. This element is trivial only if $\Gamma\!=\!\{v\}$ is a one-vertex graph with no edge and $I_v\!=\!\eset$. This establishes the last claim.
\ePf

\subsection{Gluing parameters}\label{Gluing_ss}
The last step in describing the deformation theory and establishing $(\star)$  is to prove a gluing theorem for smoothing the nodes (i.e. deformations normal to each stratum). In this section, we describe the space of gluing parameters for each $\Gamma\!\in\!\tn{DG}(g,\mfs,A)$ and show that it is essentially an affine toric variety. We sketch our idea for the construction of gluing map and defer to a future work \cite{FT3} for the details.\\

\noindent
For a classical nodal $J$-holomorphic map with $|\E|$ nodes, the space of gluing parameters is a neighborhood of the zero in $\C^\E$.  For a log map $f$ as in (\ref{fplogSetUp_e}), the gluing procedure involves a simultaneous smoothing of the nodes, together with pushing $u_v$ out in the direction of $\ze_{v,i}$ for some $v\!\in\!\V$ and $i\!\in\!I_v$.
Thus, a priori, the space of gluing parameters could be quite complicated
and the log moduli spaces (\ref{LogCpt_e}) are not always virtually smooth. For example, the log moduli space of Example~\ref{d1Rd22Pt_ex} below has an $A_1$-singularity along some stratum. For the log moduli spaces, the the space of gluing parameters along $\cM_{g,\mfs}(X,D,A)_\Gamma$ belongs to (a neighborhood of the origin in finitely many copies of) the affine toric variety $Y_{\si(\Gamma)}$ constructed from the toric fan $\si(\Gamma)\!\subset\!\K_\R$. In other words, the kernel of (\ref{DtoT_e}) gives the gluing deformation and, by (\ref{GammaLogStratum_e}), the cokernel of that gives the obstruction space for smoothability of such pre-log maps.\\

\noindent
In the following example, we describe a tuple $(X,D,g,\mfs, A,\Gamma)$ where $\cM_{g,\mfs}(X,D,A)_\Gamma$ is a point and $Y_{\si}$ has an $A_1$-singularity at its center. In this example,  the relative moduli space $\ov\cM^{\tn{rel}}_{g,\mfs}(X,D,A)$ replaces the $A_1$-singularity with a small resolution of that.\vskip.2in

\bEx{d1Rd22Pt_ex}
Suppose $X\!=\!\P^3$, $D\!\cong\!\P^1\!\times\!\P^1$ is a smooth degree $2$ hypersurface, 
$$
g\!=\!1,\qquad A\!=\![2]\!\in\!H_2(\P^3,\Z)\!\cong\!\Z,\qquad \tn{and}\quad \mfs\!=\!(0,0,4).
$$ 
By \cite[Lmm 4.2]{IP1} and (\ref{ExpectedGamma_e}), both $\ov\cM^{\tn{rel}}_{g,\mfs}(X,D,A)$ and $\ov\cM^{\log}_{g,\mfs}(X,D,A)$ are of the expected complex dimension $7$. 
Let $\cM^{\tn{rel}}_{g,\mfs}(X,D,A)_\Gamma$ be the stratum of maps in the expanded degeneration $X[2]$ with connected components:
\bIt
\item a degree $1$ map $u_0\colon \P^1\!\to\!X$ (a line) that intersects $D$ at two distinct points (with multiplicity $1$), 
\item a map $u_3\colon \P^1\!\to\!\P_XD$ in the second layer $\{2\}\!\times\P_XD$ of $X[2]$ which is made of a degree $1$ map $\ov{u}_3\colon \P^1\!\to\!D$ and a meromorphic section $\ze$ of $\ov{u}_3^*\cN_XD\!\cong\!\cO_{\P^1}(2)$ with a zero of order $4$ and $2$ poles of order one, and
\item  two maps $u_1,u_2\colon \P^1\!\to\!\P_XD$ in the first layer $\{1\}\!\times\P_XD$ of $X[2]$ carrying the first and the second marked point, respectively,  which are degree $1$ covers of fibers of $\P_XD$ connecting $u_0$ and $u_4$;
\eIt
see Figure~\ref{2ptsink_fg}-Left.
While the stratum $\cM^{\tn{rel}}_{g,\mfs}(X,D,A)_\Gamma$ is of virtual $\C$-codimension $2$, by (\ref{ExpectedGamma_e}), its image 
$$
\cM_{g,\mfs}(X,D,A)_\Gamma=\pi(\cM^{\tn{rel}}_{g,\mfs}(X,D,A)_\Gamma)
$$ 
in the log moduli space, given by the projection map $\pi$ of Proposition~\ref{ReltoLog_prp} below, is of virtual $\C$-codimension $3$. 
\begin{figure}
\begin{pspicture}(28,-1)(15,2.8)
\psset{unit=.3cm}
\psline[linewidth=.1](52,-2)(76,-2)
\psline[linewidth=.1](52,2)(76,2)
\psline[linewidth=.1](52,6)(76,6)

\pscircle*(63,-2){.2}\rput(63,-3){\small{order $4$ contact}}
\rput(53,8){$X$}
\rput(53,4){$\P_XD$}
\rput(53,0){$\P_XD$}

\psccurve(63,10)(59,8)(60,6)(63,7.5)(66,6)(67,8)\pscircle*(60,6){.2}\pscircle*(66,6){.2}
\psccurve(63,0)(62,0.25)(60,2)(59,0)(63,-2)(67,0)(66,2)(64,0.25)\pscircle*(60,6){.2}\pscircle*(66,6){.2}

\pscircle*(60,2){.2}\pscircle*(66,2){.2}
\pscircle(60,4){2}\pscircle(66,4){2}\pscircle*(60,4){.2}\pscircle*(66,4){.2}

\pscircle*(90,8){.2}\pscircle*(95,4){.2}\pscircle*(90,0){.2}\pscircle*(85,4){.2}
\psline(90,8)(95,4)(90,0)(85,4)(90,8)
\psline(90,0)(90,-2)\psline{->}(90,0)(90,-1)\rput(87,-1){\small{order$=(4)$}}
\psline(85,4)(83,4)
\psline(95,4)(97,4)
\psline[linewidth=.15]{->}(90,8)(92.5,6)\rput(95.7,7){\small{$\uvec{e}_2:~ s_{\uvec{e}_2}=(1)$}}
\psline[linewidth=.15]{->}(90,8)(87.5,6)\rput(84,7){\small{$\uvec{e}_1:~ s_{\uvec{e}_1}=(1)$}}
\psline[linewidth=.15]{->}(95,4)(92.5,2)\rput(97,2){\small{$\uvec{e}_4:~ s_{\uvec{e}_4}=(1)$}}
\psline[linewidth=.15]{->}(85,4)(87.5,2)\rput(83,2){\small{$\uvec{e}_3:~ s_{\uvec{e}_3}=(1)$}}
\rput(90,9){\small{$v_0$}}
\rput(93.5,4){\small{$v_2$}}
\rput(90,1){\small{$v_3$}}
\rput(86.5,4){\small{$v_1$}}

\end{pspicture}
\caption{On left, a nodal $2$-marked $g\!=\!1$ relative map in $X[2]$.  On right, the decorated dual graph of the image log map.}
\label{2ptsink_fg}
\end{figure}
In fact, with the labeling and the choice of orientation on the edges of the associated decorated dual graph $\Gamma$ in  Figure~\ref{2ptsink_fg}-Right, we have
$$
\aligned
&\vr\colon \Z^\E\oplus \bigoplus_{i\in [3]}\Z^{I_{v_i}}\cong\Z^{\{e_1,e_2,e_3,e_4\}} \oplus \Z^{\{v_1,v_2,v_3\}} \lra \bigoplus_{i\in [4]} \Z^{I_{e_i}}\cong\Z^{\{e_1,e_2,e_3,e_4\}}\\
&\vr(1_{e_i})=1_{e_i}\quad \forall~i\!\in\![4],\quad  \vr(1_{v_1})=-1_{e_1}+1_{e_3},  \quad \vr(1_{v_2})=-1_{e_2}+1_{e_4}, \quad\tn{and}~~ \vr(1_{v_3})=-1_{e_3}-1_{e_4}.
\endaligned
$$
Therefore, 
$$
\si\!=\!\tn{ker}(\vr_\R)\!\cap\!\big(\R_{\geq 0}^{\{e_1,e_2,e_3,e_4\}} \oplus \R_{\geq 0}^{\{v_1,v_2,v_3\}}\big)
$$ 
is the cone generated by the set of $4$ vectors
$$
\al_1\!=\!1_{v_3}+1_{e_3}+1_{e_4},\quad \al_2\!=\!1_{v_1}+1_{v_3}+1_{e_1}+1_{e_4}, \quad \al_3\!=\!1_{v_2}+1_{v_3}+1_{e_2}+1_{e_3},\quad \al_4\!=\!1_{v_1}+1_{v_2}+1_{v_3}+1_{e_1}+1_{e_2}.
$$
Since the only relation among $\al_i$ is $\al_1\!+\!\al_4\!=\!\al_2\!+\!\al_3$, the associated toric variety $Y_\si$ is isomorphic to the $3$-dimensional affine sub-variety 
$$(x_1x_4\!-\!x_2x_3=0)\!\subset\! \C^4.$$
\qed
\eEx

\noindent
For every log curve $f\!\in\!\cM_{g,\mfs}(X,D,A)_{\Gamma}$ choose a representative 
\bEq{fplogSetUp_e2}
\big(u_v ,\{\ze_{v,i}\}_{i\in I_v},C_v\equiv (\Si_v,\mfj_v,\vec{z}_v)\big)_{v\in \V}
\eEq
and a set of local coordinates $\{z_{\uvec{e}}\}_{\uvec{e}\in \uvec{\E}}$ around the nodes. 
Since $f$ is $\mc{G}$-unobstructed, by the definition or $\eta$ in (\ref{Totaleta_e}), we can choose $\ze_{v,i}$ and $z_{\uvec{e}}$ such that the leading coefficient vectors $\eta_{\uvec{e}}$ in (\ref{etae_e}) satisfy
\bEq{smoothableChoice_e}
\eta_{\uvec{e}}=\eta_{\scz{\ucev{e}}} \qquad \forall~e\!\in\!\E.
\eEq
For every $v\!\in\!\V$ and $i\!\in\![N]\!-\!I_v$, let $t_{v,i}\!=\!1$ in (\ref{GluignEquation_e}).
Then the \textbf{space of gluing parameters} for $f$ is a sufficiently small neighborhood of the origin in the complex sub-variety
\bEq{GluignEquation_e}
\aligned
\cN_\Gamma=\bigg\{\big((\ve_e)_{e\in \E},&(t_{v,i})_{v\in \V, i\in I_v}\big)\!\in\!\C^\E\times \prod_{v\in \V}\C^{I_v}\colon 
~~\ve_e^{s_{\uvec{e},i}}t_{v,i}=t_{v',i}\\
&\forall~v,v'\!\in\!\V,~e\!\in\!\E_{v,v'},i\!\in\!I_e,~ \uvec{e}~~\tn{s.t}~~s_{\uvec{e},i}\geq 0\bigg\}\subset \C^\E\times \prod_{v\in \V}\C^{I_v}.
\endaligned
\eEq 
The complex numbers $\ve_e$ are the gluing parameters for the nodes of $\Si$ and $t_{v,i}$ are the parameters for pushing $u_v$ out in the direction of $\ze_{v,i}$. In the gluing construction outlined below, given a set of representatives $\big(\{z_{\uvec{e}}\}_{\uvec{e}\in\uvec{\E}},\{\ze_{v,i}\}_{v\in \V, i\in I_v}\big)$ satisfying (\ref{smoothableChoice_e}) and a sufficiently small
$$
(\ve,t)\!\equiv\! \big((\ve_e)_{e\in \E},(t_{v,i})_{v\in \V, i\in I_v}\big)\!\in\!\cN_\Gamma,
$$ 
we will construct a pre-gluing log map $\wt{f}_{\ve,t}$. Then we must show that there is an actual log $J$-holomorphic map ``close" to it.\\

\noindent
Let 
$$
\T^\vee\cong\bigoplus_{e\in \E} \Z^{I_e}\stackrel{\vr^\vee}{\xrightarrow{\hspace*{1.5cm}}} 
\D^\vee\cong \Z^\E\oplus \bigoplus_{v\in \V} \Z^{I_v}  
$$
be the dual of $\Z$-linear map $\vr$ associated to $\Gamma$ in (\ref{DtoT_e}) (for a fixed choice of orientation $O$ on $\E$). With the kernel subspace $\K\!=\!\tn{ker}(\vr)\!\subset\!\D$ as in (\ref{cLGamma_e}), let 
$$
\K^\perp=\{m\in \D^\vee\colon \left\langle m,\al\right\rangle =0 \quad \forall~\al\!\in\!\K\}\subset \D^\vee.
$$ 
Then $\tn{Im}(\vr^\vee)\!\subset\!\K^\perp$ with the finite quotient 
$$
\K^{\perp}/\tn{image}(\vr^\vee).
$$

\bPr{GluignSpace_prp}
The space of gluing parameters $\cN_\Gamma$ in (\ref{GluignEquation_e}) is a possibly non-irreducible and non-reduced affine toric sub-variety of $\C^\E\!\times\!\prod_{v\in \V} \C^{I_v}$  that is isomorphic to $\abs{\K^{\perp}/\tn{Im}(\vr^\vee)}$ copies of the irreducible reduced affine toric variety $Y_{\si(\Gamma)}$ (with toric fan $\si$), counting with multiplicities\footnote{we do not know of any example, arising from such dual graphs, such that the multiplicities are bigger than $1$.}.
Replacing $\{z_{\uvec{e}}\}_{\uvec{e}\in \uvec{\E}}$ and $\{\ze_{v,i}\}_{v\in \V, i\in I_v}$ with another choice satisfying (\ref{smoothableChoice_e}) corresponds to a torus action on $\cN_\Gamma$.
\ePr

\bPf
Let us start with some general facts about toric varieties.
For $n\!\in\!\Z_+$, every vector $m\!\in\!\Z^n$ has a unique presentation $m\!=\!m_+\!-\!m_-$ such that $m_+,m_-\!\in\!(\Z_{\geq 0})^n$. Every $m\!=\!(a_1,\ldots,a_n)\!\in\!(\Z_{\geq 0})^n$ corresponds to the monomial 
$$
x^m\equiv x_1^{a_1}\cdots x_n^{a_n}\!\in\!\C[x_1,\ldots,x_n].
$$
For every arbitrary $m\!\in\!\Z^n$, the binomial corresponding to $m$ is the expression
$$
x^{m,\pm}\equiv x^{m_+}-x^{m_-}\!\in\!\C[x_1,\ldots,x_n].
$$
For example, if $m\!=\!0$, then $x^{m,\pm}\!=\!1-1\!=\!0$.
A binomial ideal\footnote{for more general binomial ideals see~\cite{ES}.} $I$ in $\C[x_1,\ldots,x_n]$ is an ideal generated by a finite set of binomials $x^{m_1,\pm},\ldots, x^{m_\ell,\pm}$. \\

\noindent
Suppose $\K^\vee\cong \Z^{\ell}$ is a lattice and $\Z^n\!\lra\!\K^\vee$
is a surjective $\Z$-linear map. Let $\R^n\!\lra\!\K_\R^\vee$ be the corresponding $\R$-linear projection map and $\si^\vee$ be the image of  the cone $\R_{\geq 0}^n$ in $\K^\vee_\R$. Then the dual map $\iota\colon \K\hookrightarrow \Z^n$ 
is an embedding and the dual of $\si^\vee$ is the toric fan
$$
\si\!=\!\K_\R \cap \iota^{-1}( \R_{\geq 0}^n).
$$
In this situation, by \cite[Prp 1.1.9]{CLS}, the toric variety $Y_{\si}$ associated to the toric fan $\si$ is the zero set of the binomial ideal
\bEq{Ideal_e}
I=\{ x^{m,\pm}\colon~~m\!\in\!\K^\perp \subset \Z^n\}.
\eEq
\vskip.1in

\noindent
With $\Z^n\!=\! \Z^\E\oplus \bigoplus_{v\in \V} \Z^{I_v}$, $\K$ in (\ref{kernel_e}), and $\si\!=\!\si(\Gamma)$ as in (\ref{SiCone_e}), the previous argument implies that $Y_{\si(\Gamma)}$ is the zero 
set of the binomial ideal~(\ref{Ideal_e}).\\

\noindent
Let $I'\!\subset\!I$ be the binomial sub-ideal generated by the elements of $\tn{Im}(\vr^\vee)\!\subset\!\K^\perp$. By definition of $\vr$ and (\ref{GluignEquation_e}), the space of gluing parameters $\cN_\Gamma$ is the zero set (scheme) of $I'$. Therefore $Y_{\si(\Gamma)}\!\subset\!\cN_\Gamma$. Note that $Y_{\si(\Gamma)}$ is the Zariski closure of the irreducible subgroup
$$
\big\{t\in (\C^*)^n\colon t^m=1\qquad \forall~m\!\in\!\K^\perp\big\}\subset (\C^*)^n
$$ 
and $\cN_\Gamma$ is the Zariski closure of possibly non-irreducible subgroup
\bEq{TorusofGl_e}
 \big\{t\in (\C^*)^n\colon t^m=1\qquad \forall~m\!\in\!\tn{Im}(\vr^\vee)\big\}\subset (\C^*)^n;
\eEq
see \cite[Dfn 1.1.7]{CLS}. Therefore, all the irreducible components of $\cN_\Gamma$ are isomorphic to $Y_{\si(\Gamma)}$. Since 
\bEq{my_e}
\abs{I/I'}=m_\Gamma:=\abs{\K^{\perp}/\tn{Im}(\vr^\vee)},
\eEq
$\cN_\Gamma$ is isomorphic to $m_\Gamma$ copies of $Y_{\si(\Gamma)}$, counting with multiplicities.
The last statement in Proposition~\ref{GluignSpace_prp} follows from the way subgroup (\ref{TorusofGl_e}) acts on (\ref{GluignEquation_e}).
\ePf

\bEx{ToricEx_e}
Suppose $N\!=\!2$ and $\Gamma$ is the decorated dual graph with two vertices $\V\!=\!\{v_1,v_2\}$ and two edges $e_1$ and $e_2$ connecting them. Choose $\uvec{e}_1$ and $\uvec{e}_2$ to be the orientations starting at $v_1$. Suppose 
$$
I_{v_1}=\{1\}, \quad I_{v_2}=\{2\}, \quad s_{\uvec{e}_1}=s_{\uvec{e}_2}=(-2,2).
$$
Then the linear map 
$$
\vr\colon \Z^\E\oplus \Z^{I_{v_1}}\oplus \Z^{I_{v_2}} \equiv \Z_{e_1}\oplus \Z_{e_2}\oplus \Z_{v_1}\oplus \Z_{v_2} \lra \Z^{\{1,2\}}_{e_1} \oplus \Z^{\{1,2\}}_{e_2}
$$
is given by 
$$
\aligned
&\vr(1_{e_1}) = ((-2,2)_{e_1}, (0,0)_{e_2}), \quad &\vr(1_{e_2})=((0,0)_{e_1}, (-2,2)_{e_2}),\\
&\vr(1_{v_1})= ((1,0)_{e_1}, (1,0)_{e_2}),\quad &\vr(1_{v_2})= ((0,-1)_{e_1}, (0,-1)_{e_2}).
\endaligned
$$
It is straightforward to check that $\tn{Ker}(\vr)$ is one-dimensional and is generated by 
$$
1_{e_1}\!+\!1_{e_2}\!+\!2\cdot1_{v_1}\!+\!2\cdot 1_{v_2}\,;
$$ 
i.e. $Y_{\si(\Gamma)}\cong \C$. On the other hand, 
$\cN_\Gamma$ is the sub-variety cut-out by
$$
\ve_1^2=t_{v_2}, \quad \ve_2^2=t_{v_2},\quad \ve_1^2=t_{v_1}, \quad \ve_2^2=t_{v_1}.
$$
This is isomorphic to $2$ copies of $\C$, the component $Y_{\si(\Gamma)}$ is the image of $t\lra (t,t,t^2,t^2)$ and the other one is the image of $t\lra (t,-t,t^2,t^2)$. It is straightforward to see that 
$$
\tn{Ker}(\vr)^\perp/\tn{Im}(\vr^{\vee})
$$
is isomorphic to $\Z_2$ and is generated by the class of $[1_{e_1}^\vee-1_{e_2}^\vee]$.
\qed
\eEx

\noindent
Given a log $J$-holomorphic map $f\!=\!(u, [\ze], \Si,z^1,\ldots,z^k)$ in $\ov\cM^{\log}_{g,\mfs}(X,D,A)$ with nodal domain (\ref{nodalcurve_e}), a set of local coordinates $\{z_{\uvec{e}}\}_{\uvec{e}\in \uvec{\E}}$ around the nodes such that (\ref{smoothableChoice_e}) holds, and a gluing parameter 
$(\ve,t)\!\equiv\!\big((\ve_e)_{e\in \E},(t_{v,i})_{v\in \V, i\in I_v}\big)$ in (\ref{GluignEquation_e}), the gluing construction can/will be done in the following way.\\

\noindent 
Consider for example a node $q_e$ connecting $\Si_v$ and $\Si_{v'}$ with $\tn{ord}_{q_{\uvec{e}}}(u,D_i)\!=\!s_{\uvec{e},i}\!>\!0$. Then the log tuple on $\Si_{v'}$ includes a section $\ze_{v',i}$ of $u_{v'}^*\cN_XD_i$ with a pole of order $s_{\uvec{e},i}$ at the nodal point $q_{\scz\ucev{e}}\!\in\!\Si_{v'}$. Near $q_{\uvec{e}}$, $u_{v}$ has the product form 
$$
u_{v}(z_{\uvec{e}})=(\eta_{\uvec{e},i}z_{\uvec{e}}^{s_{\uvec{e},i}}, \ov{u}_v) \in \C\times D_i.
$$
On the other hand, $\ze_{v',i}$ has a local expansion $\ze_{v',i}(z_{\scz\ucev{e}})=\eta_{\scz\ucev{e},i}z_{\scz\ucev{e}}^{-s_{\uvec{e},i}}+ \cdots$. By (\ref{smoothableChoice_e}) and (\ref{GluignEquation_e}), we have
\bEq{EE_e}
\ve_e^{s_{\uvec{e},i}} t_{v,i}\eta_{\uvec{e},i}= t_{v',i}\eta_{\scz\ucev{e},i}
\eEq
\textit{at all the nodes, simultaneously}. The smoothing of $\Si$ is given by smoothing the nodes $q_{e}$ via the equation $z_{\uvec{e}}z_{\scz\ucev{e}}=\ve_e$. The identity (\ref{EE_e}) means that the expression
\bEq{NE_e}
\eta_{\uvec{e},i}t_{v,i}z_{\uvec{e}}^{s_{\uvec{e},i}}= \eta_{\scz\ucev{e},i}t_{v',i}z_{\scz\ucev{e}}^{{s_{\scz{\ucev{e}},i}}}
\eEq
defines a function from the neck region into $\cN_{X}D_i$. We then construct the approximate-gluing log map $\wt{f}_{\ve,t}$ in the following way. On each neck region (\textit{unlike in the classical gluing construction were the approximate-gluing map is defined to be constant}), we define the approximate-gluing map to be (\ref{NE_e}) in the $i$-th direction. Away from the nodes, $\wt{f}_{\ve,t}$ is defined to be the push out\footnote{via the regularization maps in $\cR$.} of $u_v$ via the section $\sum_{i\in I_v} t_{v,i}\ze_{v,i}$ on the $v$-th component. The latter is $J$-holomorphic due to some properties of $\tn{AK}(X,D)$. In between the two regions, $\wt{f}_{\ve,t}$ interpolates between the two maps. Then, with $\tn{D}^{\log}\dbar$ in place of $\tn{D}\dbar$ in \cite[Ch~10]{MS2}, an argument similar to the classical argument allows us to find a log $J$-holomorphic map close to $\wt{f}$.

\vskip.1in
\noindent
{\it MacLean Hall, The University of Iowa, Iowa City, IA 52242\\
mohammad-tehrani@uiowa.edu}
\small
%------------------------------------------------------------------------------------------------------

\end{document}